\documentclass{article}
\usepackage[utf8]{inputenc}

\usepackage{eurosym}
\usepackage{amsfonts}
\usepackage{amsmath}
\usepackage[left=2cm, right=2cm, top=2cm, bottom=2cm]{geometry}
\usepackage{amssymb}
\usepackage{mathrsfs} 
\usepackage{physics}
\usepackage{siunitx}
\usepackage[dvipsnames]{xcolor}
\usepackage{tikz}

\usepackage{caption}

\usepackage{titling}
\usepackage{graphicx}
\usepackage{duckuments}
\usepackage{caption}
\usepackage[onehalfspacing]{setspace}
\usepackage{hyperref}
    \hypersetup{
        colorlinks=true,
        citecolor=Green,
        linkcolor=black,
        filecolor=magenta,      
        urlcolor=cyan,
        }

\usepackage{subcaption}
\usepackage{booktabs}
\newcommand{\ra}[1]{\renewcommand{\arraystretch}{#1}}
\usepackage[mathlines]{lineno}
\usepackage[color=yellow, colorinlistoftodos]{todonotes}
\usepackage{framed}
\usepackage[titletoc, title]{appendix}
    
\usepackage{scalerel}
\usepackage{pict2e}
\usepackage{tkz-euclide}
\usepackage{anyfontsize}
\usetikzlibrary{calc}
\usetikzlibrary{patterns, arrows.meta, decorations.pathreplacing}
\usetikzlibrary{shadows}
\usetikzlibrary{external}

\usepackage{pgfplots}
\pgfplotsset{compat=newest}
\usepgfplotslibrary{statistics}
\usepgfplotslibrary{fillbetween}

\definecolor{WaterBlue}{rgb}{0.11, 0.01 , 0.88}
\definecolor{AliceBlue}{rgb}{0.94, 0.972, 1.0}
\definecolor{AntiqueWhite}{rgb}{0.98, 0.92, 0.844}
\definecolor{AntiqueWhite1}{rgb}{1.0, 0.936, 0.86}
\definecolor{AntiqueWhite2}{rgb}{0.932, 0.875, 0.8}
\definecolor{AntiqueWhite3}{rgb}{0.804, 0.752, 0.69}
\definecolor{AntiqueWhite4}{rgb}{0.545, 0.512, 0.47}
\definecolor{Aqua}{rgb}{0.0, 1.0, 1.0}
\definecolor{Aquamarine}{rgb}{0.498, 1.0, 0.83}
\definecolor{Aquamarine1}{rgb}{0.498, 1.0, 0.83}
\definecolor{Aquamarine2}{rgb}{0.464, 0.932, 0.776}
\definecolor{Aquamarine3}{rgb}{0.4, 0.804, 0.668}
\definecolor{Aquamarine4}{rgb}{0.27, 0.545, 0.455}
\definecolor{Azure}{rgb}{0.94, 1.0, 1.0}
\definecolor{Azure1}{rgb}{0.94, 1.0, 1.0}
\definecolor{Azure2}{rgb}{0.88, 0.932, 0.932}
\definecolor{Azure3}{rgb}{0.756, 0.804, 0.804}
\definecolor{Azure4}{rgb}{0.512, 0.545, 0.545}
\definecolor{Beige}{rgb}{0.96, 0.96, 0.864}
\definecolor{Bisque}{rgb}{1.0, 0.894, 0.77}
\definecolor{Bisque1}{rgb}{1.0, 0.894, 0.77}
\definecolor{Bisque2}{rgb}{0.932, 0.835, 0.716}
\definecolor{Bisque3}{rgb}{0.804, 0.716, 0.62}
\definecolor{Bisque4}{rgb}{0.545, 0.49, 0.42}
\definecolor{Black}{rgb}{0.0, 0.0, 0.0}
\definecolor{BlanchedAlmond}{rgb}{1.0, 0.92, 0.804}
\definecolor{Blue}{rgb}{0.0, 0.0, 1.0}
\definecolor{Blue1}{rgb}{0.0, 0.0, 1.0}
\definecolor{Blue2}{rgb}{0.0, 0.0, 0.932}
\definecolor{Blue3}{rgb}{0.0, 0.0, 0.804}
\definecolor{Blue4}{rgb}{0.0, 0.0, 0.545}
\definecolor{BlueViolet}{rgb}{0.54, 0.17, 0.888}
\definecolor{Brown}{rgb}{0.648, 0.165, 0.165}
\definecolor{Brown1}{rgb}{1.0, 0.25, 0.25}
\definecolor{Brown2}{rgb}{0.932, 0.23, 0.23}
\definecolor{Brown3}{rgb}{0.804, 0.2, 0.2}
\definecolor{Brown4}{rgb}{0.545, 0.136, 0.136}
\definecolor{BurlyWood}{rgb}{0.87, 0.72, 0.53}
\definecolor{Burlywood1}{rgb}{1.0, 0.828, 0.608}
\definecolor{Burlywood2}{rgb}{0.932, 0.772, 0.57}
\definecolor{Burlywood3}{rgb}{0.804, 0.668, 0.49}
\definecolor{Burlywood4}{rgb}{0.545, 0.45, 0.332}
\definecolor{CadetBlue}{rgb}{0.372, 0.62, 0.628}
\definecolor{CadetBlue1}{rgb}{0.596, 0.96, 1.0}
\definecolor{CadetBlue2}{rgb}{0.556, 0.898, 0.932}
\definecolor{CadetBlue3}{rgb}{0.48, 0.772, 0.804}
\definecolor{CadetBlue4}{rgb}{0.325, 0.525, 0.545}
\definecolor{Chartreuse}{rgb}{0.498, 1.0, 0.0}
\definecolor{Chartreuse1}{rgb}{0.498, 1.0, 0.0}
\definecolor{Chartreuse2}{rgb}{0.464, 0.932, 0.0}
\definecolor{Chartreuse3}{rgb}{0.4, 0.804, 0.0}
\definecolor{Chartreuse4}{rgb}{0.27, 0.545, 0.0}
\definecolor{Chocolate}{rgb}{0.824, 0.41, 0.116}
\definecolor{Chocolate1}{rgb}{1.0, 0.498, 0.14}
\definecolor{Chocolate2}{rgb}{0.932, 0.464, 0.13}
\definecolor{Chocolate3}{rgb}{0.804, 0.4, 0.112}
\definecolor{Chocolate4}{rgb}{0.545, 0.27, 0.075}
\definecolor{Coral}{rgb}{1.0, 0.498, 0.312}
\definecolor{Coral1}{rgb}{1.0, 0.448, 0.336}
\definecolor{Coral2}{rgb}{0.932, 0.415, 0.312}
\definecolor{Coral3}{rgb}{0.804, 0.356, 0.27}
\definecolor{Coral4}{rgb}{0.545, 0.244, 0.185}
\definecolor{CornflowerBlue}{rgb}{0.392, 0.585, 0.93}
\definecolor{Cornsilk}{rgb}{1.0, 0.972, 0.864}
\definecolor{Cornsilk1}{rgb}{1.0, 0.972, 0.864}
\definecolor{Cornsilk2}{rgb}{0.932, 0.91, 0.804}
\definecolor{Cornsilk3}{rgb}{0.804, 0.785, 0.694}
\definecolor{Cornsilk4}{rgb}{0.545, 0.532, 0.47}
\definecolor{Crimson}{rgb}{0.864, 0.08, 0.235}
\definecolor{Cyan}{rgb}{0.0, 1.0, 1.0}
\definecolor{Cyan1}{rgb}{0.0, 1.0, 1.0}
\definecolor{Cyan2}{rgb}{0.0, 0.932, 0.932}
\definecolor{Cyan3}{rgb}{0.0, 0.804, 0.804}
\definecolor{Cyan4}{rgb}{0.0, 0.545, 0.545}
\definecolor{DarkBlue}{rgb}{0.0, 0.0, 0.545}
\definecolor{DarkCyan}{rgb}{0.0, 0.545, 0.545}
\definecolor{DarkGoldenrod}{rgb}{0.72, 0.525, 0.044}
\definecolor{DarkGoldenrod1}{rgb}{1.0, 0.725, 0.06}
\definecolor{DarkGoldenrod2}{rgb}{0.932, 0.68, 0.055}
\definecolor{DarkGoldenrod3}{rgb}{0.804, 0.585, 0.048}
\definecolor{DarkGoldenrod4}{rgb}{0.545, 0.396, 0.03}
\definecolor{DarkGray}{rgb}{0.664, 0.664, 0.664}
\definecolor{DarkGreen}{rgb}{0.0, 0.392, 0.0}
\definecolor{DarkGrey}{rgb}{0.664, 0.664, 0.664}
\definecolor{DarkKhaki}{rgb}{0.74, 0.716, 0.42}
\definecolor{DarkMagenta}{rgb}{0.545, 0.0, 0.545}
\definecolor{DarkOliveGreen}{rgb}{0.332, 0.42, 0.185}
\definecolor{DarkOliveGreen1}{rgb}{0.792, 1.0, 0.44}
\definecolor{DarkOliveGreen2}{rgb}{0.736, 0.932, 0.408}
\definecolor{DarkOliveGreen3}{rgb}{0.635, 0.804, 0.352}
\definecolor{DarkOliveGreen4}{rgb}{0.43, 0.545, 0.24}
\definecolor{DarkOrange}{rgb}{1.0, 0.55, 0.0}
\definecolor{DarkOrange1}{rgb}{1.0, 0.498, 0.0}
\definecolor{DarkOrange2}{rgb}{0.932, 0.464, 0.0}
\definecolor{DarkOrange3}{rgb}{0.804, 0.4, 0.0}
\definecolor{DarkOrange4}{rgb}{0.545, 0.27, 0.0}
\definecolor{DarkOrchid}{rgb}{0.6, 0.196, 0.8}
\definecolor{DarkOrchid1}{rgb}{0.75, 0.244, 1.0}
\definecolor{DarkOrchid2}{rgb}{0.698, 0.228, 0.932}
\definecolor{DarkOrchid3}{rgb}{0.604, 0.196, 0.804}
\definecolor{DarkOrchid4}{rgb}{0.408, 0.132, 0.545}
\definecolor{DarkRed}{rgb}{0.545, 0.0, 0.0}
\definecolor{DarkSalmon}{rgb}{0.912, 0.59, 0.48}
\definecolor{DarkSeaGreen}{rgb}{0.56, 0.736, 0.56}
\definecolor{DarkSeaGreen1}{rgb}{0.756, 1.0, 0.756}
\definecolor{DarkSeaGreen2}{rgb}{0.705, 0.932, 0.705}
\definecolor{DarkSeaGreen3}{rgb}{0.608, 0.804, 0.608}
\definecolor{DarkSeaGreen4}{rgb}{0.41, 0.545, 0.41}
\definecolor{DarkSlateBlue}{rgb}{0.284, 0.24, 0.545}
\definecolor{DarkSlateGray}{rgb}{0.185, 0.31, 0.31}
\definecolor{DarkSlateGray1}{rgb}{0.592, 1.0, 1.0}
\definecolor{DarkSlateGray2}{rgb}{0.552, 0.932, 0.932}
\definecolor{DarkSlateGray3}{rgb}{0.475, 0.804, 0.804}
\definecolor{DarkSlateGray4}{rgb}{0.32, 0.545, 0.545}
\definecolor{DarkSlateGrey}{rgb}{0.185, 0.31, 0.31}
\definecolor{DarkTurquoise}{rgb}{0.0, 0.808, 0.82}
\definecolor{DarkViolet}{rgb}{0.58, 0.0, 0.828}
\definecolor{DeepPink}{rgb}{1.0, 0.08, 0.576}
\definecolor{DeepPink1}{rgb}{1.0, 0.08, 0.576}
\definecolor{DeepPink2}{rgb}{0.932, 0.07, 0.536}
\definecolor{DeepPink3}{rgb}{0.804, 0.064, 0.464}
\definecolor{DeepPink4}{rgb}{0.545, 0.04, 0.312}
\definecolor{DeepSkyBlue}{rgb}{0.0, 0.75, 1.0}
\definecolor{DeepSkyBlue1}{rgb}{0.0, 0.75, 1.0}
\definecolor{DeepSkyBlue2}{rgb}{0.0, 0.698, 0.932}
\definecolor{DeepSkyBlue3}{rgb}{0.0, 0.604, 0.804}
\definecolor{DeepSkyBlue4}{rgb}{0.0, 0.408, 0.545}
\definecolor{DimGray}{rgb}{0.41, 0.41, 0.41}
\definecolor{DimGrey}{rgb}{0.41, 0.41, 0.41}
\definecolor{DodgerBlue}{rgb}{0.116, 0.565, 1.0}
\definecolor{DodgerBlue1}{rgb}{0.116, 0.565, 1.0}
\definecolor{DodgerBlue2}{rgb}{0.11, 0.525, 0.932}
\definecolor{DodgerBlue3}{rgb}{0.094, 0.455, 0.804}
\definecolor{DodgerBlue4}{rgb}{0.064, 0.305, 0.545}
\definecolor{FireBrick}{rgb}{0.698, 0.132, 0.132}
\definecolor{Firebrick1}{rgb}{1.0, 0.19, 0.19}
\definecolor{Firebrick2}{rgb}{0.932, 0.172, 0.172}
\definecolor{Firebrick3}{rgb}{0.804, 0.15, 0.15}
\definecolor{Firebrick4}{rgb}{0.545, 0.1, 0.1}
\definecolor{FloralWhite}{rgb}{1.0, 0.98, 0.94}
\definecolor{ForestGreen}{rgb}{0.132, 0.545, 0.132}
\definecolor{Fuchsia}{rgb}{1.0, 0.0, 1.0}
\definecolor{Gainsboro}{rgb}{0.864, 0.864, 0.864}
\definecolor{GhostWhite}{rgb}{0.972, 0.972, 1.0}
\definecolor{Gold}{rgb}{1.0, 0.844, 0.0}
\definecolor{Gold1}{rgb}{1.0, 0.844, 0.0}
\definecolor{Gold2}{rgb}{0.932, 0.79, 0.0}
\definecolor{Gold3}{rgb}{0.804, 0.68, 0.0}
\definecolor{Gold4}{rgb}{0.545, 0.46, 0.0}
\definecolor{Goldenrod}{rgb}{0.855, 0.648, 0.125}
\definecolor{Goldenrod1}{rgb}{1.0, 0.756, 0.145}
\definecolor{Goldenrod2}{rgb}{0.932, 0.705, 0.132}
\definecolor{Goldenrod3}{rgb}{0.804, 0.608, 0.112}
\definecolor{Goldenrod4}{rgb}{0.545, 0.41, 0.08}
\definecolor{Gray}{rgb}{0.5, 0.5, 0.5}
\definecolor{Gray0}{rgb}{0.745, 0.745, 0.745}
\definecolor{Green}{rgb}{0.0, 0.5, 0.0}
\definecolor{Green0}{rgb}{0.0, 1.0, 0.0}
\definecolor{Green1}{rgb}{0.0, 1.0, 0.0}
\definecolor{Green2}{rgb}{0.0, 0.932, 0.0}
\definecolor{Green3}{rgb}{0.0, 0.804, 0.0}
\definecolor{Green4}{rgb}{0.0, 0.545, 0.0}
\definecolor{GreenYellow}{rgb}{0.68, 1.0, 0.185}
\definecolor{Grey}{rgb}{0.5, 0.5, 0.5}
\definecolor{Grey0}{rgb}{0.745, 0.745, 0.745}
\definecolor{Honeydew}{rgb}{0.94, 1.0, 0.94}
\definecolor{Honeydew1}{rgb}{0.94, 1.0, 0.94}
\definecolor{Honeydew2}{rgb}{0.88, 0.932, 0.88}
\definecolor{Honeydew3}{rgb}{0.756, 0.804, 0.756}
\definecolor{Honeydew4}{rgb}{0.512, 0.545, 0.512}
\definecolor{HotPink}{rgb}{1.0, 0.41, 0.705}
\definecolor{HotPink1}{rgb}{1.0, 0.43, 0.705}
\definecolor{HotPink2}{rgb}{0.932, 0.415, 0.655}
\definecolor{HotPink3}{rgb}{0.804, 0.376, 0.565}
\definecolor{HotPink4}{rgb}{0.545, 0.228, 0.385}
\definecolor{IndianRed}{rgb}{0.804, 0.36, 0.36}
\definecolor{IndianRed1}{rgb}{1.0, 0.415, 0.415}
\definecolor{IndianRed2}{rgb}{0.932, 0.39, 0.39}
\definecolor{IndianRed3}{rgb}{0.804, 0.332, 0.332}
\definecolor{IndianRed4}{rgb}{0.545, 0.228, 0.228}
\definecolor{Indigo}{rgb}{0.294, 0.0, 0.51}
\definecolor{Ivory}{rgb}{1.0, 1.0, 0.94}
\definecolor{Ivory1}{rgb}{1.0, 1.0, 0.94}
\definecolor{Ivory2}{rgb}{0.932, 0.932, 0.88}
\definecolor{Ivory3}{rgb}{0.804, 0.804, 0.756}
\definecolor{Ivory4}{rgb}{0.545, 0.545, 0.512}
\definecolor{Khaki}{rgb}{0.94, 0.9, 0.55}
\definecolor{Khaki1}{rgb}{1.0, 0.965, 0.56}
\definecolor{Khaki2}{rgb}{0.932, 0.9, 0.52}
\definecolor{Khaki3}{rgb}{0.804, 0.776, 0.45}
\definecolor{Khaki4}{rgb}{0.545, 0.525, 0.305}
\definecolor{Lavender}{rgb}{0.9, 0.9, 0.98}
\definecolor{LavenderBlush}{rgb}{1.0, 0.94, 0.96}
\definecolor{LavenderBlush1}{rgb}{1.0, 0.94, 0.96}
\definecolor{LavenderBlush2}{rgb}{0.932, 0.88, 0.898}
\definecolor{LavenderBlush3}{rgb}{0.804, 0.756, 0.772}
\definecolor{LavenderBlush4}{rgb}{0.545, 0.512, 0.525}
\definecolor{LawnGreen}{rgb}{0.488, 0.99, 0.0}
\definecolor{LemonChiffon}{rgb}{1.0, 0.98, 0.804}
\definecolor{LemonChiffon1}{rgb}{1.0, 0.98, 0.804}
\definecolor{LemonChiffon2}{rgb}{0.932, 0.912, 0.75}
\definecolor{LemonChiffon3}{rgb}{0.804, 0.79, 0.648}
\definecolor{LemonChiffon4}{rgb}{0.545, 0.536, 0.44}
\definecolor{LightBlue}{rgb}{0.68, 0.848, 0.9}
\definecolor{LightBlue1}{rgb}{0.75, 0.936, 1.0}
\definecolor{LightBlue2}{rgb}{0.698, 0.875, 0.932}
\definecolor{LightBlue3}{rgb}{0.604, 0.752, 0.804}
\definecolor{LightBlue4}{rgb}{0.408, 0.512, 0.545}
\definecolor{LightCoral}{rgb}{0.94, 0.5, 0.5}
\definecolor{LightCyan}{rgb}{0.88, 1.0, 1.0}
\definecolor{LightCyan1}{rgb}{0.88, 1.0, 1.0}
\definecolor{LightCyan2}{rgb}{0.82, 0.932, 0.932}
\definecolor{LightCyan3}{rgb}{0.705, 0.804, 0.804}
\definecolor{LightCyan4}{rgb}{0.48, 0.545, 0.545}
\definecolor{LightGoldenrod}{rgb}{0.933, 0.867, 0.51}
\definecolor{LightGoldenrod1}{rgb}{1.0, 0.925, 0.545}
\definecolor{LightGoldenrod2}{rgb}{0.932, 0.864, 0.51}
\definecolor{LightGoldenrod3}{rgb}{0.804, 0.745, 0.44}
\definecolor{LightGoldenrod4}{rgb}{0.545, 0.505, 0.298}
\definecolor{LightGoldenrodYellow}{rgb}{0.98, 0.98, 0.824}
\definecolor{LightGray}{rgb}{0.828, 0.828, 0.828}
\definecolor{LightGreen}{rgb}{0.565, 0.932, 0.565}
\definecolor{LightGrey}{rgb}{0.828, 0.828, 0.828}
\definecolor{LightPink}{rgb}{1.0, 0.712, 0.756}
\definecolor{LightPink1}{rgb}{1.0, 0.684, 0.725}
\definecolor{LightPink2}{rgb}{0.932, 0.635, 0.68}
\definecolor{LightPink3}{rgb}{0.804, 0.55, 0.585}
\definecolor{LightPink4}{rgb}{0.545, 0.372, 0.396}
\definecolor{LightSalmon}{rgb}{1.0, 0.628, 0.48}
\definecolor{LightSalmon1}{rgb}{1.0, 0.628, 0.48}
\definecolor{LightSalmon2}{rgb}{0.932, 0.585, 0.448}
\definecolor{LightSalmon3}{rgb}{0.804, 0.505, 0.385}
\definecolor{LightSalmon4}{rgb}{0.545, 0.34, 0.26}
\definecolor{LightSeaGreen}{rgb}{0.125, 0.698, 0.668}
\definecolor{LightSkyBlue}{rgb}{0.53, 0.808, 0.98}
\definecolor{LightSkyBlue1}{rgb}{0.69, 0.888, 1.0}
\definecolor{LightSkyBlue2}{rgb}{0.644, 0.828, 0.932}
\definecolor{LightSkyBlue3}{rgb}{0.552, 0.712, 0.804}
\definecolor{LightSkyBlue4}{rgb}{0.376, 0.484, 0.545}
\definecolor{LightSlateBlue}{rgb}{0.518, 0.44, 1.0}
\definecolor{LightSlateGray}{rgb}{0.468, 0.532, 0.6}
\definecolor{LightSlateGrey}{rgb}{0.468, 0.532, 0.6}
\definecolor{LightSteelBlue}{rgb}{0.69, 0.77, 0.87}
\definecolor{LightSteelBlue1}{rgb}{0.792, 0.884, 1.0}
\definecolor{LightSteelBlue2}{rgb}{0.736, 0.824, 0.932}
\definecolor{LightSteelBlue3}{rgb}{0.635, 0.71, 0.804}
\definecolor{LightSteelBlue4}{rgb}{0.43, 0.484, 0.545}
\definecolor{LightYellow}{rgb}{1.0, 1.0, 0.88}
\definecolor{LightYellow1}{rgb}{1.0, 1.0, 0.88}
\definecolor{LightYellow2}{rgb}{0.932, 0.932, 0.82}
\definecolor{LightYellow3}{rgb}{0.804, 0.804, 0.705}
\definecolor{LightYellow4}{rgb}{0.545, 0.545, 0.48}
\definecolor{Lime}{rgb}{0.0, 1.0, 0.0}
\definecolor{LimeGreen}{rgb}{0.196, 0.804, 0.196}
\definecolor{Linen}{rgb}{0.98, 0.94, 0.9}
\definecolor{Magenta}{rgb}{1.0, 0.0, 1.0}
\definecolor{Magenta1}{rgb}{1.0, 0.0, 1.0}
\definecolor{Magenta2}{rgb}{0.932, 0.0, 0.932}
\definecolor{Magenta3}{rgb}{0.804, 0.0, 0.804}
\definecolor{Magenta4}{rgb}{0.545, 0.0, 0.545}
\definecolor{Maroon}{rgb}{0.5, 0.0, 0.0}
\definecolor{Maroon0}{rgb}{0.69, 0.19, 0.376}
\definecolor{Maroon1}{rgb}{1.0, 0.204, 0.7}
\definecolor{Maroon2}{rgb}{0.932, 0.19, 0.655}
\definecolor{Maroon3}{rgb}{0.804, 0.16, 0.565}
\definecolor{Maroon4}{rgb}{0.545, 0.11, 0.385}
\definecolor{MediumAquamarine}{rgb}{0.4, 0.804, 0.668}
\definecolor{MediumBlue}{rgb}{0.0, 0.0, 0.804}
\definecolor{MediumOrchid}{rgb}{0.73, 0.332, 0.828}
\definecolor{MediumOrchid1}{rgb}{0.88, 0.4, 1.0}
\definecolor{MediumOrchid2}{rgb}{0.82, 0.372, 0.932}
\definecolor{MediumOrchid3}{rgb}{0.705, 0.32, 0.804}
\definecolor{MediumOrchid4}{rgb}{0.48, 0.215, 0.545}
\definecolor{MediumPurple}{rgb}{0.576, 0.44, 0.86}
\definecolor{MediumPurple1}{rgb}{0.67, 0.51, 1.0}
\definecolor{MediumPurple2}{rgb}{0.624, 0.475, 0.932}
\definecolor{MediumPurple3}{rgb}{0.536, 0.408, 0.804}
\definecolor{MediumPurple4}{rgb}{0.365, 0.28, 0.545}
\definecolor{MediumSeaGreen}{rgb}{0.235, 0.7, 0.444}
\definecolor{MediumSlateBlue}{rgb}{0.484, 0.408, 0.932}
\definecolor{MediumSpringGreen}{rgb}{0.0, 0.98, 0.604}
\definecolor{MediumTurquoise}{rgb}{0.284, 0.82, 0.8}
\definecolor{MediumVioletRed}{rgb}{0.78, 0.084, 0.52}
\definecolor{MidnightBlue}{rgb}{0.098, 0.098, 0.44}
\definecolor{MintCream}{rgb}{0.96, 1.0, 0.98}
\definecolor{MistyRose}{rgb}{1.0, 0.894, 0.884}
\definecolor{MistyRose1}{rgb}{1.0, 0.894, 0.884}
\definecolor{MistyRose2}{rgb}{0.932, 0.835, 0.824}
\definecolor{MistyRose3}{rgb}{0.804, 0.716, 0.71}
\definecolor{MistyRose4}{rgb}{0.545, 0.49, 0.484}
\definecolor{Moccasin}{rgb}{1.0, 0.894, 0.71}
\definecolor{NavajoWhite}{rgb}{1.0, 0.87, 0.68}
\definecolor{NavajoWhite1}{rgb}{1.0, 0.87, 0.68}
\definecolor{NavajoWhite2}{rgb}{0.932, 0.81, 0.63}
\definecolor{NavajoWhite3}{rgb}{0.804, 0.7, 0.545}
\definecolor{NavajoWhite4}{rgb}{0.545, 0.475, 0.37}
\definecolor{Navy}{rgb}{0.0, 0.0, 0.5}
\definecolor{NavyBlue}{rgb}{0.0, 0.0, 0.5}
\definecolor{OldLace}{rgb}{0.992, 0.96, 0.9}
\definecolor{Olive}{rgb}{0.5, 0.5, 0.0}
\definecolor{OliveDrab}{rgb}{0.42, 0.556, 0.136}
\definecolor{OliveDrab1}{rgb}{0.752, 1.0, 0.244}
\definecolor{OliveDrab2}{rgb}{0.7, 0.932, 0.228}
\definecolor{OliveDrab3}{rgb}{0.604, 0.804, 0.196}
\definecolor{OliveDrab4}{rgb}{0.41, 0.545, 0.132}
\definecolor{Orange}{rgb}{1.0, 0.648, 0.0}
\definecolor{Orange1}{rgb}{1.0, 0.648, 0.0}
\definecolor{Orange2}{rgb}{0.932, 0.604, 0.0}
\definecolor{Orange3}{rgb}{0.804, 0.52, 0.0}
\definecolor{Orange4}{rgb}{0.545, 0.352, 0.0}
\definecolor{OrangeRed}{rgb}{1.0, 0.27, 0.0}
\definecolor{OrangeRed1}{rgb}{1.0, 0.27, 0.0}
\definecolor{OrangeRed2}{rgb}{0.932, 0.25, 0.0}
\definecolor{OrangeRed3}{rgb}{0.804, 0.215, 0.0}
\definecolor{OrangeRed4}{rgb}{0.545, 0.145, 0.0}
\definecolor{Orchid}{rgb}{0.855, 0.44, 0.84}
\definecolor{Orchid1}{rgb}{1.0, 0.512, 0.98}
\definecolor{Orchid2}{rgb}{0.932, 0.48, 0.912}
\definecolor{Orchid3}{rgb}{0.804, 0.41, 0.79}
\definecolor{Orchid4}{rgb}{0.545, 0.28, 0.536}
\definecolor{PaleGoldenrod}{rgb}{0.932, 0.91, 0.668}
\definecolor{PaleGreen}{rgb}{0.596, 0.985, 0.596}
\definecolor{PaleGreen1}{rgb}{0.604, 1.0, 0.604}
\definecolor{PaleGreen2}{rgb}{0.565, 0.932, 0.565}
\definecolor{PaleGreen3}{rgb}{0.488, 0.804, 0.488}
\definecolor{PaleGreen4}{rgb}{0.33, 0.545, 0.33}
\definecolor{PaleTurquoise}{rgb}{0.688, 0.932, 0.932}
\definecolor{PaleTurquoise1}{rgb}{0.732, 1.0, 1.0}
\definecolor{PaleTurquoise2}{rgb}{0.684, 0.932, 0.932}
\definecolor{PaleTurquoise3}{rgb}{0.59, 0.804, 0.804}
\definecolor{PaleTurquoise4}{rgb}{0.4, 0.545, 0.545}
\definecolor{PaleVioletRed}{rgb}{0.86, 0.44, 0.576}
\definecolor{PaleVioletRed1}{rgb}{1.0, 0.51, 0.67}
\definecolor{PaleVioletRed2}{rgb}{0.932, 0.475, 0.624}
\definecolor{PaleVioletRed3}{rgb}{0.804, 0.408, 0.536}
\definecolor{PaleVioletRed4}{rgb}{0.545, 0.28, 0.365}
\definecolor{PapayaWhip}{rgb}{1.0, 0.936, 0.835}
\definecolor{PeachPuff}{rgb}{1.0, 0.855, 0.725}
\definecolor{PeachPuff1}{rgb}{1.0, 0.855, 0.725}
\definecolor{PeachPuff2}{rgb}{0.932, 0.796, 0.68}
\definecolor{PeachPuff3}{rgb}{0.804, 0.688, 0.585}
\definecolor{PeachPuff4}{rgb}{0.545, 0.468, 0.396}
\definecolor{Peru}{rgb}{0.804, 0.52, 0.248}
\definecolor{Pink}{rgb}{1.0, 0.752, 0.796}
\definecolor{Pink1}{rgb}{1.0, 0.71, 0.772}
\definecolor{Pink2}{rgb}{0.932, 0.664, 0.72}
\definecolor{Pink3}{rgb}{0.804, 0.57, 0.62}
\definecolor{Pink4}{rgb}{0.545, 0.39, 0.424}
\definecolor{Plum}{rgb}{0.868, 0.628, 0.868}
\definecolor{Plum1}{rgb}{1.0, 0.732, 1.0}
\definecolor{Plum2}{rgb}{0.932, 0.684, 0.932}
\definecolor{Plum3}{rgb}{0.804, 0.59, 0.804}
\definecolor{Plum4}{rgb}{0.545, 0.4, 0.545}
\definecolor{PowderBlue}{rgb}{0.69, 0.88, 0.9}
\definecolor{Purple}{rgb}{0.5, 0.0, 0.5}
\definecolor{Purple0}{rgb}{0.628, 0.125, 0.94}
\definecolor{Purple1}{rgb}{0.608, 0.19, 1.0}
\definecolor{Purple2}{rgb}{0.57, 0.172, 0.932}
\definecolor{Purple3}{rgb}{0.49, 0.15, 0.804}
\definecolor{Purple4}{rgb}{0.332, 0.1, 0.545}
\definecolor{Red}{rgb}{1.0, 0.0, 0.0}
\definecolor{Red1}{rgb}{1.0, 0.0, 0.0}
\definecolor{Red2}{rgb}{0.932, 0.0, 0.0}
\definecolor{Red3}{rgb}{0.804, 0.0, 0.0}
\definecolor{Red4}{rgb}{0.545, 0.0, 0.0}
\definecolor{RosyBrown}{rgb}{0.736, 0.56, 0.56}
\definecolor{RosyBrown1}{rgb}{1.0, 0.756, 0.756}
\definecolor{RosyBrown2}{rgb}{0.932, 0.705, 0.705}
\definecolor{RosyBrown3}{rgb}{0.804, 0.608, 0.608}
\definecolor{RosyBrown4}{rgb}{0.545, 0.41, 0.41}
\definecolor{RoyalBlue}{rgb}{0.255, 0.41, 0.884}
\definecolor{RoyalBlue1}{rgb}{0.284, 0.464, 1.0}
\definecolor{RoyalBlue2}{rgb}{0.264, 0.43, 0.932}
\definecolor{RoyalBlue3}{rgb}{0.228, 0.372, 0.804}
\definecolor{RoyalBlue4}{rgb}{0.152, 0.25, 0.545}
\definecolor{SaddleBrown}{rgb}{0.545, 0.27, 0.075}
\definecolor{Salmon}{rgb}{0.98, 0.5, 0.448}
\definecolor{Salmon1}{rgb}{1.0, 0.55, 0.41}
\definecolor{Salmon2}{rgb}{0.932, 0.51, 0.385}
\definecolor{Salmon3}{rgb}{0.804, 0.44, 0.33}
\definecolor{Salmon4}{rgb}{0.545, 0.298, 0.224}
\definecolor{SandyBrown}{rgb}{0.956, 0.644, 0.376}
\definecolor{SeaGreen}{rgb}{0.18, 0.545, 0.34}
\definecolor{SeaGreen1}{rgb}{0.33, 1.0, 0.624}
\definecolor{SeaGreen2}{rgb}{0.305, 0.932, 0.58}
\definecolor{SeaGreen3}{rgb}{0.264, 0.804, 0.5}
\definecolor{SeaGreen4}{rgb}{0.18, 0.545, 0.34}
\definecolor{Seashell}{rgb}{1.0, 0.96, 0.932}
\definecolor{Seashell1}{rgb}{1.0, 0.96, 0.932}
\definecolor{Seashell2}{rgb}{0.932, 0.898, 0.87}
\definecolor{Seashell3}{rgb}{0.804, 0.772, 0.75}
\definecolor{Seashell4}{rgb}{0.545, 0.525, 0.51}
\definecolor{Sienna}{rgb}{0.628, 0.32, 0.176}
\definecolor{Sienna1}{rgb}{1.0, 0.51, 0.28}
\definecolor{Sienna2}{rgb}{0.932, 0.475, 0.26}
\definecolor{Sienna3}{rgb}{0.804, 0.408, 0.224}
\definecolor{Sienna4}{rgb}{0.545, 0.28, 0.15}
\definecolor{Silver}{rgb}{0.752, 0.752, 0.752}
\definecolor{SkyBlue}{rgb}{0.53, 0.808, 0.92}
\definecolor{SkyBlue1}{rgb}{0.53, 0.808, 1.0}
\definecolor{SkyBlue2}{rgb}{0.494, 0.752, 0.932}
\definecolor{SkyBlue3}{rgb}{0.424, 0.65, 0.804}
\definecolor{SkyBlue4}{rgb}{0.29, 0.44, 0.545}
\definecolor{SlateBlue}{rgb}{0.415, 0.352, 0.804}
\definecolor{SlateBlue1}{rgb}{0.512, 0.435, 1.0}
\definecolor{SlateBlue2}{rgb}{0.48, 0.404, 0.932}
\definecolor{SlateBlue3}{rgb}{0.41, 0.35, 0.804}
\definecolor{SlateBlue4}{rgb}{0.28, 0.235, 0.545}
\definecolor{SlateGray}{rgb}{0.44, 0.5, 0.565}
\definecolor{SlateGray1}{rgb}{0.776, 0.888, 1.0}
\definecolor{SlateGray2}{rgb}{0.725, 0.828, 0.932}
\definecolor{SlateGray3}{rgb}{0.624, 0.712, 0.804}
\definecolor{SlateGray4}{rgb}{0.424, 0.484, 0.545}
\definecolor{SlateGrey}{rgb}{0.44, 0.5, 0.565}
\definecolor{Snow}{rgb}{1.0, 0.98, 0.98}
\definecolor{Snow1}{rgb}{1.0, 0.98, 0.98}
\definecolor{Snow2}{rgb}{0.932, 0.912, 0.912}
\definecolor{Snow3}{rgb}{0.804, 0.79, 0.79}
\definecolor{Snow4}{rgb}{0.545, 0.536, 0.536}
\definecolor{SpringGreen}{rgb}{0.0, 1.0, 0.498}
\definecolor{SpringGreen1}{rgb}{0.0, 1.0, 0.498}
\definecolor{SpringGreen2}{rgb}{0.0, 0.932, 0.464}
\definecolor{SpringGreen3}{rgb}{0.0, 0.804, 0.4}
\definecolor{SpringGreen4}{rgb}{0.0, 0.545, 0.27}
\definecolor{SteelBlue}{rgb}{0.275, 0.51, 0.705}
\definecolor{SteelBlue1}{rgb}{0.39, 0.72, 1.0}
\definecolor{SteelBlue2}{rgb}{0.36, 0.675, 0.932}
\definecolor{SteelBlue3}{rgb}{0.31, 0.58, 0.804}
\definecolor{SteelBlue4}{rgb}{0.21, 0.392, 0.545}
\definecolor{Tan}{rgb}{0.824, 0.705, 0.55}
\definecolor{Tan1}{rgb}{1.0, 0.648, 0.31}
\definecolor{Tan2}{rgb}{0.932, 0.604, 0.288}
\definecolor{Tan3}{rgb}{0.804, 0.52, 0.248}
\definecolor{Tan4}{rgb}{0.545, 0.352, 0.17}
\definecolor{Teal}{rgb}{0.0, 0.5, 0.5}
\definecolor{Thistle}{rgb}{0.848, 0.75, 0.848}
\definecolor{Thistle1}{rgb}{1.0, 0.884, 1.0}
\definecolor{Thistle2}{rgb}{0.932, 0.824, 0.932}
\definecolor{Thistle3}{rgb}{0.804, 0.71, 0.804}
\definecolor{Thistle4}{rgb}{0.545, 0.484, 0.545}
\definecolor{Tomato}{rgb}{1.0, 0.39, 0.28}
\definecolor{Tomato1}{rgb}{1.0, 0.39, 0.28}
\definecolor{Tomato2}{rgb}{0.932, 0.36, 0.26}
\definecolor{Tomato3}{rgb}{0.804, 0.31, 0.224}
\definecolor{Tomato4}{rgb}{0.545, 0.21, 0.15}
\definecolor{Turquoise}{rgb}{0.25, 0.88, 0.815}
\definecolor{Turquoise1}{rgb}{0.0, 0.96, 1.0}
\definecolor{Turquoise2}{rgb}{0.0, 0.898, 0.932}
\definecolor{Turquoise3}{rgb}{0.0, 0.772, 0.804}
\definecolor{Turquoise4}{rgb}{0.0, 0.525, 0.545}
\definecolor{Violet}{rgb}{0.932, 0.51, 0.932}
\definecolor{VioletRed}{rgb}{0.816, 0.125, 0.565}
\definecolor{VioletRed1}{rgb}{1.0, 0.244, 0.59}
\definecolor{VioletRed2}{rgb}{0.932, 0.228, 0.55}
\definecolor{VioletRed3}{rgb}{0.804, 0.196, 0.47}
\definecolor{VioletRed4}{rgb}{0.545, 0.132, 0.32}
\definecolor{Wheat}{rgb}{0.96, 0.87, 0.7}
\definecolor{Wheat1}{rgb}{1.0, 0.905, 0.73}
\definecolor{Wheat2}{rgb}{0.932, 0.848, 0.684}
\definecolor{Wheat3}{rgb}{0.804, 0.73, 0.59}
\definecolor{Wheat4}{rgb}{0.545, 0.494, 0.4}
\definecolor{White}{rgb}{1.0, 1.0, 1.0}
\definecolor{WhiteSmoke}{rgb}{0.96, 0.96, 0.96}
\definecolor{Yellow}{rgb}{1.0, 1.0, 0.0}
\definecolor{Yellow1}{rgb}{1.0, 1.0, 0.0}
\definecolor{Yellow2}{rgb}{0.932, 0.932, 0.0}
\definecolor{Yellow3}{rgb}{0.804, 0.804, 0.0}
\definecolor{Yellow4}{rgb}{0.545, 0.545, 0.0}
\definecolor{YellowGreen}{rgb}{0.604, 0.804, 0.196}

\newcommand{\vect}{\overset{\rightharpoonup}}

\begin{document}

\title{Noise calibration for the stochastic rotating shallow water model}
\author{Dan Crisan$^\ast$\and
        Oana Lang$^\ast$\and
        Alexander Lobbe$^\ast$\and
        Peter Jan van Leeuwen$^\dagger$\and
        Roland Potthast$^{\ddagger,}$ $^\diamond$}
\date{{\footnotesize{\raggedright
$^\ast$Department of Mathematics, Imperial College London, 180 Queen's Gate, SW7 2AZ, London, UK\\
$^\dagger$Department of Atmospheric Science, Colorado State University, Fort Collins, Colorado 80523 USA and Department of Meteorology, University of Reading, Earley Gate,  RG6 7BE, UK\\
$^\ddagger$Deutscher Wetterdienst (DWD), Frankfurter Strasse 135, 63067 Offenbach, Germany\\
$^\diamond$Dept.~of Mathematics and Statistics, University of Reading, Whiteknights, PO Box 220, Reading RG6 6AX, UK\\
}}
\vspace{3ex}
%\date{}
}

\maketitle
%%%%%%%%%%%%%%%%%%%%%%%%%%%%%%%%%%%%%%%%%%%%%%%%%%

\abstract{
%%%%%%%%%%%%%%%%%%%%%%%%%%%%%%%%%%%%%%%%%%%%%%%%%%
Stochastic partial differential equations have been used in a variety of contexts to model the evolution of uncertain dynamical systems. In recent years, their applications to geophysical fluid dynamics has increased massively. For a judicious usage in modelling fluid evolution, one needs to calibrate the amplitude of the noise to data. In this paper we address this requirement for the stochastic rotating shallow water (SRSW) model. This work is a continuation of \cite{lucrisanlangmemin}, where a data assimilation methodology has been introduced for the SRSW model. The noise used in \cite{lucrisanlangmemin} was introduced as an arbitrary random phase shift in the Fourier space. This is not necessarily consistent with the uncertainty induced by a model reduction procedure. In this paper, we introduce a new method of noise calibration of the SRSW model which is compatible with the model reduction technique. The method is generic and can be applied to arbitrary stochastic parametrizations. It is also agnostic as to the source of data (real or synthetic). It is based on a principal component analysis technique to generate the eigenvectors and the eigenvalues of the covariance matrix of the stochastic parametrization. For SRSW model covered in this paper, we calibrate the noise by using the elevation variable of the model, as this is an observable easily obtainable in practical application, and use synthetic data as input for the calibration procedure. 
}

\clearpage
\tableofcontents
%%%%%%%%%%%%%%%%%%%%%%%%%%%%%%%%%%%%%%%%%%%%%%%%%%

\section{Introduction}

Stochastic parameterisations model the uncertainty caused by unknown or neglected physical effects, incomplete  or inaccurate information in both observational data and the formulation of the theoretical models used for prediction. The topic 
of stochastic parameterisations has been a very active area of research in the last two decades, partly due to their use in modelling uncertainty generated by the  reduction of  high-resolution solutions to coarser scale. 
In recent years, several stochastic parametrizations have been proposed to model this type of uncertainty, see e.g. \cite{EtienneLong}, \cite{Wei1}, \cite{Holm2015}, \cite{Memin2014}, \cite{Wei3}.
        
The correct calibration of the stochastic model parameters is crucial for ensuring the effective application of the combined stochastic parameterisation and data assimilation procedure. Various numerical methods for calibration, see \cite{Wei1}, \cite{EtienneLong}, \cite{Wei3}, \cite{Resseguir2021} have been implemented to show that data driven models and state of the art data assimilation techniques can be successfully combined.    

In this paper, we introduce a new methodology to calibrate a stochastic partial differential equation where the stochasticity accounts for small-scale effects which are missed as a result of working with models run at coarse resolution. This is part of the current efforts aimed at designing data-driven models in which real uncertainty is accounted for based on input from measurements and statistically-informed initial data. 
%OK \todo[inline]{We need to choose a consistent spelling for "parametrisation".}
We give next a brief description of the stochastic parametrization framework and the calibration methodology. We highlight that the procedure developed here is generic. In other words, it can be applied to any model satisfying (\ref{det}) and any stochastic parametrization satisfying (\ref{sto}).  

We denote by $m^f$ the model state. We will assume that the evolution of $m^f$ is governed by a partial differential equation   
\begin{equation}\label{det}
\frac{dm^f}{dt}=\mathcal{A}(m^f), \ \ t\ge 0,
\end{equation}
where $\mathcal{A}$ is the model operator. We will denote by $m^c$ the coarse scale model. The effect of the \textit{unresolved scales} can be mathematically modelled by a term of the form 
\begin{equation}\label{det}
\int_0^t \mathcal{M}(m_s^c)dW_s,
\end{equation}
where $\mathcal{M}$ is a suitably chosen operator   
and $W_t=W(t,x)$ is a space-time Brownian motion. In other words, the model run on the coarse scale will 
satisfy the stochastic (partial) differential equation 
\begin{equation}\label{sto}
dm^c=\mathcal{A}(m^c)dt+\mathcal{M}(m^c)dW_t, \ \ t\ge 0,
\end{equation}
We note that, theoretically, the solutions of both the deterministic equation 
(\ref{det}) and its stochastic counter-part (\ref{sto}) live on the same physical domain ($\mathbb{R}^n$, the torus, a horizontal strip, etc). \emph{Numerically} they are approximated on different space grids: the space discretisation for (\ref{det}) is finer than the one for (\ref{sto}). If we are to refer to the numerical resolution for (\ref{det}) and (\ref{sto}), then we 
could distinguish between the model operator for (\ref{det}), and call it, say 
$\mathcal{A}^f$ and that for (\ref{sto}), and call it, say, $\mathcal{A}^c$.

%\DC{DC: How is the above paragraph Peter Jan}
%\pj{PJ: This $\mathcal{A}$ is a different operator from that in equation (1), so it needs a different symbol. 
%I suggest to use f and c superindices, i.e. $\mathcal{A}^c$, and similar of the fine model.} 

For our analysis, it is convenient to decompose $W$ as follows
\begin{equation}
W_t = \displaystyle\sum_{k=1}^{\infty} \xi_k W_t^k
\end{equation}
where the coefficients $(\xi_k)$ 
%\pj{PJ: I don't understand the two k indices, should that be just one?} \DC{DC: there is one eigenvector/eigenvalue per Brownian motion/source of noise}  
are space-dependent vector fields and the processes $(W^k)$ are one-dimensional standard Brownian motions.

The challenge that we address in this paper is to obtain an approximate representation of uncertainty that uses a finite set of individual noises $W_k$ and
calibrate the amplitude of these individual noises 
%(\pj{PJ: suggest to remove this:}) 
(in other words, the vector fields $\xi_k$). More precisely, we will model the uncertainty by using 
\begin{equation}\label{repres}
\displaystyle\int_0^t \mathcal{M}(m_s^c)dW_s^N
\end{equation}
on a sufficiently large time window $[0,T]$
where
\begin{equation}
W_t^N = \displaystyle\sum_{k=1}^N \xi_kW_t^k
\end{equation}
and both $\xi_k$ and $N$ have to be estimated from data. Once the calibration procedure is complete, we can simulate the calibrated stochastic equation (\ref{sto}) on the coarse scale instead of the deterministic equation (\ref{det}) on the refined scale. This model reduction can lead to a significant reduction of computational effort.

%Currently, we save computational effort proportional to a factor of $C^2$ where $C$ is the coarsening factor\footnote{In our analysis $C=4$.}. \pj{PJ: I suggest to remove the footnote, see above. Furthermore, the power of C depends on the dimension of the model. Typically, for a 2-D model we gain a factor close to $C^3$ because the time step is also decreased.}

The calibration procedures introduced below is agnostic to the source of the input data.
It can be real data, such as satellite observations of e.g. the ocean sea-surface height, data from re-analysis such as ERA5 (\cite{longlist}), or synthetic data from a model run of (\ref{det}) computed on a sufficiently large time window $[0,T]$. 
The model run is then mollified using a procedure that will eliminate the small/fast scales effects, for example by using a low-pass filter, Gaussian mollifier, Helmholtz projection, subsampling, etc; or combinations thereof. We will denote by $C(m^f)$ the resulting mollification of the data. Note that both $m^f$ and $C(m^f)$ live on the same space. We emphasise that $C(m^f)$ is not the solution of (\ref{sto}). However, we make the ansatz that the difference between the two processes 
$\hat{m} := m^f - C(m^f)$ has a stochastic representation given by (\ref{repres})
and, therefore, on a sufficiently small time interval $[t, t+\delta]$ we will have
\begin{equation}\label{diff2}
\hat{m}_{t+\delta} - \hat{m}_{t} = \displaystyle\int_t^{t+\delta}\mathcal{M}(C(m_s^f))dW_s^N 
 \approx \displaystyle\sum_{k=1}^N\mathcal{M}(C(m_t^f))\xi_k\Delta W_t^k.
\end{equation}

Note that in the representation (\ref{diff2}), we work with an approximation as the difference $\hat{m}_{t+\delta} - \hat{m}_{t}$ has leading term $\sum_{k=1}^n\mathcal{M}(C(m_t^f))\xi_k\Delta W_t^k$ with error term of order $\mathcal{O}(\delta)$ or no error term at all in the case when $\mathcal{M}(C(m^f))$ is constant (i.e. the additive noise case).

So for a suitably chosen partition  of the time interval $[0,T]$ of the form $0\leq t_1 < t_1+\delta < \ldots t_n < t_n + \delta \leq T$ (see details below) we need to estimate the vector fields $\xi_k$ from the data: 
\begin{equation}
\hat{m}_{t_k+\delta} - \hat{m}_{t_k}, \ \ k=2,...,n
\end{equation}
This will depend on the choice of  perturbation chosen to account for the model uncertainty. Here are some examples:
\begin{itemize}
    \item Additive noise, e.g. \cite{lucrisanlangmemin}. In this simple case, $\mathcal{M}(C(m_s^f))\equiv 1$, in other words,
    \begin{equation}\label{add}
    \hat{m}_{t+\delta} - \hat{m}_{t}   
    \approx \displaystyle\sum_{k=1}^N\xi_k\Delta W_t^k.
    \end{equation}
    and therefore $\hat{m}_{t_k+\delta} - \hat{m}_{t_k}$
    can be interpreted as samples from a multi-dimensional Gaussian random variable $M$ where 
    \begin{equation}\label{gr}
    M \sim \mathcal{N}\left( 0, \displaystyle\sum_{k=1}^N
    (\xi_k)(\xi_k)^Tt\right). 
    \end{equation}
    \item Multiplicative noise (e.g. \cite{pathiraja}). In this case,  
    $\mathcal{M}(C(m_t^f))=f(C(m_t^f))$ in other words,
    \begin{equation}\label{mult}
    \hat{m}_{t+\delta} - \hat{m}_{t} 
    \approx f(C(m_t^f)) \displaystyle\sum_{k=1}^N\xi_k\Delta W_t^k.
    \end{equation}
    in which $f(C(m_t^f))$ is a scalar function of $C(m_t^f)$, and therefore $\frac{\hat{m}_{t_k+\delta} - \hat{m}_{t_k}}{f(C(m_t^f))}$
    can be interpreted as samples from a multi-dimensional Gaussian random variable $M$ with distribution (\ref{gr}) %\pj{(PJ:Only if m is a scalar, otherwise we divide by a vector. Or do we mean element-wise? )}.\DC{DC: Yes scalar.} Of course we would need to assume that  $f(C(m^f_s))\neq 0$.
    %\item Transport noise. We describe for a particular class of transport noise stochastic parametrizations that cover the 2D incompressible Euler \cite{Wei1}, 2-layer quasi-geostrophic model \cite{Wei2} (with fluid vorticity being the state model in both cases), rotating shallow water equation \
    \item Transport noise. We describe the procedure for a particular class of transport noise stochastic parametrizations that cover the 2D incompressible Euler \cite{Wei1}, 2-layer quasi-geostrophic model \cite{Wei2} (with fluid vorticity being the model state in both cases), 2D thermal quasigeostrophic (TQG) model (with fluid buoyancy being the model state in this case), rotating shallow water model \cite{lucrisanlangmemin}, \cite{Holm2015}, \cite{lucrisanlangmemin} (with fluid elevation being the model state in this case). In all of these cases, the model state is defined over
    ${\mathbb R}^2$ (or a subset of ${\mathbb R}^2$) and the ansatz is that  
    \begin{eqnarray}
    \hat{m}_{t+\delta} - \hat{m}_{t}   
    &\approx&  \displaystyle\sum_{k=1}^N(\xi_k^1\partial^1C(m_s^f)+
    \xi_k^2\partial^2C(m_s^f)) \Delta W_t^k.\nonumber\\
    &=& \displaystyle\sum_{k=1}^N(\partial^2\psi_k\partial^1C(m_s^f)-
    \partial^1\psi_k\partial^2C(m_s^f)) \Delta W_t^k \nonumber\\
    &=& \left(-\partial^2C(m_s^f) \partial^1 + \partial^1C(m_s^f) \partial^2 \right)\left (\sum_{k=1}^N\psi_k\Delta W_t^k \right ),\label{transport}
    \end{eqnarray}
    where $\partial^{i}$ is the standard partial derivative and we imposed the additional assumption that $\xi_k\equiv\grad^\perp \psi_k$. This will ensure that $\div \xi_k=0$,  which is an assumption made in all the above models. It follows from (\ref{transport}) that the solution of the linear hyperbolic equation  
 \begin{equation}\label{ceh}
     f = \vect{q} \cdot \nabla \psi
 \end{equation}
 with $f:=-(\hat{m}_{t+\delta} - \hat{m}_{t})$ and $\vect{q}:=-\nabla^\perp C(m_t^f)$ can be interpreted as samples from the multi-dimensional Gaussian random variable
\begin{equation}
M \sim \mathcal{N}\left( 0, \displaystyle\sum_{k=1}^N(\psi_k)(\psi_k)^Tt\right).
\end{equation}
This interpretation is the basis of the calibration procedure which is detailed in Section \ref{section:calibration}.
\end{itemize}

%\noindent \DC{DC: one or two references for fluid models with additive/multiplicative noise would be good here. I added a sentence about the Lagrangian approach in section 1.1, not keen to expand on this. If we do the referees might ask us to compare the two.} 
%\pj{PJ: Would these work? Agree on Lagrangian.}    
   
%\OL{OL: Here I'd add a more extensive literature review and explain on stochastic parameterisations and compare the method with the Lagrangian approach: i.e. we don't need to simulate Lagrangian trajectories corresponding to the velocity vector field etc.}

%\OL{To simulate the evolution of the atmospheric Jetstream, the model domain corresponds to a strip situated between 30 and 60 degrees north latitude, with east-west periodic boundary conditions and free-slip boundary conditions on the northern and southern domain.}

\subsection{Calibration methodology for the Rotating Shallow Water Models}\label{cals}

We apply the new calibration procedure for a class of stochastic parametrizations for the rotating shallow water model, see equation (\ref{h}) below. Here, we use as calibration data input from the elevation variable $h$ of the rotating shallow water system and we will assume that the stochastic parametrization for $h$ is given by equation (\ref{model1}) below. This form for the $h$ equation is common to 
both \cite{lucrisanlangmemin}, \cite{lucrisanlangmemin}. The model in 
\cite{lucrisanlangmemin} is based on the Location Uncertainty paradigm first introduced in \cite{Memin2014}, whilst the model in \cite{lucrisanlangmemin} is based on the Stochastic Advection by Lie Trasport paradigm  first introduced in \cite{Holm2015}.

The rotating shallow water model is classically known (see e.g. \cite{Kalnay}, \cite{Vallis}, \cite{Zeitlinbook}) for its complex structure which captures important aspects of the oceanic and the atmospheric dynamics, such as potential vorticity and energy conservation, or the existence of gravity waves. These features make it one of the favourite systems for modelling geophysical turbulence. Further, the stochastic parametrizations introduced in \cite{Holm2015} and \cite{Memin2014}, see also \cite{lucrisanlangmemin} and \cite{lucrisanlangmemin}, offer possible ways to account for small-scale/fast scale components that remain unresolved when running the model state on coarser grids. In particular, these stochastic parametrization schemes are both mathematically derived from, as well as physically consistent with, the underlying deterministic/fine-scale dynamics.

To estimate the stochastic noise, we will use synthetic data as input for the calibration procedure. In particular, we will use a simulation of the rotating shallow water model run on a fine rectangular grid. 
%$2224\times 320$\todo[inline]{I would prefer to reserve the introduction of specific numerical parameters for the experiment section. Maybe we can be more generic here and mention the physical dimensions of the domain (it is some strip on the earth surface) and then mention two degrees of coarsening.}. 
In the above, we used the generic notation $m^f$ in the for this model run.
We will calibrate stochastic parametrizations corresponding to two coarser grids, one of size $556\times 80$ and one of size $256\times 40$. We then construct a mollified version of the fine grid trajectory using a low-pass filter. We refer to the mollified solution as the \textit{coarse grid trajectory}. Alternative coarsenings may be performed using Gaussian mollifiers or subsampling techniques. In the above, we used the generic notation $C(m^f)$ in the for this mollified solution. This is done by using the following steps: 
\begin{itemize}
\item [1.] Compute the time-increments of the discrepancy between the fine resolution and the coarse resolution trajectories.
\item [2.] Choose a \textit{calibration time grid} in such a way that the data is decorrelated. For this, we first estimate the \textit{decorrelation time} of the data corresponding to the fine grid trajectory.
\item [3.] Compute the sample noises corresponding to each of the times on the calibration time grid.  In our case, this amounts to solving a sequence of hyperbolic equations of the form (\ref{ceh}), see Section \ref{section:calibration} for details.

%\todo[inline]{No need to mention firedrake here, in my opinion}{}(the hyperbolic equations are solved using Firedrake).

\item [4.] Extract a basis for the stochastic noise together with the corresponding eigenvalues that explains a sufficiently large part of the variance in the data. For this we use the Principal Components Analysis algorithm used in \cite{Wei1} . 
\end{itemize}

% DC: I commented out this paragraph as this is now covered extensively in the introduction. We prove that this approach can be effectively utilized to calibrate stochastic shallow water models. Our implementation focuses on a SRSW model which can be derived within the SALT framework. However, the calibration method we propose does not rely on the particular structure of the stochastic term which drives the signal process. Actually, it does not rely on the particular form of the PDE itself, as long as the stochastic parameters can be written using a potential function. 

As explained above, the calibration procedure presented here is generic and can be used for a variety of stochastic parametrizations. In \cite{Wei1} and \cite{Wei2} a different calibration procedure was implemented. The calibration procedure in these works uses the specificity of the stochastic advection by Lie transport (SALT) models and cannot be applied to non-SALT stochastic parametrizations.

The methodology described above is used as a means of accounting for the resolution error in numerical simulations.  We evaluate it by performing a number of uncertainty quantification tests, see Section \ref{spc} for details. 
 
We perform the uncertainty quantification (UQ) analysis by substituting the parameters estimated above and simulating an ensemble of  \textit{particles} associated with the SPDE trajectory. The aim is to use the ensemble generated with the chosen
stochastic parametrization to quantify the uncertainty introduce when moving from the fine grid to the two coarser grids (the grids used to run the ensemble, and which will be eventually used to run the Data Assimilation methodology).

\subsection{Structure of the Paper} \label{spc}

We complete the introduction with a description of the contents of the paper: 
\begin{itemize}
\item In Section \ref{sect:srswmodel} we introduce the stochastic rotating shallow water model in both deterministic and stochastic setting. The input data for the calibration procedure will be extracted from the deterministic version run on a fine resolution grid.  

\item Section \ref{section:calibration} contains a detailed description of the calibration methodology for the stochastic model, which is tested on the stochastic rotating shallow water equations. The calibration is performed following four main steps which are presented thoroughly: mollification of the fine grid data, data generation, solving the calibration equation, extracting a basis for the stochastic noise. 
\item In Section \ref{nr} we perform a series of uncertainty quantification tests to validate the calibration procedure and the stochastic parametrisation. The tests are run for elevation, zonal velocity and meridional velocity, and they show that this calibration methodology can be efficiently used to estimate spatial correlations from synthetic data. The ensemble of stochastic trajectories captures the large-scale behaviour of the underlying \textit{truth}.
\item We conclude the paper with Section \ref{sect:conclusions} which summarizes our results and introduces possible further research directions, and an Appendix where we include technical details related to the numerical discretisation and implementation of the stochastic model.  
\end{itemize}

\section{The Rotating Shallow Water Model}\label{sect:srswmodel}
%%%%%%%%%%%%%%%%%%%%%%%%%%%%%%%%%%%%%%%%%%%%%%%%%%
%%%%%%%%%%%%%%%%%%%%%%%%%%%%%%%%%%%%%%%%%%%%%%%%%%
%%%%%%%%%%%%%%%%%%%%%%%%%%%%%%%%%%%%%%%%%%%%%%%%%%
The rotating shallow water model describes the evolution of a
compressible rotating fluid below a free surface. The typical vertical length scale is assumed to be much smaller than the horizontal one, which gives the \textit{shallow} aspect, as one can see in Figure \ref{fig:sw-variables}. The inviscid version of the model is given by a horizontal momentum equation and a mass continuity equation
which, in the presence of rotation, can be described as follows (see \cite{Vallis}):

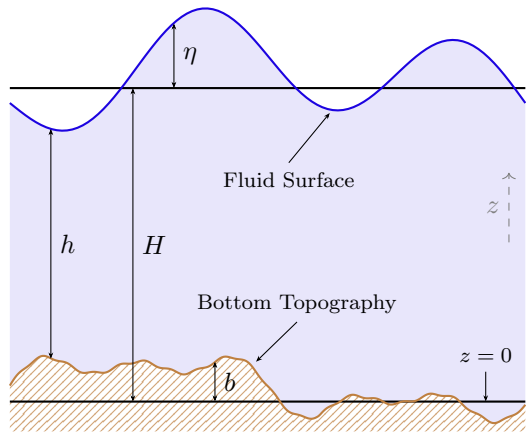
\begin{figure}[h]
    \centering
    \begin{tikzpicture}[] 
    \begin{axis}[xmin = 0.0, xmax = 2*pi, ymin = -0.3, ymax = 4, trig format = rad,
    axis y line = middle, axis x line = none, 
    %y axis line style={color=gray}, 
    ticks=none,  
    axis y line = none,
    %ylabel=$z$, y label style={at={(current axis.above origin)},anchor=south east, color=gray}
    ]
     	\addplot[domain={0.0:2*pi}, samples=2, thick]{0};
        \addplot[color=black, domain={0.0:2*pi}, samples=2, thick]{5.15-2};
        \addplot[color=WaterBlue, domain={0.0:2*pi}, samples=100, thick, name path = fluid]{(0.1*sin(x)+sin(x/2)+sin(2*x+pi))/2.1+5-2};
        \addplot[color=brown, domain={0.0:2*pi}, samples=300, thick, name path = bottom]{(0.2*abs(x-1.2*pi))+(sin(0.1*x)+sin(0.5*x)+sin(x)+sin(5*x)+sin(3*x))/10-0.1-((x-pi)^2)/20+sin(20*x)/70};
        \addplot[fill=WaterBlue, fill opacity = 0.1] fill between[of=fluid and bottom, soft clip = {domain={0:2*pi}}];
        \path[name path=xAxis] (axis cs:0,-0.3) -- (axis cs:2*pi,-0.3);
        \addplot[pattern color=brown!70!white, pattern=north east lines] fill between[of=bottom and xAxis, soft clip = {domain={0:2*pi}}];
        \draw[{Stealth[scale=0.5]}-{Stealth[scale=0.5]}] (axis cs:1.5,0) -- (axis cs:1.5,3.15) node[midway, right]{$H$};
        \draw[{Stealth[scale=0.5]}-{Stealth[scale=0.5]}] (axis cs:0.5,0.43) -- (axis cs:0.5,2.74) node[midway, right]{$h$};
        \draw[{Stealth[scale=0.5]}-{Stealth[scale=0.5]}] (axis cs:2,3.15) -- (axis cs:2,3.80) node[midway, right]{$\eta$};  
        \draw[{Stealth[scale=0.5]}-{Stealth[scale=0.5]}] (axis cs:2.5,0) -- (axis cs:2.5,0.4) node[midway, right]{$b$}; 
        \draw[-{Stealth[scale=0.5]}, thin]  (axis cs:3.4,2.4) node[below]{\footnotesize Fluid Surface} -- (axis cs:3.9,2.9) ;
        \draw[-{Stealth[scale=0.5]}, thin]  (axis cs:3.5,0.8) node[above]{\footnotesize Bottom Topography} -- (axis cs:3.0,0.4) ;
        \draw[->, color=gray, dashed] (axis cs:2*pi-0.2,1.6) -- (axis cs:2*pi-0.2,2.3) node[midway, left]{$z$} ;
        \draw[-{Stealth[scale=0.5]}, thin] (axis cs:5.8,0.3) node[above]{\footnotesize $z=0$} -- (axis cs:5.8,0.0)  ;
    \end{axis}
\end{tikzpicture}
    \caption{Illustration of the variables in the shallow water model. The scalar function $h$ denotes the height of the fluid column and the scalar value $H$ is the average height of the fluid column over the domain. We denote by $\eta$ the scalar function that gives the elevation of the fluid surface relative to $H$ and the scalar function $b$ is the bottom topography. Thus, the $z$-coordinate of the fluid surface is given by $H+\eta=h+b$.}
    \label{fig:sw-variables}
\end{figure}

    \begin{subequations}\label{h}
    \begin{equation} 
    \frac{D}{Dt}u_t + f\hat{z}\times u_{t}+ g\nabla h_{t} =0 
    \end{equation}
    \begin{equation}
    \frac{\partial h_{t}}{\partial t}+\nabla \cdot (h_{t}u_{t}) =0 ,
    \end{equation}
    \end{subequations}
    where
    \begin{itemize}
    \item $\frac{D}{Dt} := \frac{\partial}{\partial t} + u \cdot \nabla$ is the
    material derivative.
    
    \item $u=(u^1,u^2)^T$ is the horizontal fluid velocity vector field
    \item $h$ is the thickness of the fluid column (total depth)
    \item $f$ is the Coriolis parameter, $f=2\Theta \sin\varphi$ where $\Theta$
    is the rotation rate of the Earth and $\varphi$ is the latitude; $f\hat{z}
    \times u = (-fu^2, fu^1)^T$, where $\hat{z}$ is a unit vector pointing away
    form the centre of the Earth
    \item $g$ is the gravitational acceleration.
    \end{itemize}

    We can formally re-write a viscous version of the RSW system: denote by $X:=\left(u,h\right)^T $ and then \footnote{%
    We use here the differential notation to match the stochastic version (\ref%
    {shortsrsw}).}
    \begin{equation}
    dX_{t}+F\left( X_{t}\right)dt =0,
    \end{equation}%
    where $F\left(X_{t}\right) $ denotes 
    \begin{equation}\label{advectiveF}
    F\left( 
    \begin{array}{c}
    u \\ 
    h%
    \end{array}
    \right) =\left( 
    \begin{array}{c}
    u \cdot \nabla u + f\hat{z} \times u + g\nabla h -
    \nu \Delta u\\ 
    \nabla \cdot (hu) -
    \eta \Delta h
    \end{array}
    \right).
    \end{equation}

Different methods for introducing stochasticity to model uncertainty have been analysed in the literature, see e.g. \cite{Palmer}, \cite{MajdaDA} or \cite{Buizza}. In general, the noise is introduced into the forcing part of the signal process. The stochastic model used in this paper has been introduced in \cite{Holm2015}, 
using a new approach to subgrid transport modelling, called Stochastic Advection by Lie Transport
(SALT). In SALT, stochasticity is introduced into the advective part of the model equation so that the resulting stochastic system models the uncertain \emph{transport} behaviour:
    
    \begin{subequations} \label{model1}
     \begin{equation} %\label{model1} 
     du_t + \big[u_t \cdot \nabla u_t + f\hat{z} \times u_t + g\nabla h_t \big]dt + \displaystyle\sum_{i=1}^{\infty}\big[(\mathcal{L}_i + \mathcal{A}_i)u_t\big] \circ dW_t^i = \nu\Delta u_t dt
     \end{equation}
     \begin{equation} \label{model2} 
      dh_t + \nabla \cdot (h_t u_t)dt+\displaystyle\sum_{i=1}^{\infty}\big[\nabla \cdot (\xi_i h_t)\big]\circ dW_t^i= \eta\Delta h_t dt
     \end{equation}
    \end{subequations}
where $\xi _{i}$ are divergence-free and time-independent vector fields,
    $\mathcal{L}_iu := \xi_i \cdot \nabla u, \ \mathcal{A}_iu:= u_j\nabla \xi_i^j = \displaystyle\sum_{j=1}^2 u_j\nabla\xi_i^j$.
    
Similar to the deterministic case, we can recast the stochastic model as:
    \begin{equation}  \label{shortsrsw}
    dX_{t}+F\left( X_{t}\right)dt +\sum_{i=1}^{\infty }\mathcal{G}_{i}\left(
    X_{t}\right) \circ dW_{t}^{i}=\gamma \Delta X_{t}dt
    \end{equation}% 
where $\gamma = (\nu,\eta) $ is positive and corresponds to the eddyviscosity. Different levels of viscosity for the different components of $X$ can be
treated in the same manner.  $W^{i}$ are independent Brownian motions, $F$ is the nonlinear advective term defined in \eqref{advectiveF}, and $\mathcal{G}_{i}$ are differential
    operators:
    \begin{equation*}
\mathcal{G}_{i}(X)=\mathcal{G}_{i}\left( 
\begin{array}{c}
u \\ 
h%
\end{array}%
\right) =\left( 
\begin{array}{c}
\mathcal{L}_i u + \mathcal{A}_i u \\ 
\mathcal{L}_i h%
\end{array}%
\right).
\end{equation*}%
The integrals in (\ref{shortsrsw}) are
of Stratonovitch type. 
In the SALT approach, the random parameter typically multiplies the gradient of the solution, but for the SRSW model it contains also a zero-order operator denoted here by $\mathcal{A}_i$. The amplitude of the noise is therefore modulated by the $\mathcal{G}_i=\mathcal{L}_i + \mathcal{A}_i$ operators. This is specifically designed to describe the (otherwise un-modelled) effect of the small-scale components on the large-scale components of the fluid, see \cite{Holm2015}. A data assimilation application for this (still uncalibrated) model has first been implemented in \cite{lucrisanlangmemin}. Numerical implementations and particle filter algorithms for other SALT models (2D Euler, SQG) have been developed in \cite{Wei1} and \cite{Weifinal}.

\section{Calibration of the noise correlation}\label{section:calibration}
%%%%%%%%%%%%%%%%%%%%%%%%%%%%%%%%%%%%%%%%%%%%%%%%%%
%%%%%%%%%%%%%%%%%%%%%%%%%%%%%%%%%%%%%%%%%%%%%%%%%%
%%%%%%%%%%%%%%%%%%%%%%%%%%%%%%%%%%%%%%%%%%%%%%%%%%
%\todo[inline]{Calibration Numerics (Alex)}

%\subsection{Generic Method}
%%%%%%%%%%%%%%%%%%%%%%%%%%%%%%%%%%%%%%%%%%%%%%%%%%
%%%%%%%%%%%%%%%%%%%%%%%%%%%%%%%%%%%%%%%%%%%%%%%%%%

\subsection{Case study: Rotating Shallow Water equation}
%%%%%%%%%%%%%%%%%%%%%%%%%%%%%%%%%%%%%%%%%%%%%%%%%%
%%%%%%%%%%%%%%%%%%%%%%%%%%%%%%%%%%%%%%%%%%%%%%%%%%

\begin{figure}

\centering
\begin{subfigure}{0.325\textwidth}
    \centering
    \includegraphics[width=\textwidth]{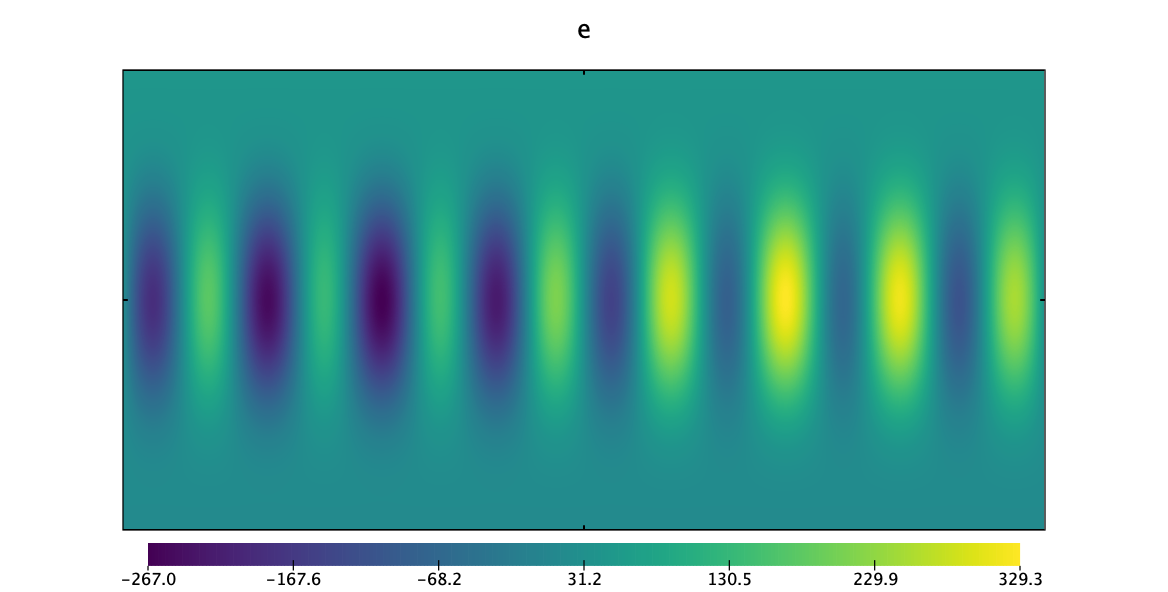}
    \caption{Surface Elevation $\eta$.}
    \label{fig:e_start}
\end{subfigure}
\begin{subfigure}{0.325\textwidth}
    \centering
    \includegraphics[width=\textwidth]{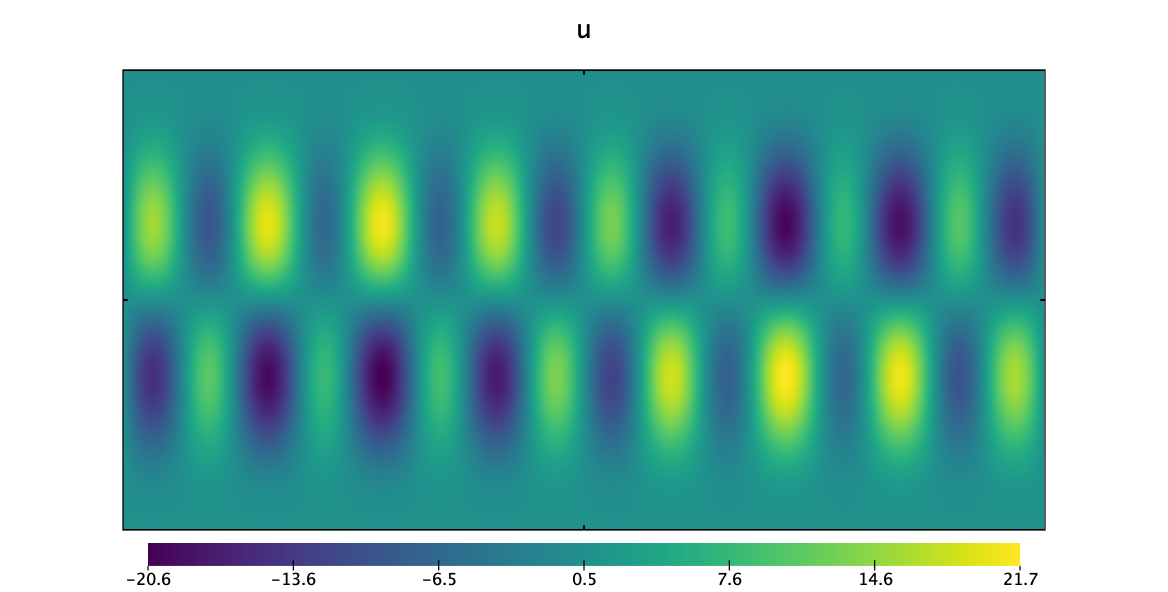}
    \caption{Zonal Velocity $u$.}
    \label{fig:u_start}
\end{subfigure}
\begin{subfigure}{0.325\textwidth}
    \centering
    \includegraphics[width=\textwidth]{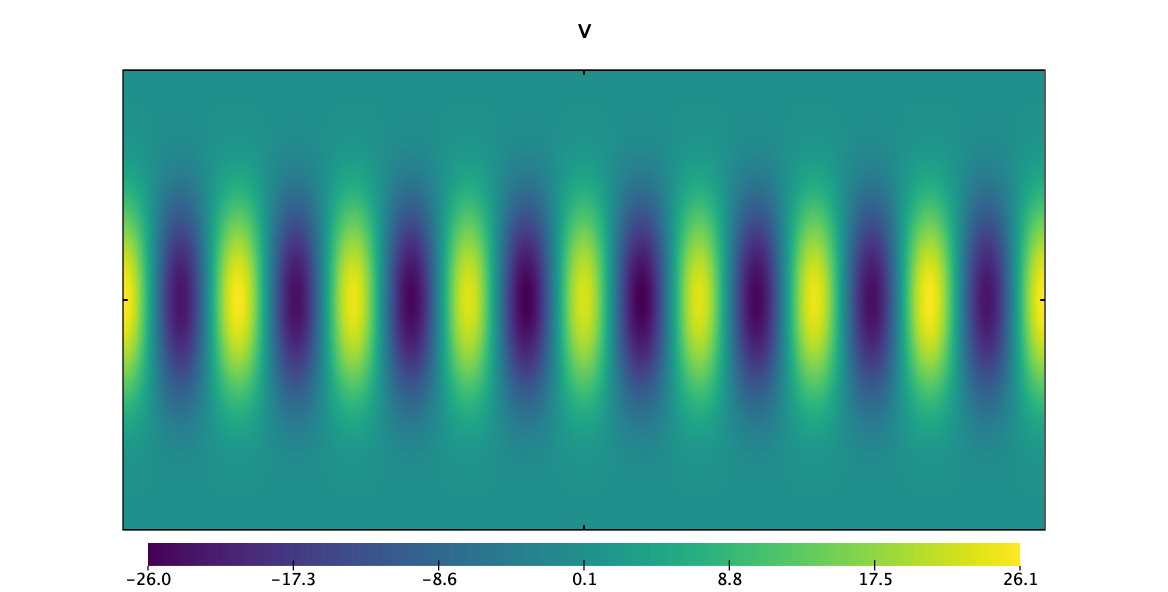}
    \caption{Meridional Velocity $v$.}
    \label{fig:v_start}
\end{subfigure}
\caption{Starting Condition \eqref{eq:start_cond} of the Rotating Shallow Water model from which we spin up the system. The fields are given on the computational domain $\Omega$ discretised using the fine grid of size $N_x = 2224$ and $N_y=320$.}
\label{fig:rsw_init}
 
 \end{figure}

\begin{figure}[ht!]

\centering
\begin{subfigure}{0.325\textwidth}
    \centering
    \includegraphics[width=\textwidth]{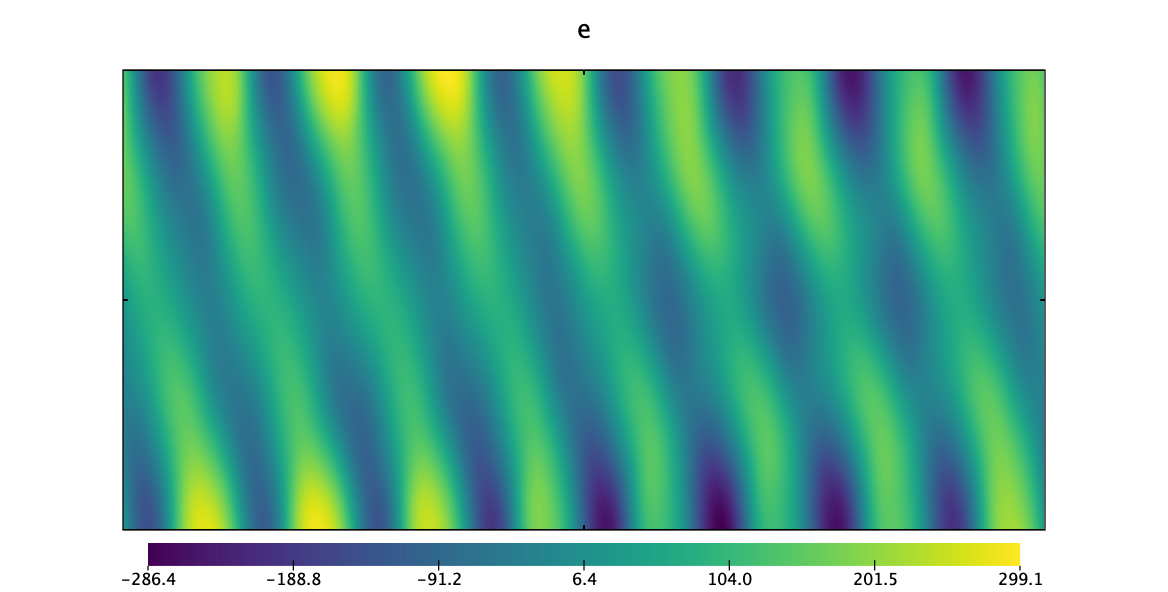}
    \caption{Elevation $\eta$ -- Fine Grid.}
    \label{fig:init_fine_e}
\end{subfigure}
\begin{subfigure}{0.325\textwidth}
    \centering
    \includegraphics[width=\textwidth]{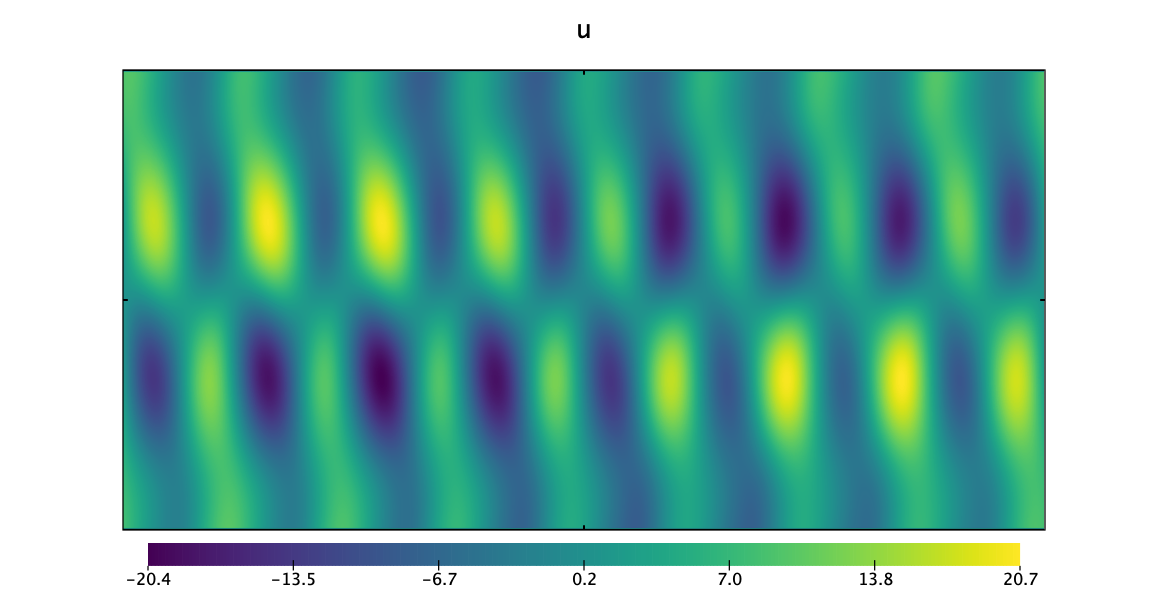}
    \caption{Zonal Vel.~$u$ -- Fine Grid.}
    \label{fig:init_fine_u}
\end{subfigure}
\begin{subfigure}{0.325\textwidth}
    \centering
    \includegraphics[width=\textwidth]{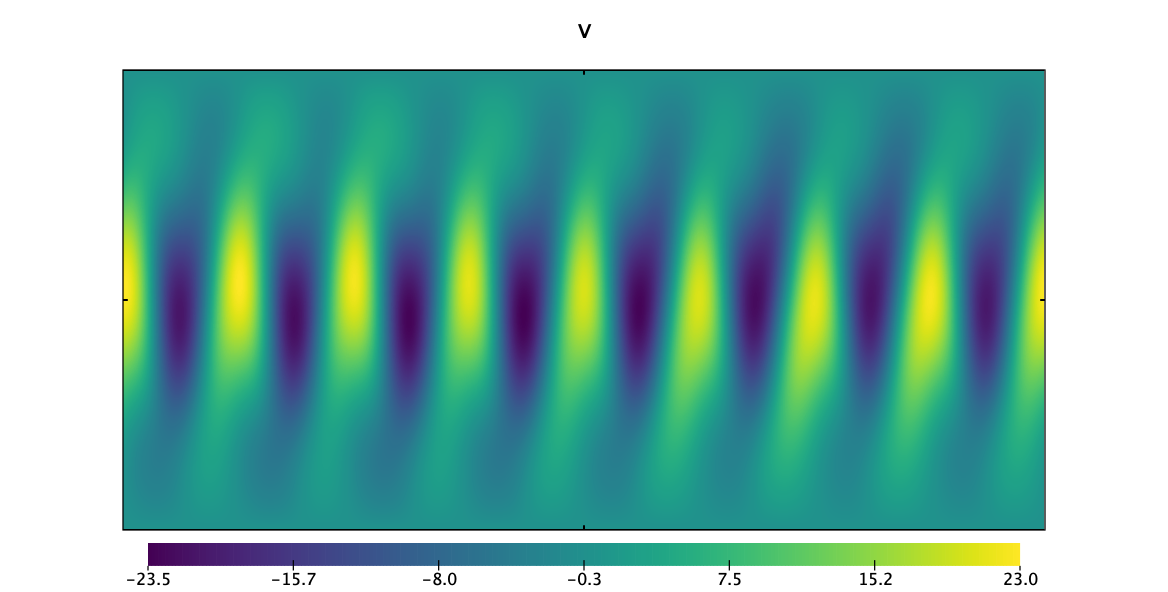}
    \caption{Meridional Vel.~$v$ -- Fine Grid.}
    \label{fig:init_fine_v}
\end{subfigure}
\hfill
\begin{subfigure}{0.325\textwidth}
    \centering
    \includegraphics[width=\textwidth]{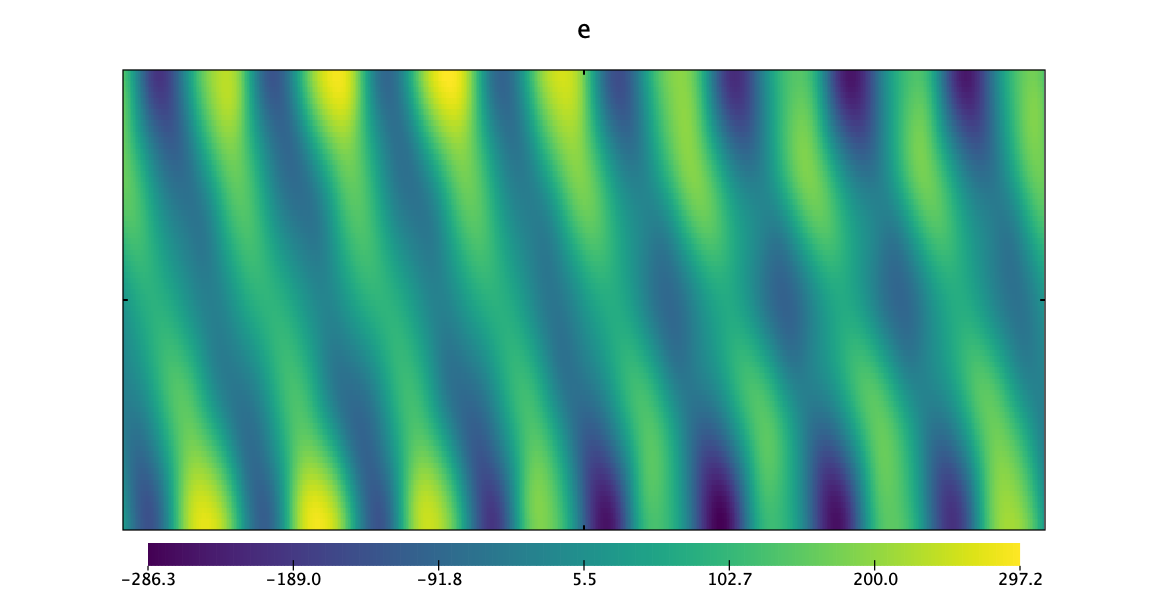}
    \caption{Elevation $\eta$ -- Coarsening $c=4$.}
    \label{fig:init_coarse4_e}
\end{subfigure}
\begin{subfigure}{0.325\textwidth}
    \centering
    \includegraphics[width=\textwidth]{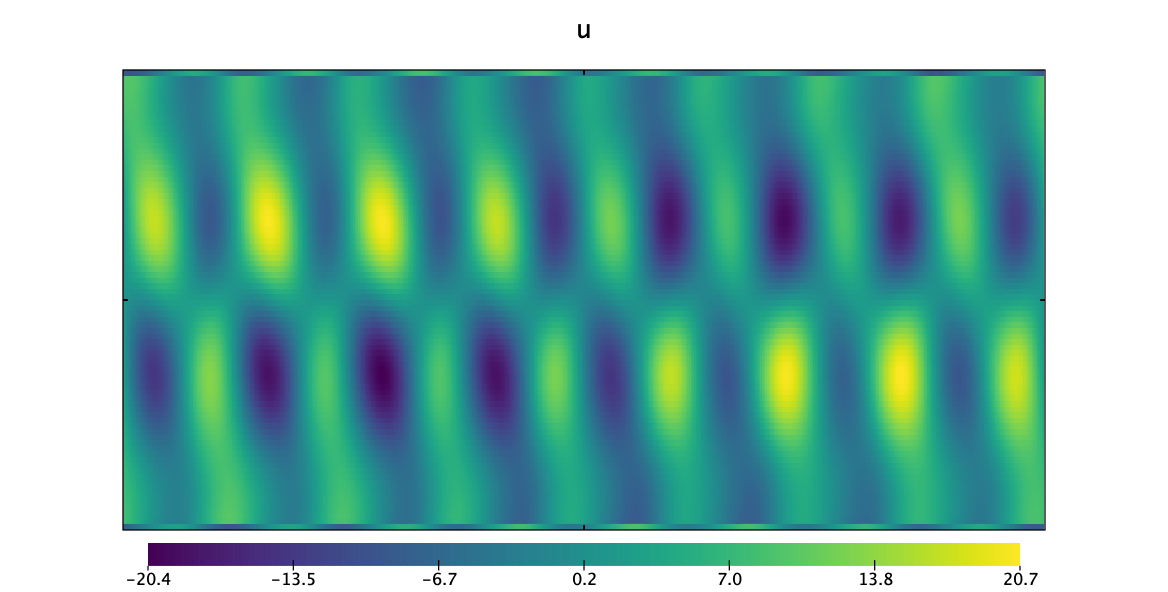}
    \caption{Zonal Vel.~$u$ -- Coarsening $c=4$.}
    \label{fig:init_coarse4_u}
\end{subfigure}
\begin{subfigure}{0.325\textwidth}
    \centering
    \includegraphics[width=\textwidth]{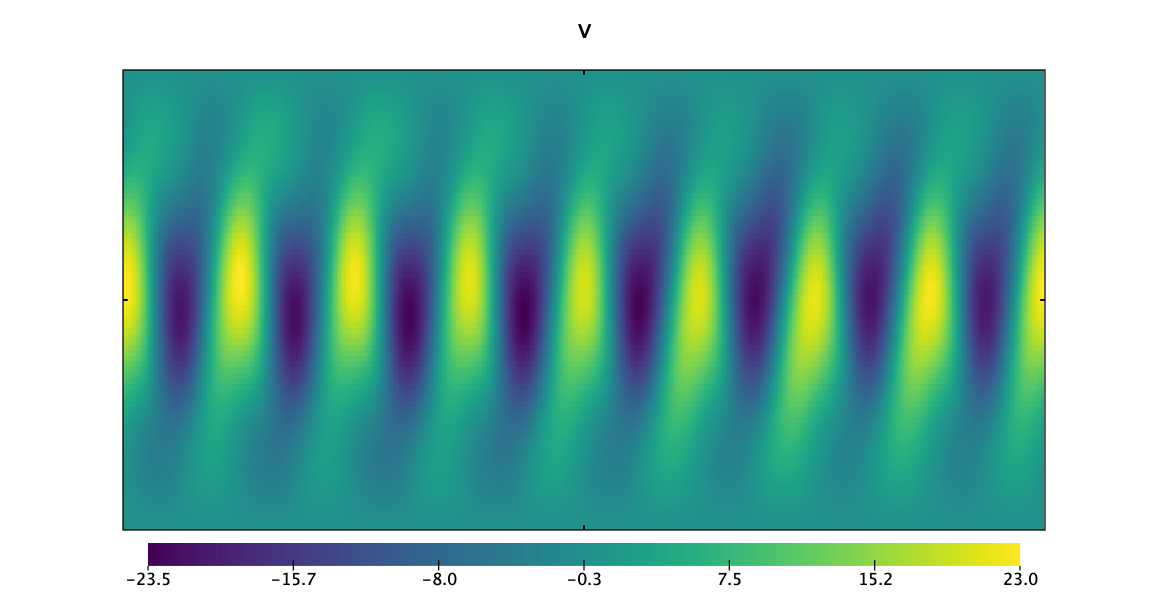}
    \caption{Meridional Vel.~$v$ -- Coarsening $c=4$.}
    \label{fig:init_coarse4_v}
\end{subfigure}
\hfill
\begin{subfigure}{0.325\textwidth}
    \centering
    \includegraphics[width=\textwidth]{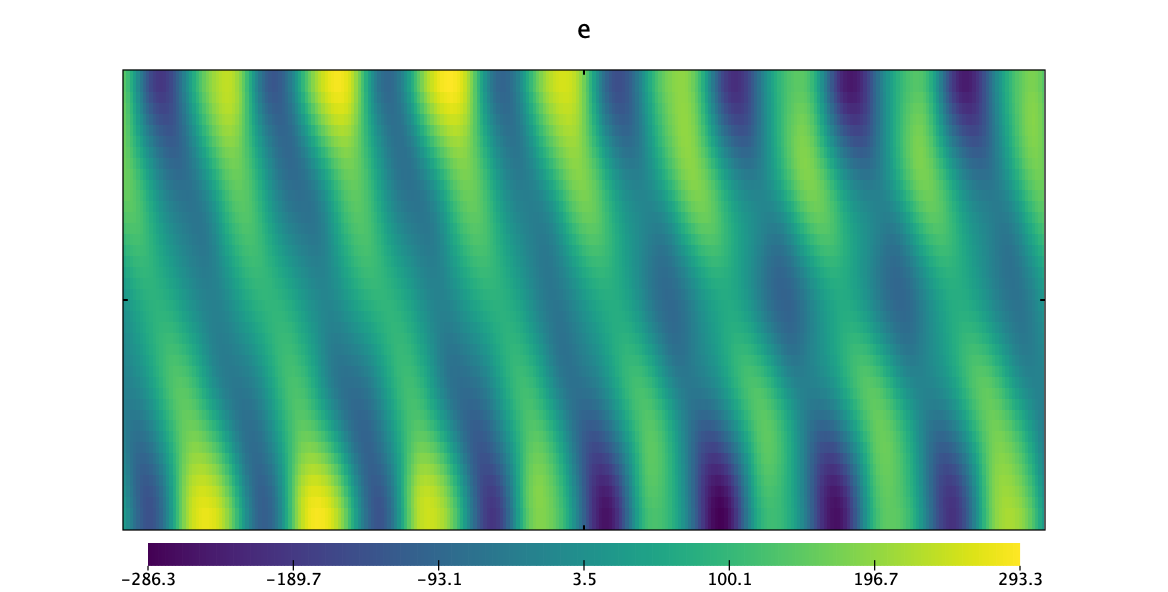}
    \caption{Elevation $\eta$ -- Coarsening $c=8$.}
    \label{fig:init_coarse8_e}
\end{subfigure}
\begin{subfigure}{0.325\textwidth}
    \centering
    \includegraphics[width=\textwidth]{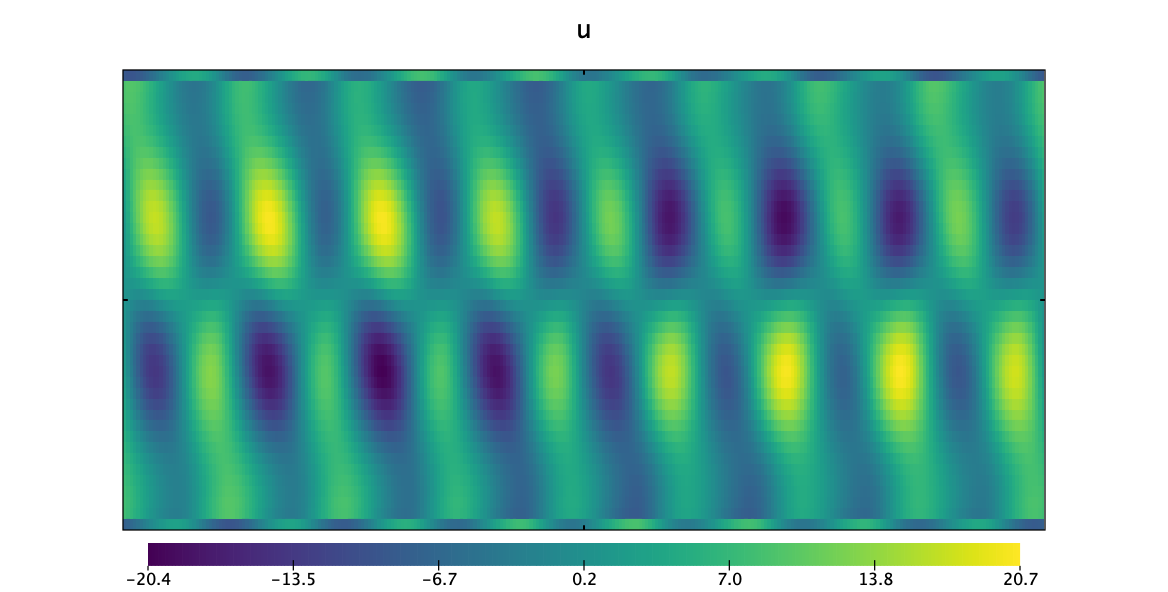}
    \caption{Zonal Vel.~$u$ -- Coarsening $c=8$.}
    \label{fig:init_coarse8_u}
\end{subfigure}
\begin{subfigure}{0.325\textwidth}
    \centering
    \includegraphics[width=\textwidth]{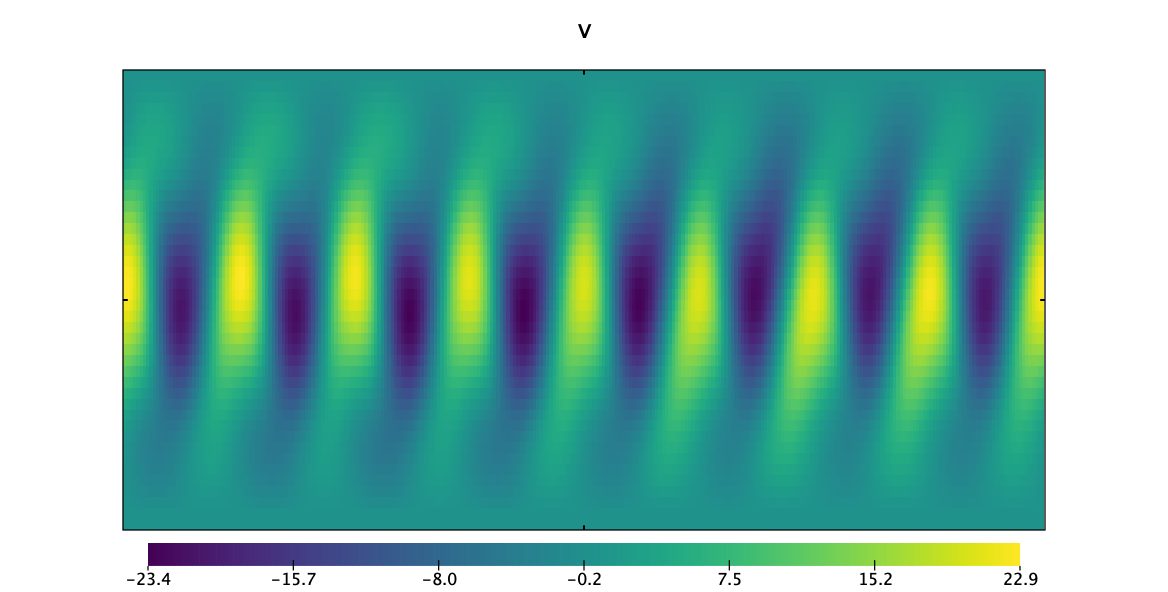}
    \caption{Meridional Vel.~$v$ -- Coarsening $c=8$.}
    \label{fig:init_coarse8_v}
\end{subfigure}
\caption{Initial condition obtained from the starting condition after the burn-in period. The top row shows the original, fine grid, solution of the Rotating Shallow Water PDE. The middle and bottom rows show the same solution projected on the coarsened grids with coarsening $c=4$ and $c=8$, respectively.} 
\label{fig:init}

\end{figure}

We are calibrating the Rotating Shallow Water System on the domain $\Omega=[0,L_x]\times[0,L_y]$ with $L_x = \SI{27787.5}{\km}$ and $L_y = \SI{3975}{\km}$. As starting condition we have chosen the function
\begin{equation}
    \eta_0(x,y) = -a\arctan{\left(0.05\left(\frac{y}{L_y}-0.5\right)\pi\right)} + \Bigg[ a \sin{\left(16\pi\frac{x}{L_x}\right)} + 0.5 a \sin{\left(2\pi\frac{x}{L_x}\right)}\Bigg] \sin{\left(\pi\frac{y}{L_y}\right)}^4,
\label{eq:start_cond}
\end{equation}
where the parameter $a$ was chosen to be equal to $100$. From this initial state for the elevation, the starting velocities were computed using geostrophic balance, which is an approximate balance in the system, equating the Coriolis force with the pressure gradient.
The starting condition is depicted in Figure~\ref{fig:rsw_init}.
We use a burn-in period of $1000$ timesteps of size $\delta t=\SI{22.5}{\s}$ in order to construct the initial condition of the model state.  This will ensure a more realistic initial state of the shallow water model. The initial condition obtained after the burn-in period is shown in Figure~\ref{fig:init}.


\paragraph{Mollification of the fine grid data.}
We use the PDE for $h$ run on the fine grid to obtain $h^f$ and construct $C(h^f)$ by mollifying $h^f$.
%%%%%%%%%%%%%%%%%%%%%%%%%%%%%%%%%%%%%%%%%%%%%%%%%%
For the mollification, we chose to apply a low-pass filter. Consider the field $h^f = h^f_0$ obtained as above. Write $h^f_{ij} = h^f(x_i, y_j)$ for the values on the grid points. Then consider the normalised convolution kernel $K\in\mathbb{R}^{k_x\times k_{y}}$ which can be adapted to the desired degree of coarsening.  We obtain the \emph{coarsened} field $C(h^f)$ from $h^f$ via the discrete convolution for the interior $i=1+\lfloor k_x/2 \rfloor,\dots, N_x-\lfloor k_x/2 \rfloor$, $j=1+\lfloor k_y/2 \rfloor,\dots, N_y-\lfloor k_y/2 \rfloor$ as
  \begin{equation}
      C(h^f)_{ij} = (h^f \star_{\text{d}} K)_{ij} = \sum_{\iota=0}^{k_x-1}\sum_{\jmath=0}^{k_y-1} K_{\iota+1,\jmath+1} h^f_{i-\lfloor k_x/2 \rfloor+\iota, j - \lfloor k_y/2 \rfloor+\jmath}
  \end{equation}
  and on the East and West boundaries $i\in \{1,\dots, \lfloor k_x/2 \rfloor\} \cup \{N_x-\lfloor k_x/2 \rfloor +1,\dots, N_x\}$ we set periodic boundary conditions
  \begin{equation}
       C(h^f)_{ij} = (h^f \star_{\text{d}} K)_{ij} = \sum_{\iota=0}^{k_x-1}\sum_{\jmath=0}^{k_y-1} K_{\iota+1,\jmath+1} h^f_{i-\lfloor k_x/2 \rfloor+\iota, j - \lfloor k_y/2 \rfloor+\jmath}; \; j \in \{\lfloor k_y/2 \rfloor , \dots, N_y-\lfloor k_y/2 \rfloor \}
  \end{equation}
  $i\in \{1,\dots, \lfloor k_x/2 \rfloor\} \cup \{N_x-\lfloor k_x/2 \rfloor +1,\dots, N_x\}$, $j\in \{1,\dots, \lfloor k_y/2 \rfloor\} \cup \{N_y-\lfloor k_y/2 \rfloor +1,\dots, N_y\}$ whereas on the North and South Boundaries we simply set
  \begin{equation}
      C(h^f)_{ij} = h^f_{ij}
  \end{equation}
  so that we conserve the PDE boundary conditions as well as the integral
  \begin{equation}
      \sum_{ij}  C(h^f)_{ij} = \sum_{ij} h^f_{ij}.
  \end{equation}
  As mentioned above, this procedure is known as a low-pass filter in signal processing. In our context, the low-pass filter is a principled choice, since the high frequency effects present in the field will be smoothed out which corresponds to the fact that we expect coarser simulations to misrepresent higher frequency effects as well. Alternative coarsenings may be performed using Gaussian mollifiers or subsampling. 
  
\paragraph{Data generation.} Here we perform steps 1. and 2. as described in  
Section \ref{cals}. More precisely, we prepare the data required to estimate 
the noise parametrization.  We use the time-increments of the discrepancy between the high resolution fields and coarsened fields with a calibration timestep $\delta t_\text{calib}\ll \delta t^f_{\text{PDE}}$, which corresponds to an Euler-Maruyama step of the stochastic integral. Thus we consider the data over a calibration time-grid $t^n_{\text{calib}}$ obtained as a sub-grid of the PDE time grid as described below.
 \begin{equation}
     \Delta_{t^n_{\text{calib}}} = C(h^f_{t^n_{\text{calib}}+\delta t_\text{calib}})-h^f_{t^n_{\text{calib}}+\delta t_\text{calib}} - (C(h^f_{t^n_{\text{calib}}})-h^f_{t^n_{\text{calib}}})
 \end{equation}
 \vspace{3mm}

We denote the calibration partition by  $t_0<t_1<t_2<...<t_n$.
 The partition has to be sufficiently sparse so that data is decorrelated. 
We choose the calibration time-grid $t^n_{\text{calib}}$ by estimating a \emph{decorrelation time} of the data on the full PDE time grid $\Delta_{t^f_i}$.
This is done as follows. We compute the auto-correlation of the data time-series $\Delta_{t^f_i}$ pointwise in space for different lags $\ell$ as
    \begin{equation}
        \mathcal{C}(\ell) = \frac{\frac{1}{N_{steps}-\ell}\sum_{i=0}^{N_{steps}-\ell} \left[\Delta_{t^f_i}\odot\Delta_{t^f_{i}+\ell\delta t^f_{\text{PDE}}} - \Delta_{t^f_i}^{\odot 2}\right]}{\frac{1}{N_{steps}-\ell-1}\sum_{i=0}^{N_{steps}-\ell}\left[ \Delta_{t^f_i} - \frac{1}{N_{steps}-\ell}\sum_{i=0}^{N_{steps}-\ell} \Delta_{t^f_i}\right]^{\odot2}} \in\mathbb{R}^{N_x\times N_y},
    \end{equation}
where $\odot$ denotes the element-wise product of matrices and the division is also taken to be elementwise. 
    We choose as the decorrelation time the lag $\ell_{\text{decorr}}$ as the minimum value of $\ell$ for which the mean absolute correlation in space
    \begin{equation}
        \bar{\mathcal{C}}(\ell) = \frac{1}{N_x N_y}\sum_{ij} \left|\mathcal{C}_{ij}(\ell)\right|
    \end{equation}
    stabilizes to a value $\bar{\mathcal{C}}(\ell_{\text{decorr}})\sim \alpha_{\text{decorr}}$ where we can choose the decorrelation threshold $\alpha_{\text{decorr}}$. The calibration time grid is then 
    \begin{equation}
        t^n_{\text{calib}} = t^f_0 + n \ell_{\text{decorr}} \delta t^f_{\text{PDE}}, \qquad n=1,\dots, \lfloor N/\ell_\text{decorr}\rfloor.
    \end{equation}

%\pj{(What is n?) Is that from the time partitioning? I seem to miss something.}\todo[inline]{Yes, n counts the timepoints. Need to clarify.  
Next, we prepare the input data.
%%%%%%%%%%%%%%%%%%%%%%%%%%%%%%%%%%%%%%%%%%%%%%%%%%
Let 
$\hat{h}_i:=C(h^f_{t_i})-h^f_{t_i}$. Choose $\delta<\min_i{(t_i-t_{i-1})}$ and compute 
$\delta \hat{h}_i:=\hat{h}_{t_i+\delta}-\hat{h}_{t_i}$. 
The input data for the calibration is then $\{\Delta_{t^n_{\text{calib}}}; n= 1,\dots\}$ obtained from the above method.

\paragraph{Solving the calibration equation.} Here we perform step 3. as described in  Section \ref{cals}.
%%%%%%%%%%%%%%%%%%%%%%%%%%%%%%%%%%%%%%%%%%%%%%%%%%
Given the data generated as above, we solve the sequence of hyperbolic equations  
 \begin{equation}\label{calibc}
 \delta \hat{h}_i = C(h_x^f)\frac{\partial \psi_i}{\partial y}-C(h_y^f)\frac{\partial \psi_i}{\partial x}.
 \end{equation}
for each calibration time, with periodic boundary condition in the East-West direction and in the North-South direction we impose the free-slip constraint
 \[
 \frac{\partial \psi_i}{\partial x} = 0 
 \]
In order to implement the solution of this equation, we use a first order discontinuous Galerkin finite element method, accessible from the Firedrake\footnote{\href{https://firedrakeproject.org/}{https://firedrakeproject.org/}} package.

To this end we rewrite the equations as follows
 \begin{equation}
     f = \vect{q} \cdot \nabla \psi
 \end{equation}
 with $f:=\partial \hat{h}_i$ and $\vect{q}:=-\nabla^\perp C(h^f)$
 %\todo[inline]{no bar}. 
 The boundary is periodic in East-West, and in North and South we set $\psi = 0$. 
 
 The weak form on each cell is
 \begin{equation}
     \int_e f \phi_e \,dx = \int_e \phi_e\vect{q}\cdot\nabla\psi \,dx
 \end{equation}
 so that after integration by parts
 \begin{equation}
     \int_e f \phi_e \,dx  = -\int_e \psi \nabla\cdot(\phi_e\vect{q}) \, dx + \int_{\partial e} \phi_e\psi\vect{q}\cdot\vect{n}_e\, dS
 \end{equation}
 There are three different types of boundary facets: Interior, Exterior inflow, and Exterior outflow, where the latter two are determined by the velocity vector field $\vect{q}$. The boundary condition of the PDE applies for the exterior inflow facets.
Summing up the elementwise weak form over the domain and decomposing the boundary according to the types of facets we get
 \begin{align}
     \int_\Omega f \phi \,dx  &= -\int_\Omega \psi \nabla\cdot(\phi_e\vect{q}) \, dx + \int_{\Gamma_{\text{int}}} \tilde{\psi}(\phi_+\vect{q}\cdot\vect{n}_+ + \phi_-\vect{q}\cdot\vect{n}_-)\, dS\\ &+ \int_{\Gamma_{\text{ext, in}}} \phi\psi_0\vect{q}\cdot\vect{n}\, dS +  \int_{\Gamma_{\text{ext, out}}}  \phi\psi\vect{q}\cdot\vect{n}\, dS
 \end{align}
%\pj{(I'm not sure I understand where his is going.)}\todo[inline]{Not necessarily going anywhere, this is meant to be a clarification of the how the calibration equation is implemented in terms of its weak form using FEM. This may be better placed in the Appendix. I don't mind either way.}
 
 \paragraph{Extract a basis for the stochastic noise.} This is Step 4. as described in Section \ref{cals}.
 
 %%%%%%%%%%%%%%%%%%%%%%%%%%%%%%%%%%%%%%%%%%%%%%%%%%
 The solutions the sequence of calibration equations can be thought of as \emph{stream functions} for the perturbation fields $\tilde{u}$, $\tilde{v}$, so that $(\tilde{u},\tilde{v})^\top = \grad^\perp \psi$. This grid based data $(\tilde{u}_{ij},\tilde{v}_{ij})$ is vectorized as $\Psi = ((\tilde{u}_{ij}), (\tilde{v}_{ij}))^\top$ and represented in the following form 
\begin{equation*}
\frac{\Psi _{i}-\bar{\Psi}}{\sqrt{\delta }}=\sum_{j=1}^{N}\zeta ^{j}
\Delta W_{i}^{j} 
\end{equation*}%
where $\Delta W_{i}^{j}$ are i.i.d. standard normal random variables $\Delta
W_{i}^{j}$ $\sim N(0,1)$. We estimate $\zeta ^{j}$ by minimising 
\begin{equation*}
\mathbb{E}\left[ \left\vert \left\vert \sum_{i}\left( \frac{\Psi _{i}-\bar{\Psi}}{%
\sqrt{\delta }}-\sum_{j=1}^{N}\zeta ^{j}\left( \Delta W_{i}^{j}\right)
\right) \right\vert \right\vert ^{2}\right] ,
\end{equation*}%
where the choice of $N$ can be decided
by using empirical orthogonal
functions (EOFs). The EOFs can be thought of as principal components that
correspond to the spatial correlations of a field. To compute the EOFs we use the Principal Component Analysis algorithm based on the singular value decomposition (SVD).
% \begin{itemize}
% \item Write the
% data time series $\frac{\psi _{i}-\bar{\psi}}{\sqrt{\delta }}$ as a matrix $%
% \mathfrak{F}$ that will have on the $i$'th row the "serialised" data $\frac{%
% \psi _{i}-\bar{\psi}}{\sqrt{\delta }}$. That means that we need to transform
% the data $\frac{\psi _{i}-\bar{\psi}}{\sqrt{\delta }}$ that is initially a
% matrix corresponding to the grid points into a row vector. \pj{(But $%
% \mathfrak{F}$ is a matrix. Where is the row vector?)}
% \item Estimate
% the
% spatial covariance tensor by computing $\mathfrak{R}:=$ $\frac{1}{n-1}%
% \mathfrak{F}^{T}\mathfrak{F}.$
% \item Take the EOFs to be the eigenvectors of $
% \mathfrak{R}$, ranked in descending order according to the eigenvalues. \pj{Is is often faster to calculate an SVD of 
% $\mathfrak{F}$, why calculate the marix product first?)}\todo[inline]{This is true, we are using SVD (this is what is implemented in sklearn). This is an old description. We need to rewrite it according to what actually happens, i.e. stacked flattened $\tilde{u}$ and $\tilde{v}$ into a large vector. We get as many of them as we have data points from the calibration equation. So the data matrix has size $\sim 2 \cdot N_x\cdot N_y \times N_{data}$}
% \item Stop when we reach $N_\xi$\% of the total variability for $\frac{\psi _{i}-\bar{%
% \psi}}{\sqrt{\delta }}$. 

% \end{itemize}

Once the EOFs and the corresponding values are computed, the  procedure is complete. We are now in a position to assess the reliability of our calibration methodology. This is done in the next section.   

\section{Numerical Results}\label{nr}

In order to assess the reliability of our calibration methodology we conducted several numerical studies.
We present below the results for eight different scenarios summarized in Table~\ref{tab:scenarios}. In particular we analyze two different grid coarsenings, $c=4$ and $c=8$, different numbers of particles $N_p=50$ and $N_p=100$ as well as two variance thresholds $N_\xi=0.9$ and $N_\xi=0.99$. 
% \pj{(Note that $N_\xi$ is defined as a percentage, not a fraction.)}\todo[inline]{I like to keep it lose on this occasion, we can read the $\%$ sign as a fraction rather than a unit so that both are equivalent. Unless people find this outrageous...}

\begin{table*}\centering
\ra{1.3}
\begin{tabular}{@{}rrrcrr@{}}\toprule
& \multicolumn{2}{c}{$N_{\xi} = 0.90$} & \phantom{abc} & \multicolumn{2}{c}{$N_{\xi} = 0.99$} \\
\cmidrule{2-3} \cmidrule{5-6}  & $N_p=50$ & $N_p=100$  && $N_p=50$ & $N_p=100$\\ \midrule
$c=4$ & \textbf{a} & \textbf{b} && \textbf{c} & \textbf{d} \\
$c=8$ & \textbf{e} & \textbf{f} && \textbf{g} & \textbf{h} \\
\bottomrule
\end{tabular}
\caption{Scenarios for the numerical experiments.}
\label{tab:scenarios}
\end{table*}

The grid coarsening $c=4$ uses a coarse grid for the stochastic simulation which uses every fourth grid point relative to the fine PDE grid. Similarly, the $c=8$ coarsening uses every eighth grid point relative to the PDE grid. Note that the low-pass filter used to obtain the calibration data is different between the two. The $c=4$ coarsening uses the convolution kernel $K = \frac{1}{9} \mathbf{1}_{3\times 3}$, and the $c=8$ coarsening uses the $9\times 9$ kernel
\begin{equation*}
    K = \frac{1}{165}\begin{bmatrix}
    1 & 1 & 1 & 1 & 1 & 1 & 1 & 1 & 1 \\
    1 & 2 & 2 & 2 & 2 & 2 & 2 & 2 & 1 \\
    1 & 2 & 3 & 3 & 3 & 3 & 3 & 2 & 1 \\
    1 & 2 & 3 & 4 & 4 & 4 & 3 & 2 & 1 \\
    1 & 2 & 3 & 4 & 5 & 4 & 3 & 2 & 1 \\
    1 & 2 & 3 & 4 & 4 & 4 & 3 & 2 & 1 \\
    1 & 2 & 3 & 3 & 3 & 3 & 3 & 2 & 1 \\
    1 & 2 & 2 & 2 & 2 & 2 & 2 & 2 & 1 \\
    1 & 1 & 1 & 1 & 1 & 1 & 1 & 1 & 1 \\
    \end{bmatrix},
\end{equation*}
where the notation $\mathbf{1}_{m\times n}$ denotes an $m\times n$-matrix filled with ones.

\begin{figure}[ht!]

\centering
\begin{subfigure}{0.45\textwidth}
    \centering
    \includegraphics[width=\textwidth]{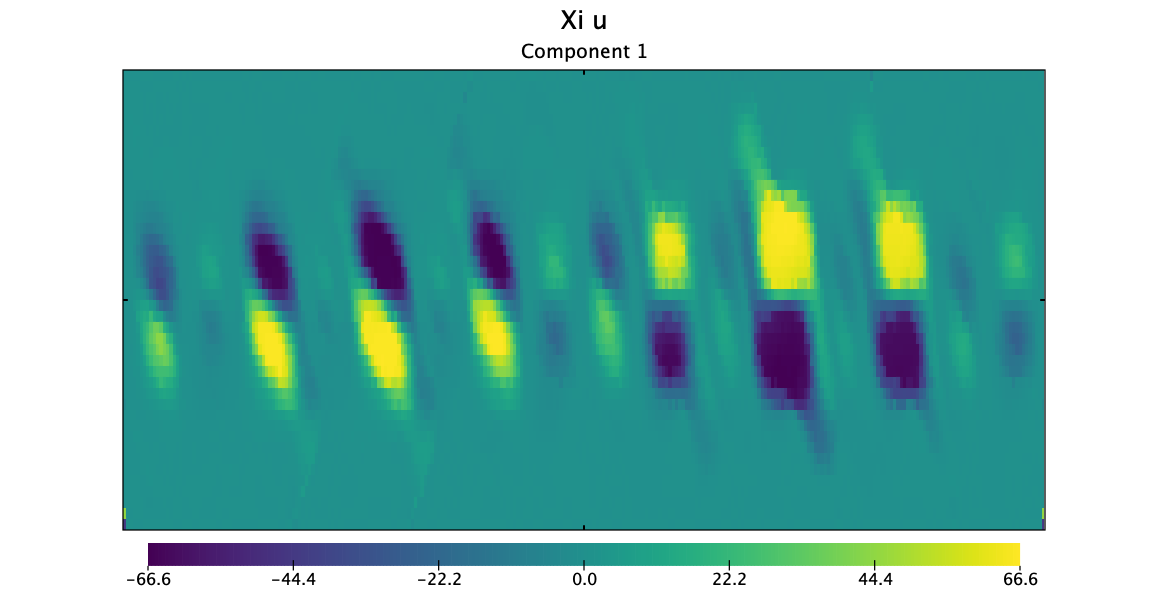}
    \caption{$\xi^u_1$.}
    \label{fig:xi_u_1}
\end{subfigure}
\begin{subfigure}{0.45\textwidth}
    \centering
    \includegraphics[width=\textwidth]{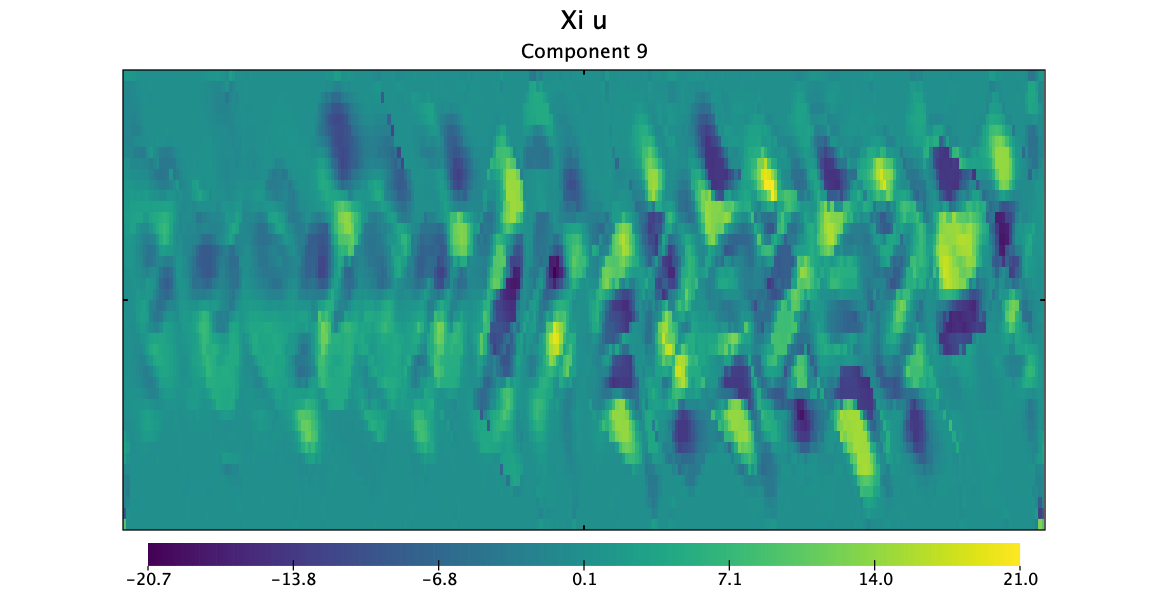}
    \caption{$\xi^u_9$.}
    \label{fig:xi_u_last}
\end{subfigure}
\hfill
\begin{subfigure}{0.45\textwidth}
    \centering
    \includegraphics[width=\textwidth]{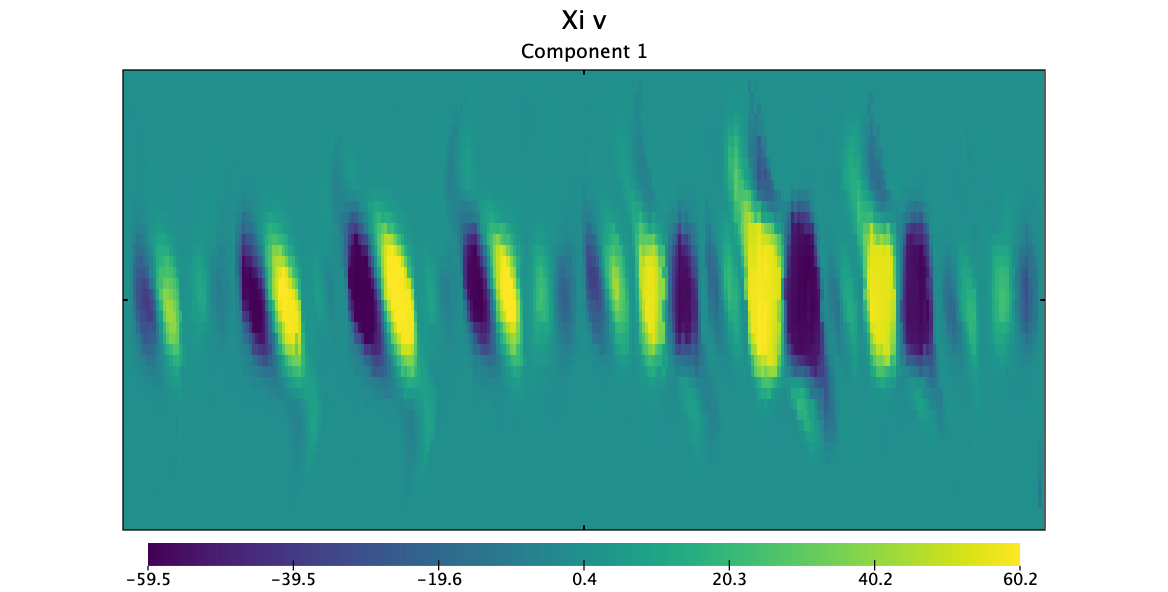}
    \caption{$\xi^v_1$.}
    \label{fig:xi_v_1}
\end{subfigure}
\begin{subfigure}{0.45\textwidth}
    \centering
    \includegraphics[width=\textwidth]{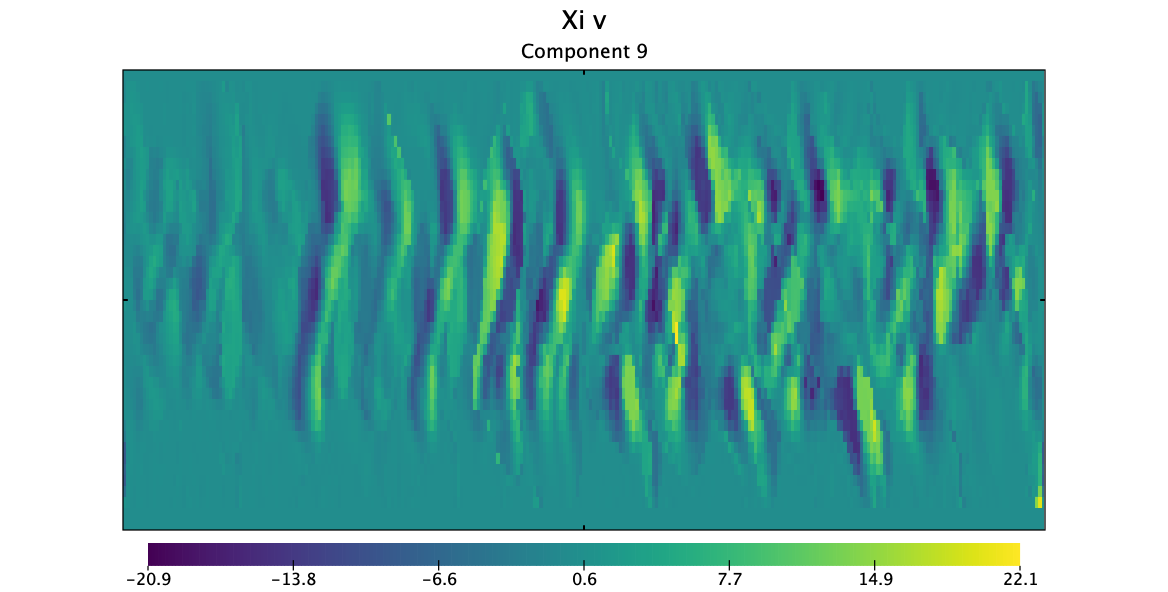}
    \caption{$\xi^v_9$.}
    \label{fig:xi_v_last}
\end{subfigure}
\caption{Estimated $\xi$'s on the $c=8$ grid. The left column depicts the first component and the right column the last component of the set explaining $99\%$ of the variance in the data obtained from the solutions of the calibration equation.} 
\label{fig:xis}

\end{figure}

A sample of the obtained EOFs ($\xi$'s) for the coarsening $c=8$ is depicted in Figure~\ref{fig:xis}. There are $9$ EOFs needed to explain $99\%$ of the variance in the data, and shown is the largest (left column) and the smallest (right column) of those EOFs. We observe that the magnitude for the higher order EOF is significantly smaller than that of the first EOF. Moreover, as expected, the first component exhibits larger scale structures and the last more fine grained structure.

% \begin{itemize}
%     \item We test the forecast reliability by comparing the forecast root mean square error (RMSE) 
%     {
%     \begin{equation*}
%         RMSE^2 = \frac{1}{N}\displaystyle\sum_{i=1}^N(x^i-x^{\star})^2
%     \end{equation*}
%     with the forecast ensemble spread (ES)
%     \begin{equation*}
%         ES^2=\frac{1}{N-1}\displaystyle\sum_{i=1}^N(x^i - \hat{E}x)^2
%     \end{equation*}
%     where
%     \begin{equation*}
%         \hat{E}x = \frac{1}{N}\displaystyle\sum_{i=1}^Nx^i
%     \end{equation*}
%     is the estimated ensemble mean and $x^{\star}$ is the truth.}
%     \item Using the adaptive tempering particle filter, the numerical results show that although the dimension of the observation is low compared to the degrees of freedom corresponding to the SPDE, there is sufficient information to allow for an accurate approximation of the truth.  
% \end{itemize}
% \begin{itemize}
% \item BIAS: ...
% \item $L^2 error:$
% \end{itemize}

In each case we plot the root mean square error (RMSE) and the ensemble spread (ES). These two are comparable, a feature highly appreciated by data assimilation practitioners as it shows that the estimate of the uncertainty as measured by the width of the ensemble is a good estimate of the actual error in the ensemble mean.

\begin{figure}[ht!]

\centering
\begin{subfigure}{0.32\textwidth}
    \centering
    \includegraphics[width=\textwidth]{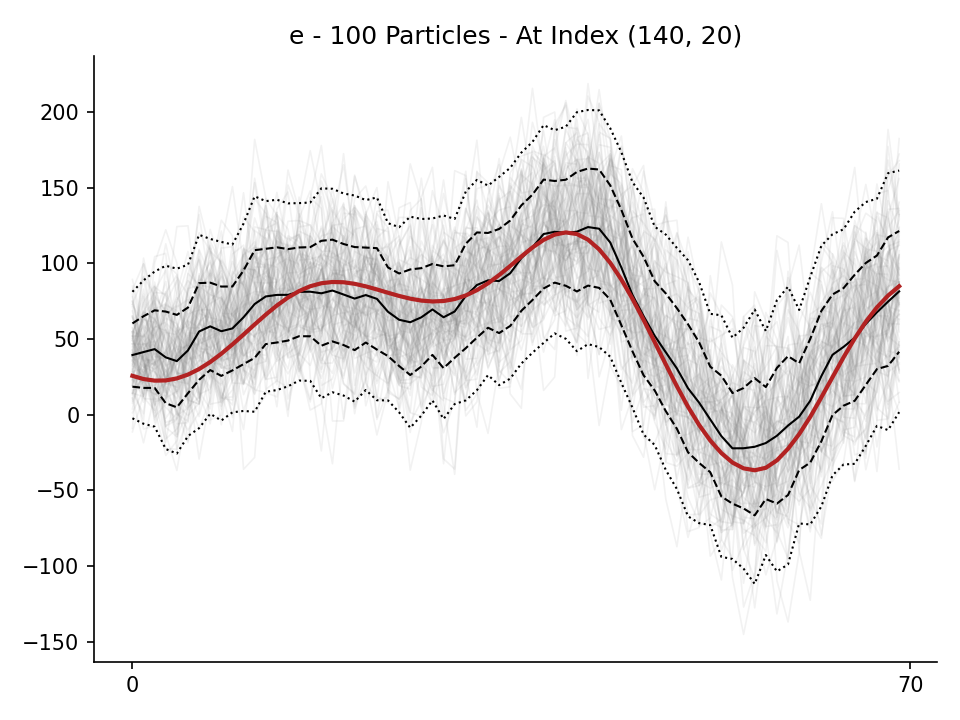}
    \caption{Elevation $\eta$ at $(140,20)$.}
    \label{fig:particle_e}
\end{subfigure}
\begin{subfigure}{0.32\textwidth}
    \centering
    \includegraphics[width=\textwidth]{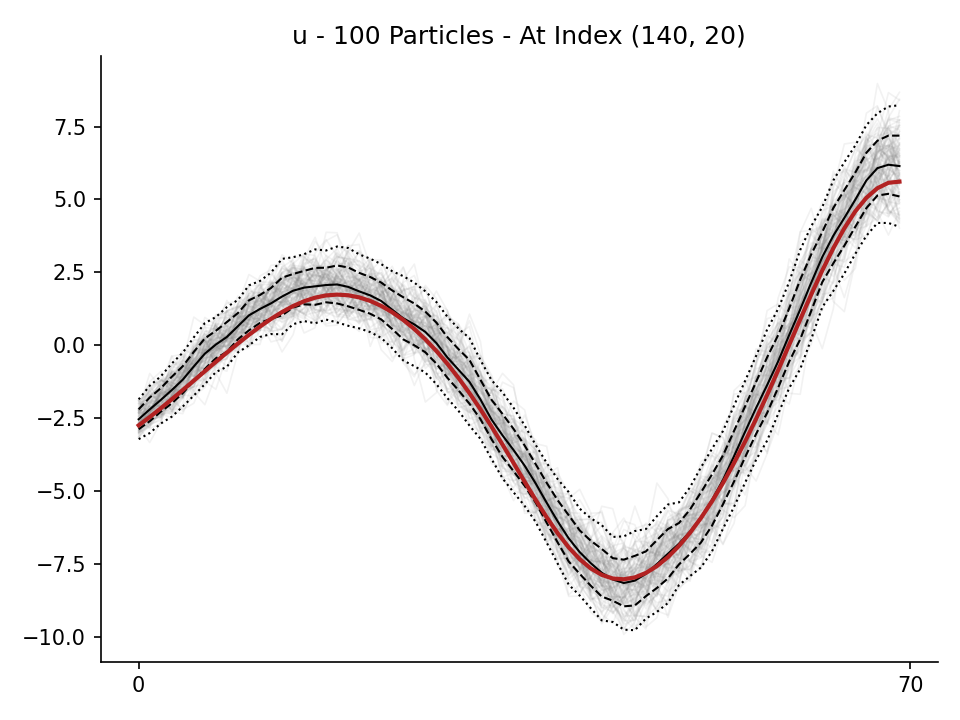}
    \caption{Zonal Vel.~$u$ at $(140,20)$.}
    \label{fig:particle_u}
\end{subfigure}
\begin{subfigure}{0.32\textwidth}
    \centering
    \includegraphics[width=\textwidth]{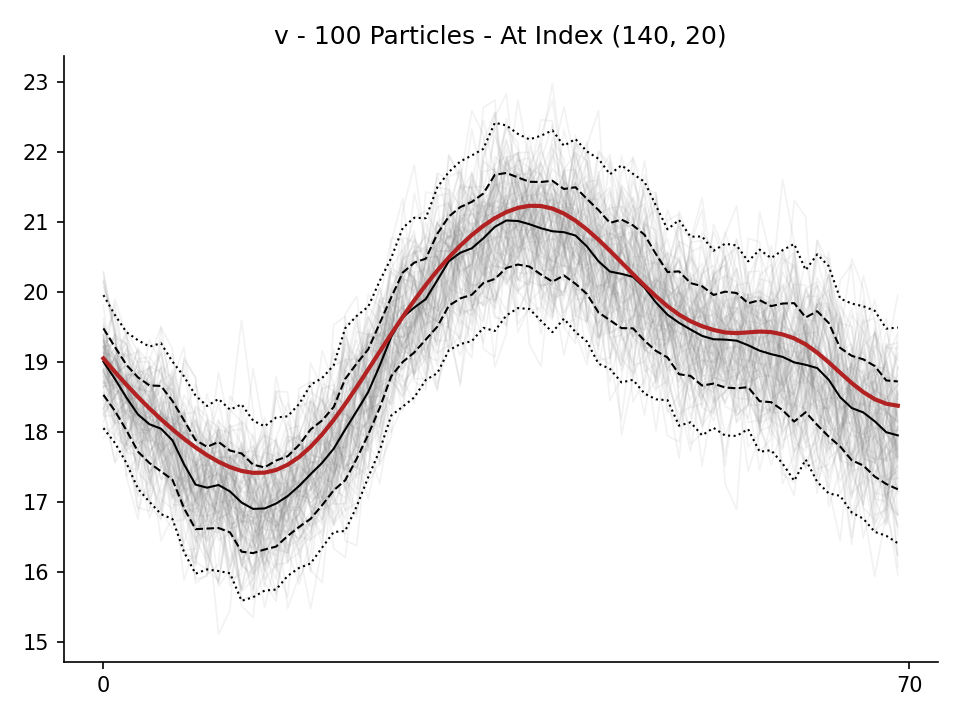}
    \caption{Meridional Vel.~$v$ at $(140,20)$.}
    \label{fig:particle_v}
\end{subfigure}
\hfill
\begin{subfigure}{0.32\textwidth}
    \centering
    \includegraphics[width=\textwidth]{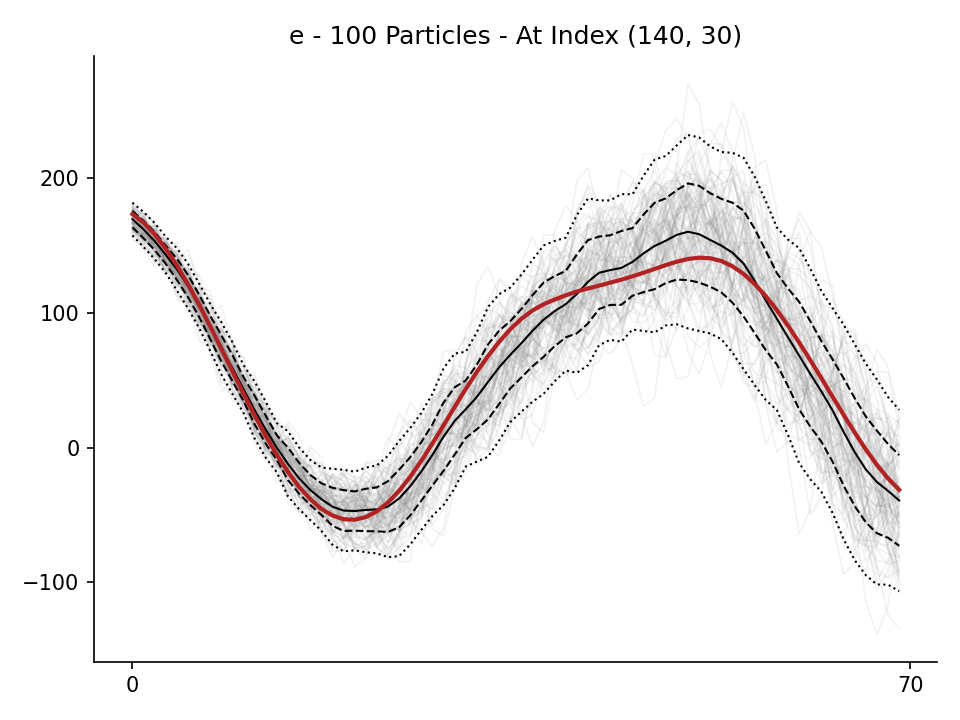}
    \caption{Elevation $\eta$ at $(140,30)$.}
    \label{fig:particle_e_2}
\end{subfigure}
\begin{subfigure}{0.32\textwidth}
    \centering
    \includegraphics[width=\textwidth]{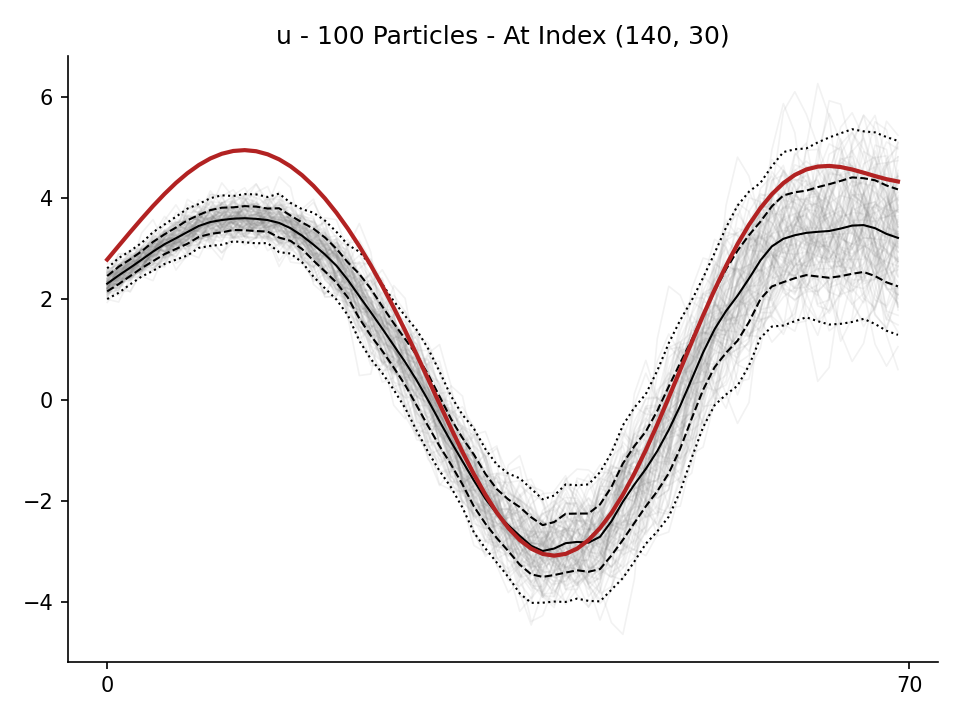}
    \caption{Zonal Vel.~$u$ at $(140,30)$.}
    \label{fig:particle_u_2}
\end{subfigure}
\begin{subfigure}{0.32\textwidth}
    \centering
    \includegraphics[width=\textwidth]{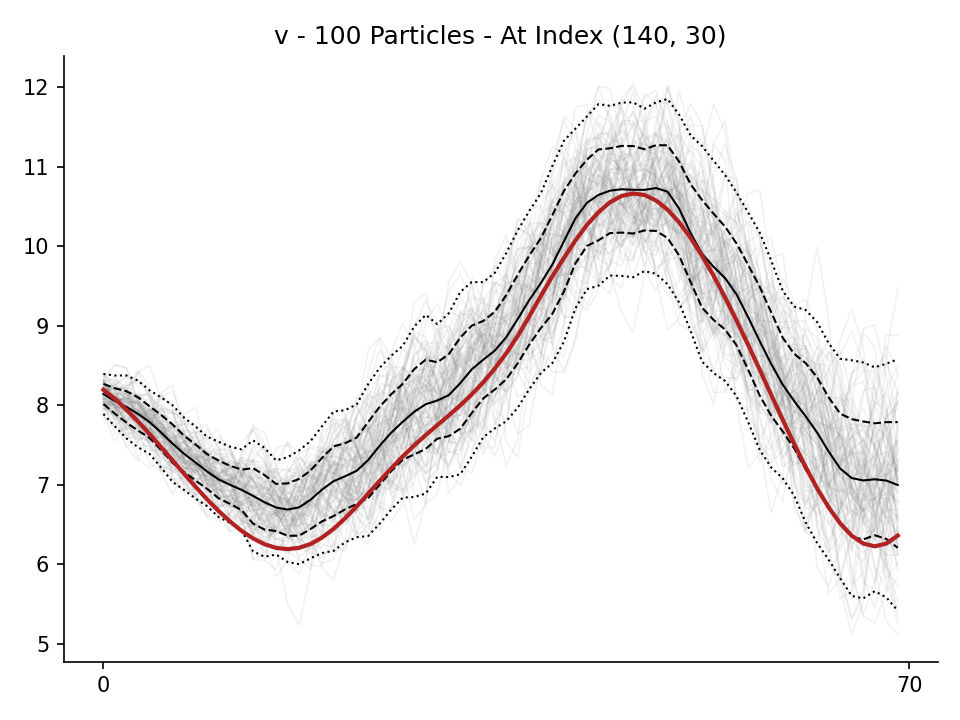}
    \caption{Meridional Vel.~$v$ at $(140,30)$.}
    \label{fig:particle_v_2}
\end{subfigure}
\hfill
\begin{subfigure}{0.32\textwidth}
    \centering
    \includegraphics[width=\textwidth]{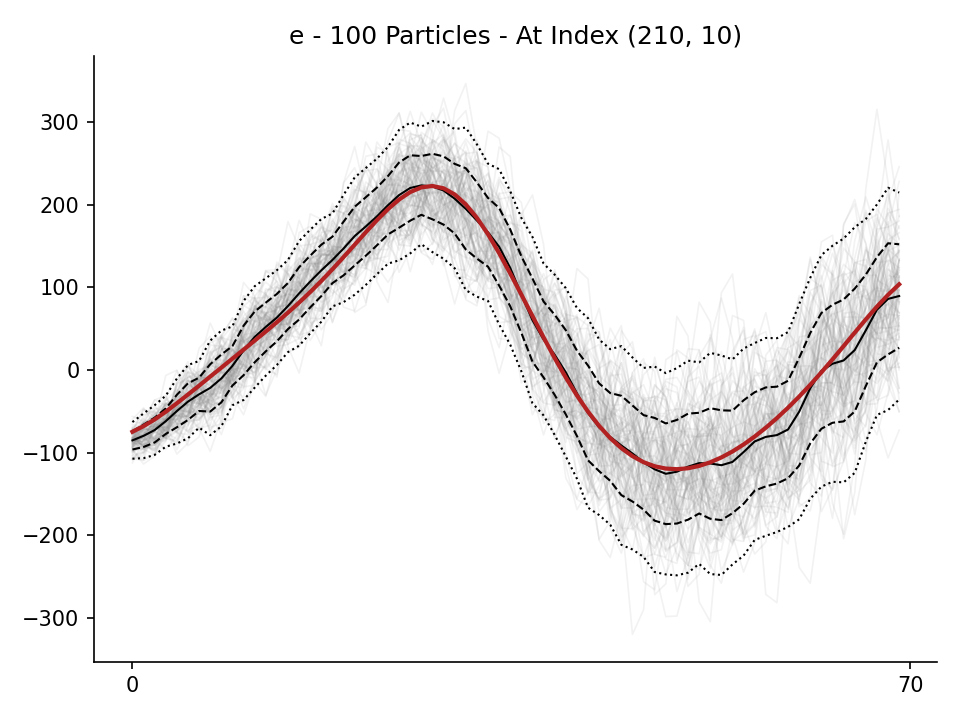}
    \caption{Elevation $\eta$ at $(210,10)$.}
    \label{fig:particle_e_3}
\end{subfigure}
\begin{subfigure}{0.32\textwidth}
    \centering
    \includegraphics[width=\textwidth]{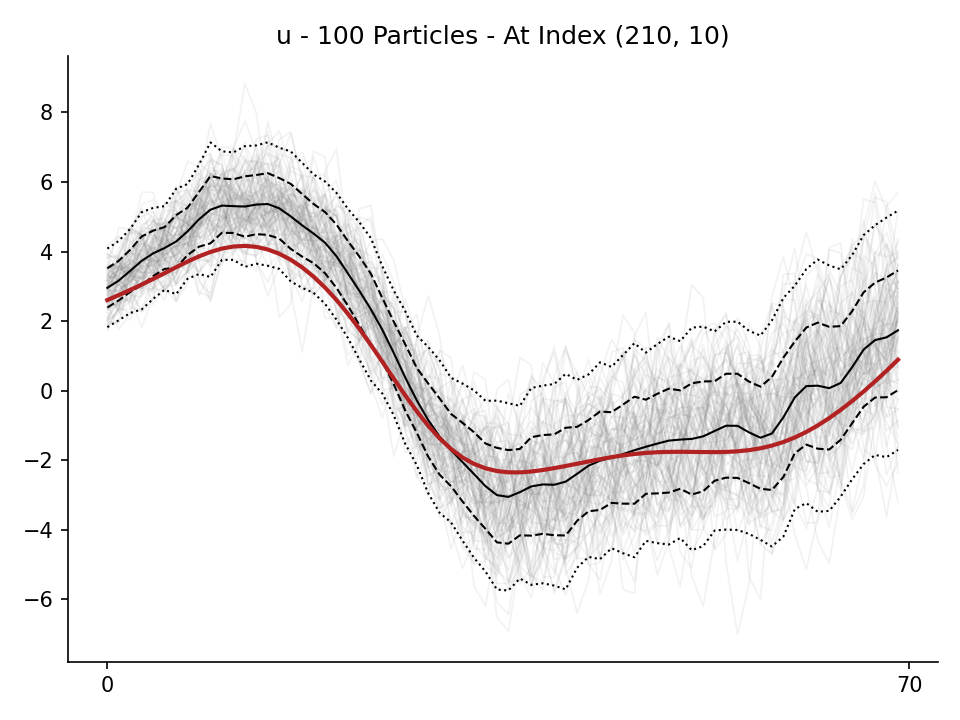}
    \caption{Zonal Vel.~$u$ at $(210,10)$.}
    \label{fig:particle_u_3}
\end{subfigure}
\begin{subfigure}{0.32\textwidth}
    \centering
    \includegraphics[width=\textwidth]{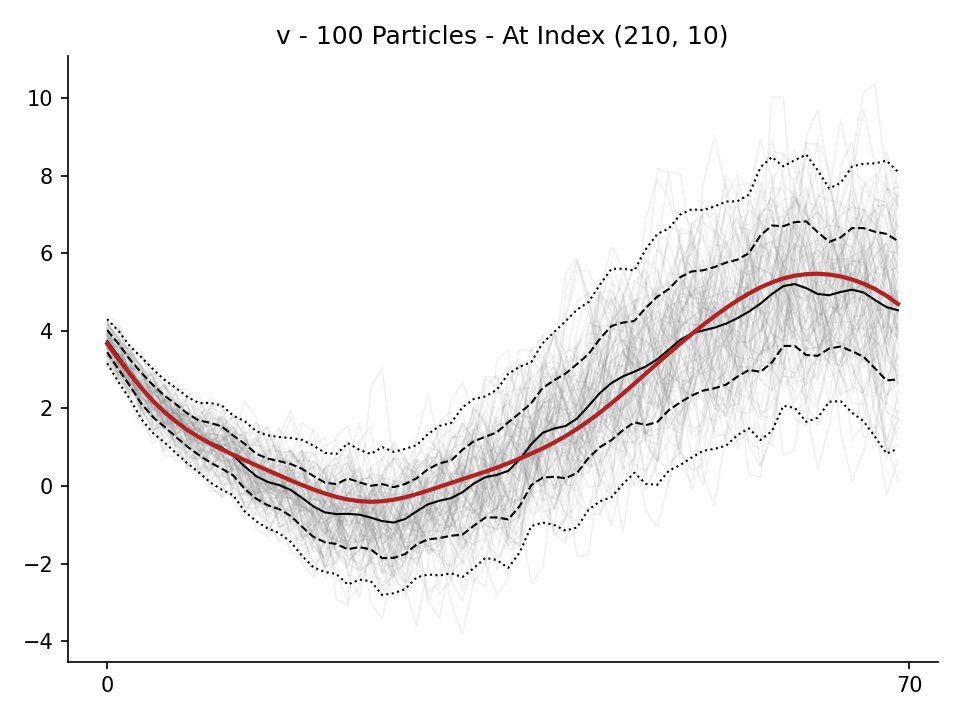}
    \caption{Meridional Vel.~$v$ at $(210,10)$.}
    \label{fig:particle_v_3}
\end{subfigure}
\caption{SPDE Particle trajectories (grey) over time (horizontal axis) on the coarse grid with $c=8$. The red line is the trajectory of the truth, i.e. the fine grid PDE projected onto the SPDE grid. The ensemble mean (black line) and the regions of one (dashed lines) and two (dotted lines) ensemble standard deviations are indicated. The time horizon is $560$ fine PDE timesteps, equivalent to $70$ steps on the coarse SPDE grid. The rows correspond to three different points chosen within the computational domain. The top row is the gridpoint with index $(140,20)$ on the coarse grid. This point lies in the center of the domain. The middle and bottom rows correspond to the points with indices $(140,30)$ and $(210,10)$.} 
\label{fig:particles}

\end{figure}

\begin{figure}[ht!]

\centering
\begin{subfigure}{0.32\textwidth}
    \centering
    \includegraphics[width=\textwidth]{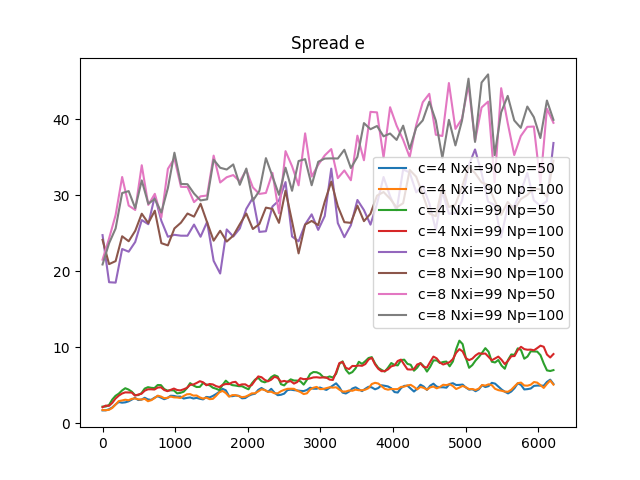}
    \caption{Elevation $\eta$.}
    \label{fig:spread_e}
\end{subfigure}
\begin{subfigure}{0.32\textwidth}
    \centering
    \includegraphics[width=\textwidth]{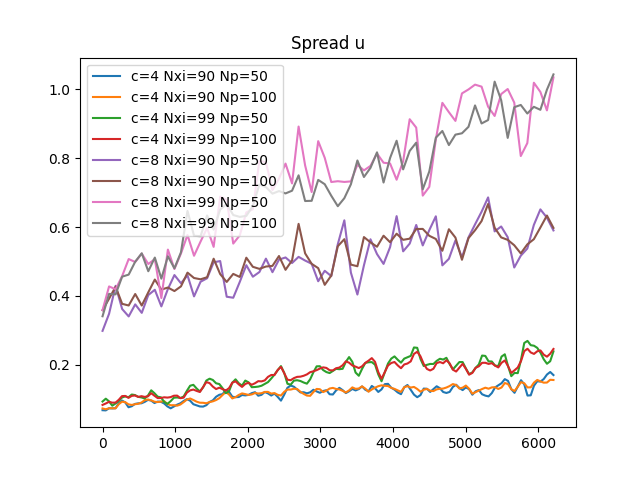}
    \caption{Zonal Velocity $u$.}
    \label{fig:spread_u}
\end{subfigure}
\begin{subfigure}{0.32\textwidth}
    \centering
    \includegraphics[width=\textwidth]{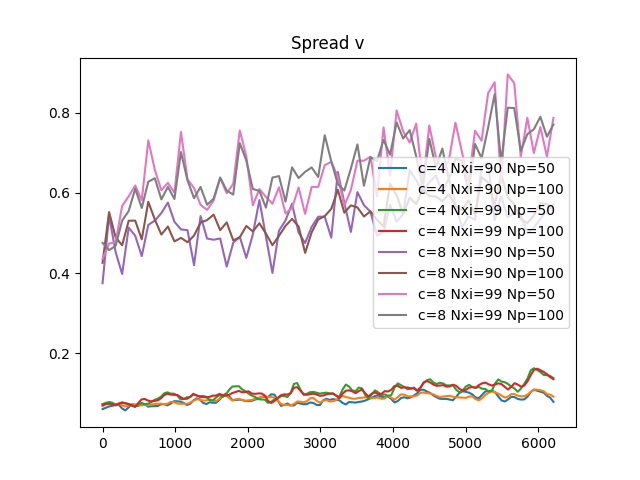}
    \caption{Meridional Velocity $v$.}
    \label{fig:spread_v}
\end{subfigure}

\caption{Ensemble spread in all scenarios on the central grid point with index $(140,20)$ in the case of $c=8$ and index $(280,40)$ in the $c=4$ case.} 
\label{fig:spread}

\end{figure}

\paragraph{Ensemble trajectories and spread.}
In a first experiment, we assess the particle trajectories of an ensemble generated from independent runs of the calibrated stochastic equation run on the coarse grid. The ensemble is compared to the true evolution which was run on the fine grid and sub-sampled, at every time-step, onto the respective coarse grid for comparison. Figure~\ref{fig:particles} shows a sample of the obtained ensemble at different grid locations. We observe that the trajectories of the truth 
are well-contained in the particle ensemble. The signal trajectory stays within the cloud of particles most of the time, consistent with the goal that the true fine resolution run and the stochastic coarse resolution runs are drawn from the same probability density function. In that case, the fine run should fall outside the stochastic ensemble range a fraction 1/(ensemble size +1), which our ensemble achieves approximately.
This experiment is a first indication of the effectiveness of our calibrated stochastic model.
% Not in this paper \todo[inline]{Scientifically, we should maybe compare this against a "random calibration" whatever that means...}
The ensemble spread at a central grid location for all considered scenarios is plotted in Figure~\ref{fig:spread} for all three components of the system. 
The spread exhibits the expected behaviour in all scenarios: there is a steady increase, and the size depends strongly on the coarsening $c$ and the ratio of explained variance $N_\xi$, but there is no dependence on the number particles. Scenario \textbf{a} exhibits the smallest spread. That is because the model reduction is minimal in this case and the variance of the system is not sufficiently explained. Scenarios \textbf{g} and \textbf{h} exhibit the highest  spread as the model reduction is a lot higher in this case and the variance of the system is now well explained. %\pj{(PJ: Perhaps show plots of RMSE and ensemble standard deviation. These should be similar.)}

\begin{figure}[ht!]

\centering
\begin{subfigure}{0.32\textwidth}
    \centering
    \includegraphics[width=\textwidth]{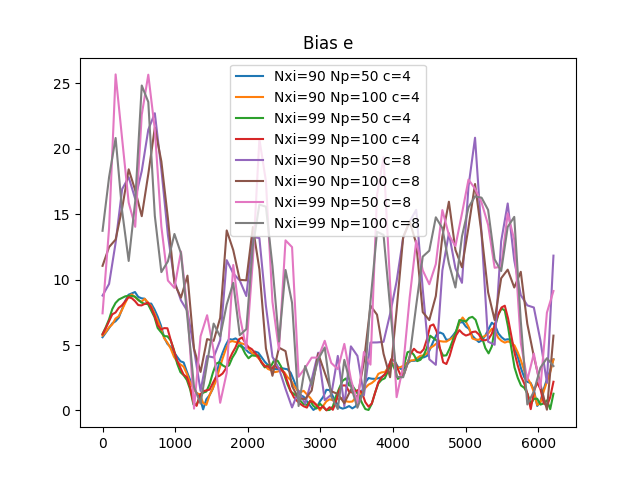}
    \caption{Elevation $\eta$.}
    \label{fig:bias_e}
\end{subfigure}
\begin{subfigure}{0.32\textwidth}
    \centering
    \includegraphics[width=\textwidth]{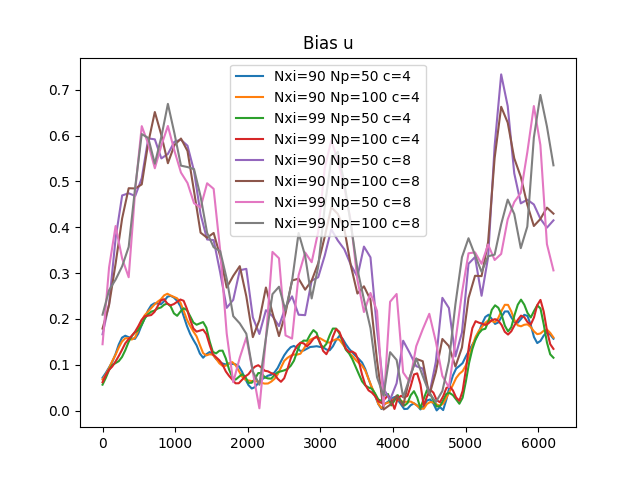}
    \caption{Zonal Velocity $u$.}
    \label{fig:bias_u}
\end{subfigure}
\begin{subfigure}{0.32\textwidth}
    \centering
    \includegraphics[width=\textwidth]{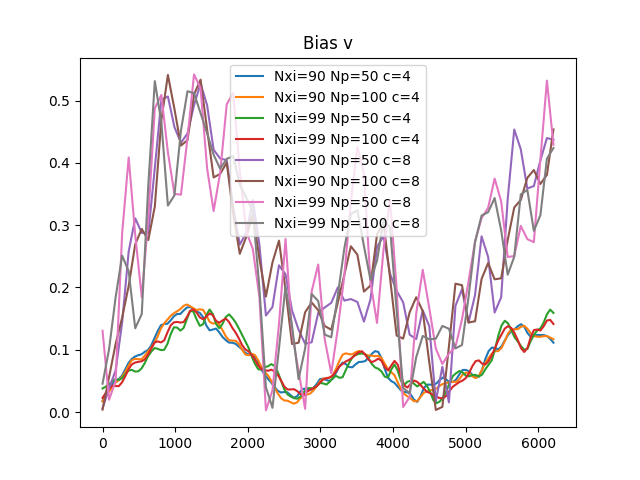}
    \caption{Meridional Velocity $v$.}
    \label{fig:bias_v}
\end{subfigure}

\caption{Ensemble bias with respect to the truth of the ensemble in all scenarios on the central grid point with index $(140,20)$ in the case of $c=8$ and index $(280,40)$ in the $c=4$ case.} 
\label{fig:bias}

\end{figure}

 \paragraph{Ensemble error.}
We assess the ensemble error with respect to the PDE (the truth) using three different widely used error metrics. 
The first is the bias
\begin{equation}
    \frac{1}{N_p}\sum_{i=1}^{N_p} x^i - x^{\text{true}},
\end{equation}
the second is the root-mean-square error (RMSE),
\begin{equation}
    \sqrt{\frac{1}{N_p}\sum_{i=1}^{N_p} (x^i - x^{\text{true}})^2}
\end{equation}
 and the last is the mean relative $L^2$ error
\begin{equation}
    \frac{1}{N_p}\sum_{i=1}^{N_p}\sqrt{  \left(\frac{\sum_{\iota,\jmath}^{N_x, N_y} (x^i_{\iota, \jmath} - x^{\text{true}}_{\iota, \jmath})^2}{ \sum_{\iota,\jmath}^{N_x, N_y} (x^{\text{true}}_{\iota, \jmath})^2}\right)}.
\end{equation}
 Note that both the bias and the RMSE are local errors in space, whereas the relative $L^2$ error accounts for the whole spatial domain.
The bias evaluated the central grid location is shown in Figure~\ref{fig:bias} for all scenarios.
The bias remains stable for all three components of the system. It has a oscillatory behavior. The most striking difference between the parameters is exhibited by the different levels of coarsening. The stronger model reduction in the $c=8$ case leads to higher levels of the bias. There is no strong influence of neither the variance level $N_\xi$ nor the number of particles.

\begin{figure}[ht!]

\centering
\begin{subfigure}{0.32\textwidth}
    \centering
    \includegraphics[width=\textwidth]{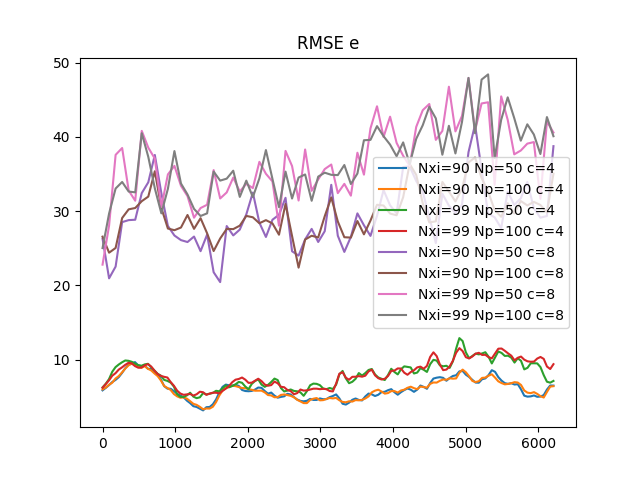}
    \caption{Elevation $\eta$.}
    \label{fig:rmse_e}
\end{subfigure}
\begin{subfigure}{0.32\textwidth}
    \centering
    \includegraphics[width=\textwidth]{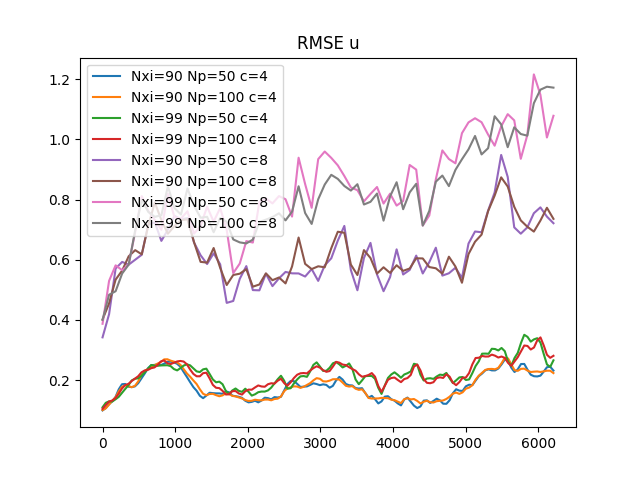}
    \caption{Zonal Velocity $u$.}
    \label{fig:rmse_u}
\end{subfigure}
\begin{subfigure}{0.32\textwidth}
    \centering
    \includegraphics[width=\textwidth]{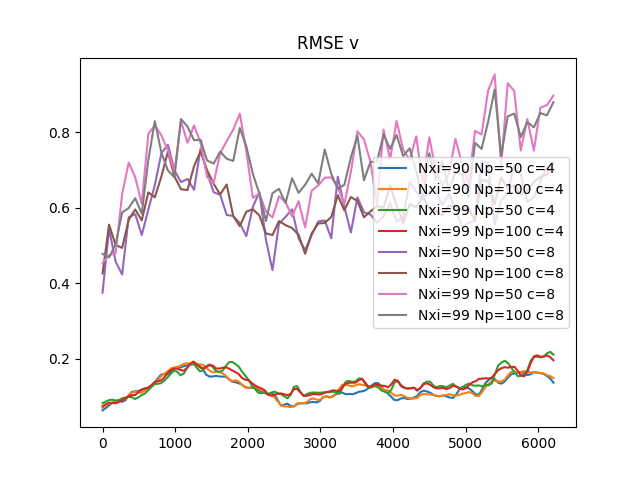}
    \caption{Meridional Velocity $v$.}
    \label{fig:rmse_v}
\end{subfigure}

\caption{Ensemble RMSE with respect to the truth of the ensemble in all scenarios on the central grid point with index $(140,20)$ in the case of $c=8$ and index $(280,40)$ in the $c=4$ case.} 
\label{fig:rmse}

\end{figure}

Moreover, the RMSE evaluated at the same central grid location is depicted in Figure~\ref{fig:rmse} for all three variables in all scenarios.
It steadily increases over time for all three components of the system. In addition to the strong influence of the coarsening on the level of the RMSE, it is also clearly influenced by the variance parameter $N_\xi$. Here, a higher variance parameter translates to more noise components in the SPDE. Importantly, the RMSE estimates are consistent with the ensemble spread estimates in Figure~\ref{fig:spread}, pointing to the accuracy of our uncertainty quantification.

The averages in time of the bias and RMSE are given in Table~\ref{tab:time_mean}.
In addition to the clear effect of the coarsening, we can see again that increasing the variance parameter in the different scenarios increases the mean bias and RMSE.

\begin{figure}[ht!]

\centering
\begin{subfigure}{0.32\textwidth}
    \centering
    \includegraphics[width=\textwidth]{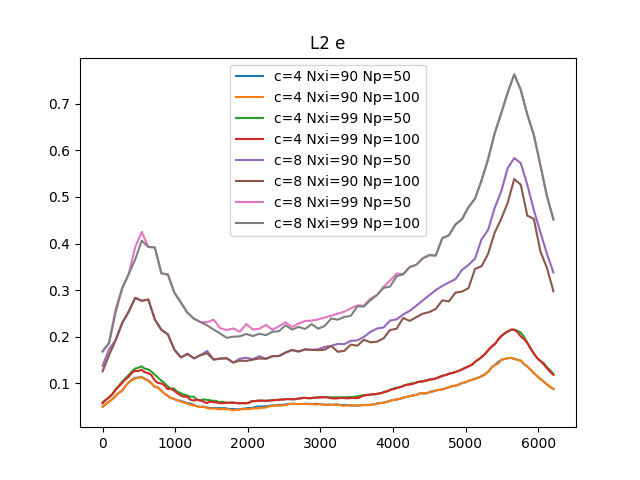}
    \caption{Elevation $\eta$.}
    \label{fig:l2_e}
\end{subfigure}
\begin{subfigure}{0.32\textwidth}
    \centering
    \includegraphics[width=\textwidth]{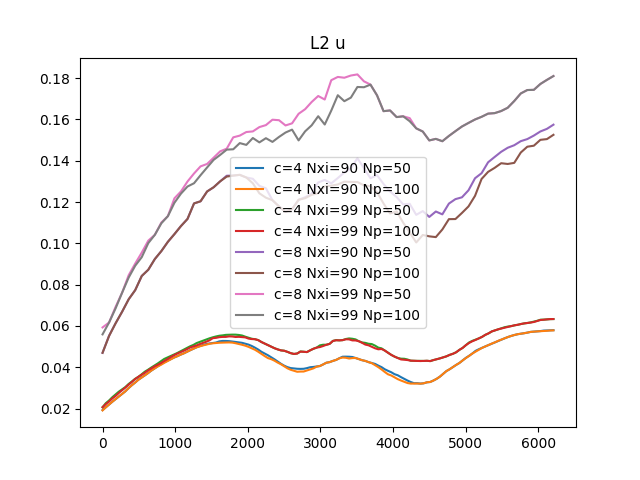}
    \caption{Zonal Velocity $u$.}
    \label{fig:l2_u}
\end{subfigure}
\begin{subfigure}{0.32\textwidth}
    \centering
    \includegraphics[width=\textwidth]{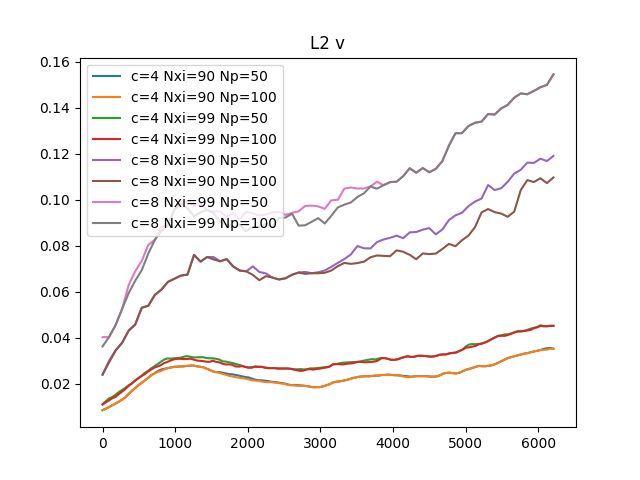}
    \caption{Meridional Velocity $v$.}
    \label{fig:l2_v}
\end{subfigure}

\caption{Relative $L^2$ error with respect to the coarsened truth averaged over the ensemble for all scenarios.} 
\label{fig:L2}

\end{figure}

Lastly, the average relative $L^2$ error is shown in Figure~\ref{fig:L2}. The $L^2$ error exhibits an upward trend.
Here again, the coarsening has the strongest effect, and variance parameter plays a secondary but still significant role.

\begin{table*}\centering
\small
\ra{1.3}
\tabcolsep=0.11cm
\begin{tabular}{@{}rrrcrrcrrcrr@{}}\toprule
& \multicolumn{5}{c}{Bias} & \phantom{x} & \multicolumn{5}{c}{RMSE}\\
\cmidrule{2-6} \cmidrule{8-12}& \multicolumn{2}{c}{$N_{\xi} = 0.90$} & \phantom{x} & \multicolumn{2}{c}{$N_{\xi} = 0.99$} & \phantom{x} & \multicolumn{2}{c}{$N_{\xi} = 0.90$} & \phantom{x} & \multicolumn{2}{c}{$N_{\xi} = 0.99$}\\
\cmidrule{2-3} \cmidrule{5-6} \cmidrule{8-9} \cmidrule{11-12}   & $N_p=50$      & $N_p=100$            &               & $N_p=50$ & $N_p=100$                 &               & $N_p=50$ & $N_p=100$                 &               & $N_p=50$ & $N_p=100$\\ \midrule
$c=4$\\
$e$ & 3.8093 & 3.6988  && 3.8272 & 3.7311 && 8.9550 & 8.9615  && 11.4424 & 11.5070 \\
$u$ & 0.2579 & 0.2576 && 0.2569 & 0.2571 && 0.3420 & 0.3428 && 0.3883 & 0.3887\\
$v$ & 0.1198 & 0.1187 && 0.1242 & 0.1212 && 0.2321 & 0.2318 && 0.2851 & 0.2863 \\
$c=8$\\
$e$ & 9.1173 & 8.5740 && 9.9699 & 8.9823 && 34.7980 & 35.0126 && 43.7557 & 44.2226\\
$u$ & 0.6194 & 0.6122 && 0.6268 & 0.6113 && 1.0713 & 1.0703 && 1.2443 & 1.2473\\
$v$ & 0.3119 & 0.3057 && 0.3228 & 0.3044 && 0.8349 & 0.8353 && 1.0562 & 1.0591\\
\bottomrule
\end{tabular}
\caption{Time averages of the Bias and RMSE at central grid location for all scenarios.}
\label{tab:time_mean}
\end{table*}

\begin{figure}[ht!]

\centering
\begin{subfigure}{0.24\textwidth}
    \centering
    \includegraphics[width=\textwidth]{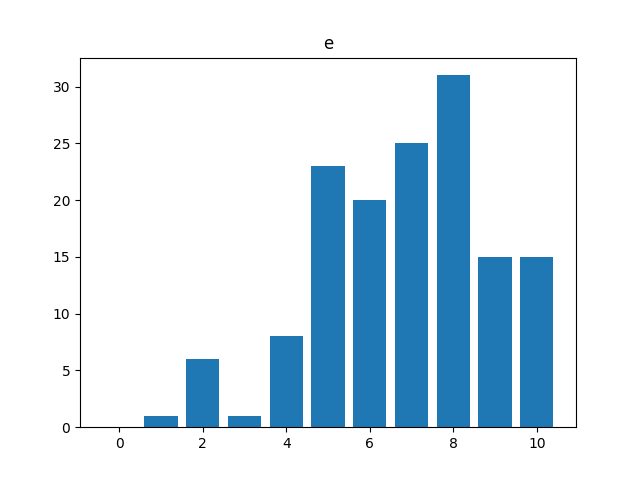}
    \caption{$\eta$: 64 Steps at $(140,20)$.}
    \label{fig:hist8_e_0}
\end{subfigure}
\begin{subfigure}{0.24\textwidth}
    \centering
    \includegraphics[width=\textwidth]{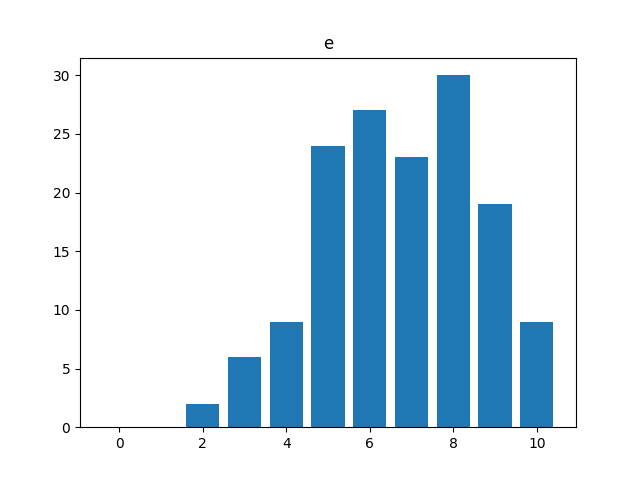}
    \caption{$\eta$: 128 Steps at $(140,20)$}
    \label{fig:hist16_e_0}
\end{subfigure}
\begin{subfigure}{0.24\textwidth}
    \centering
    \includegraphics[width=\textwidth]{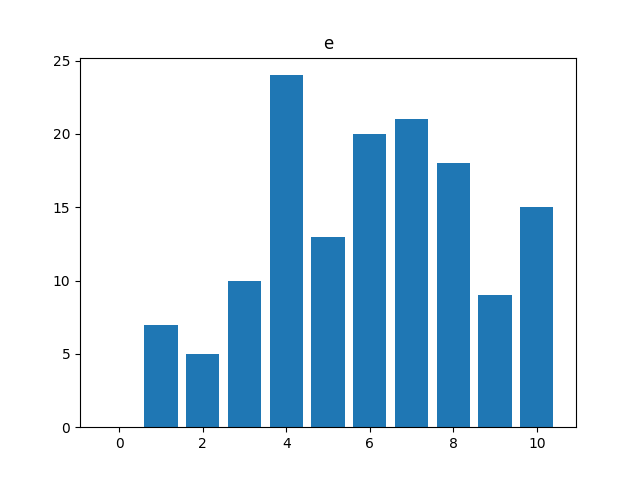}
    \caption{$\eta$: 64 Steps at $(140,30)$.}
    \label{fig:hist8_e_3}
\end{subfigure}
\begin{subfigure}{0.24\textwidth}
    \centering
    \includegraphics[width=\textwidth]{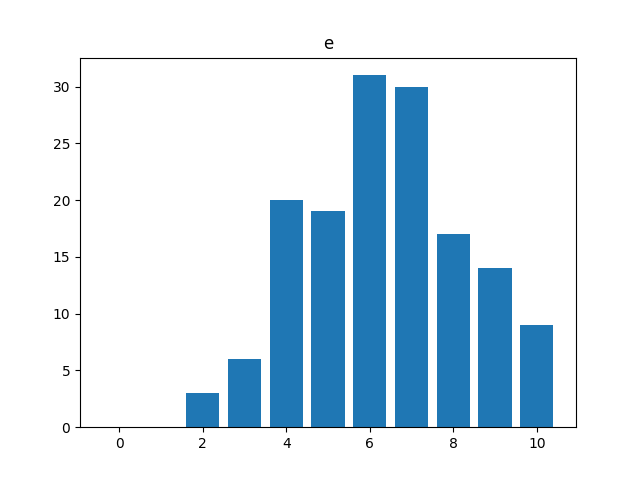}
    \caption{$\eta$: 128 Steps at $(140,30)$.}
    \label{fig:hist16_e_3}
\end{subfigure}
\hfill
\begin{subfigure}{0.24\textwidth}
    \centering
    \includegraphics[width=\textwidth]{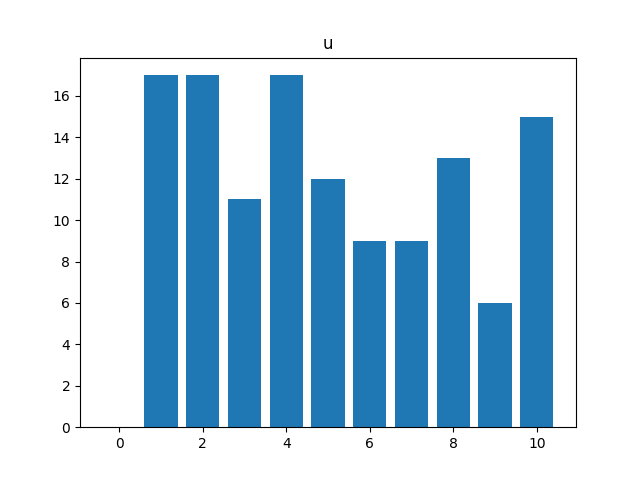}
    \caption{$u$: 64 Steps at $(140,20)$.}
    \label{fig:hist8_u_0}
\end{subfigure}
\begin{subfigure}{0.24\textwidth}
    \centering
    \includegraphics[width=\textwidth]{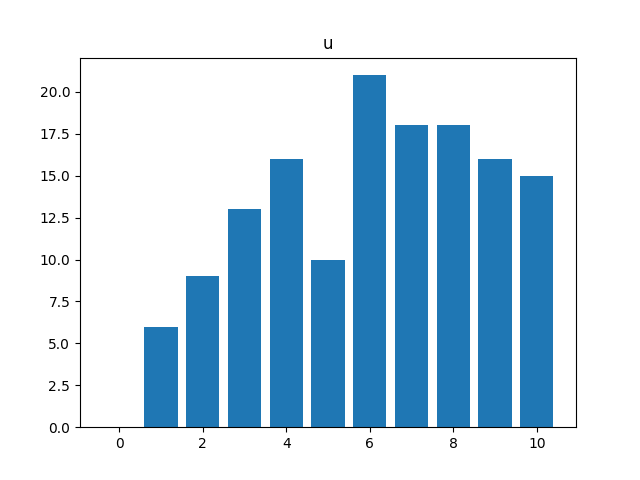}
    \caption{$u$: 128 Steps at $(140,20)$}
    \label{fig:hist16_u_0}
\end{subfigure}
\begin{subfigure}{0.24\textwidth}
    \centering
    \includegraphics[width=\textwidth]{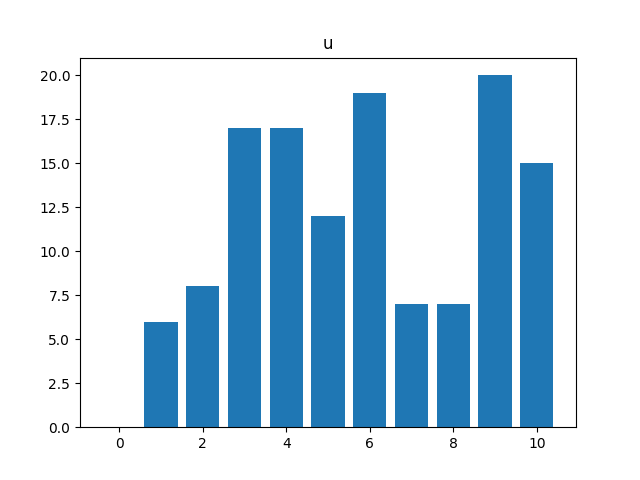}
    \caption{$u$: 64 Steps at $(140,30)$.}
    \label{fig:hist8_u_3}
\end{subfigure}
\begin{subfigure}{0.24\textwidth}
    \centering
    \includegraphics[width=\textwidth]{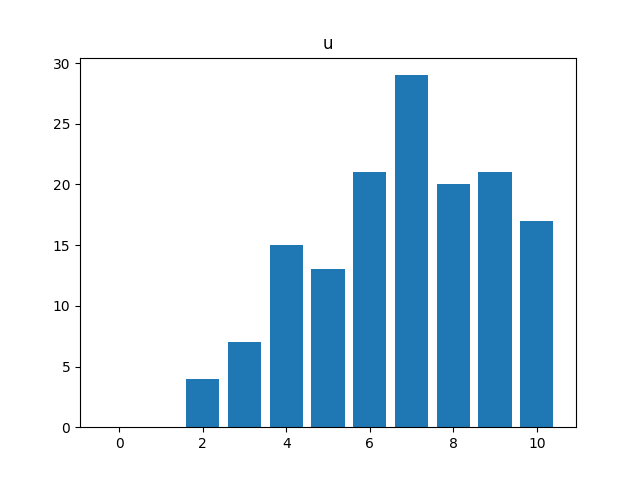}
    \caption{$u$: 128 Steps at $(140,30)$.}
    \label{fig:hist16_u_3}
\end{subfigure}
\hfill
\begin{subfigure}{0.24\textwidth}
    \centering
    \includegraphics[width=\textwidth]{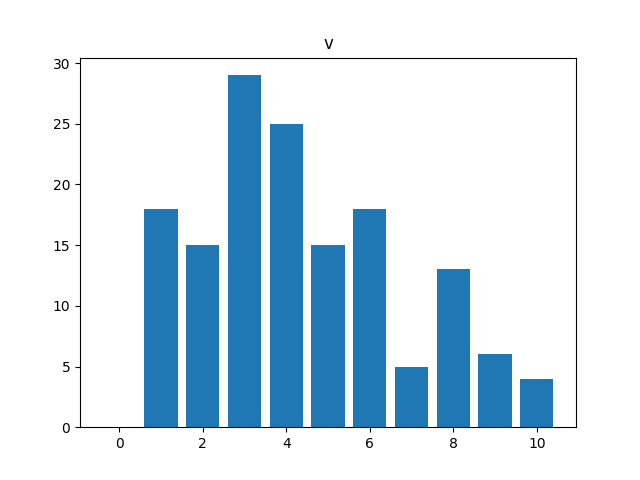}
    \caption{$v$: 64 Steps at $(140,20)$.}
    \label{fig:hist8_v_0}
\end{subfigure}
\begin{subfigure}{0.24\textwidth}
    \centering
    \includegraphics[width=\textwidth]{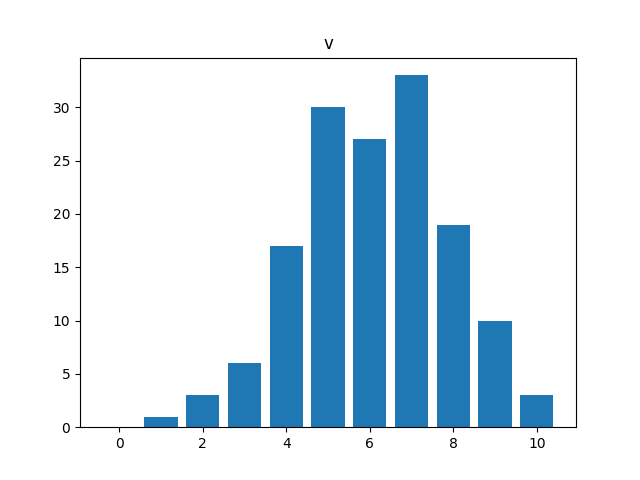}
    \caption{$v$: 128 Steps at $(140,20)$}
    \label{fig:hist16_v_0}
\end{subfigure}
\begin{subfigure}{0.24\textwidth}
    \centering
    \includegraphics[width=\textwidth]{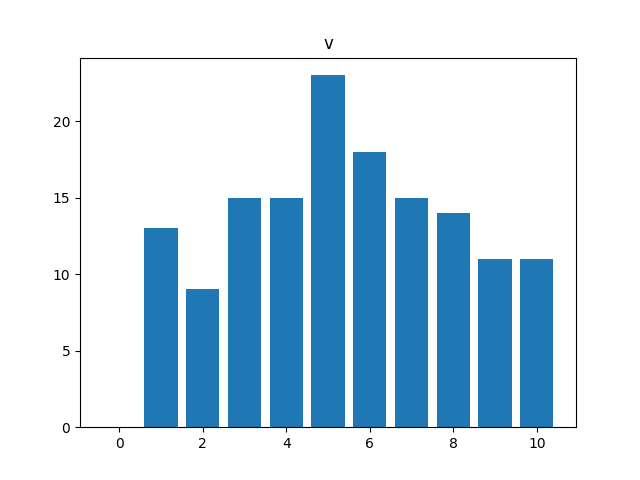}
    \caption{$v$: 64 Steps at $(140,30)$.}
    \label{fig:hist8_v_3}
\end{subfigure}
\begin{subfigure}{0.24\textwidth}
    \centering
    \includegraphics[width=\textwidth]{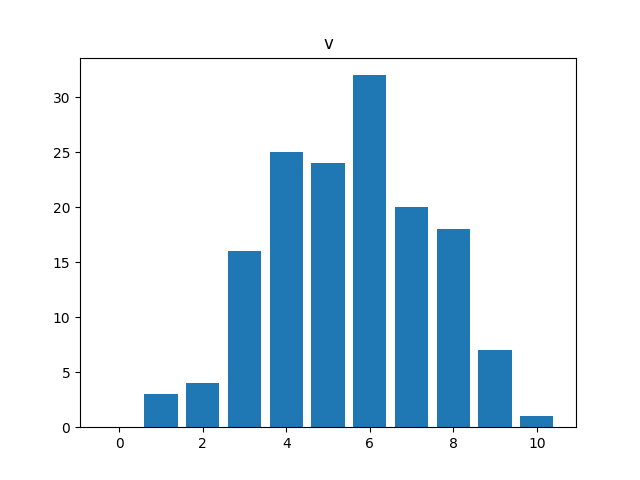}
    \caption{$v$: 128 Steps at $(140,30)$.}
    \label{fig:hist16_v_3}
\end{subfigure}

\caption{Rank Histograms for all three shallow water variables. The rank histograms are estimated using $64$ and $128$ forecast steps, respectively. Shown are the histograms at two different points in the domain. One is the central point with index $(140,20)$ and the other one is the northward shifted point $(140,30)$.} 
\label{fig:rank_histograms}

\end{figure}

\paragraph{Rank Histograms.}
In order to asses the quality of the ensemble for Data Assimilation we estimated rank histograms for an ensemble of $10$ particles on the $c=8$ coarse grid using $N_\xi=0.99$. The rank histogram ranks the true fine resolution run among the ensemble members, and the statistics is gathered over time. The results are depicted in Figure~\ref{fig:rank_histograms} for two different grid locations as well as two different time horizons. The plots show that the true run ranks closer to the ensemble mean, showing that the ensemble spread is slightly too large. However, that is considered a good sign at it means that the uncertainty estimate will be robust against extreme behavior of the true run.

% \begin{table*}\centering
% \ra{1.3}
% \begin{tabular}{@{}rrrcrr@{}}\toprule
% & \multicolumn{2}{c}{$N_{\xi} = 0.90$} & \phantom{abc} & \multicolumn{2}{c}{$N_{\xi} = 0.99$} \\
% \cmidrule{2-3} \cmidrule{5-6}  & $N_p=50$ & $N_p=100$  && $N_p=50$ & $N_p=100$\\ \midrule
% $c=4$\\
% $e$ & 3.8093 & 3.6988  && 3.8272 & 3.7311 \\
% $u$ & 0.2579 & 0.2576 && 0.2569 & 0.2571 \\
% $v$ & 0.1198 & 0.1187 && 0.1242 & 0.1212 \\
% $c=8$\\
% $e$ & 9.1173 & 8.5740 && 9.9699 & 8.9823\\
% $u$ & 0.6194 & 0.6122 && 0.6268 & 0.6113 \\
% $v$ & 0.3119 & 0.3057 && 0.3228 & 0.3044 \\
% \bottomrule
% \end{tabular}
% \caption{Bias average in time}
% \label{tab:time_mean_bias}
% \end{table*}

% \begin{table*}\centering
% \ra{1.3}
% \begin{tabular}{@{}rrrcrr@{}}\toprule
% & \multicolumn{2}{c}{$N_{\xi} = 0.90$} & \phantom{abc} & \multicolumn{2}{c}{$N_{\xi} = 0.99$} \\
% \cmidrule{2-3} \cmidrule{5-6}  & $N_p=50$ & $N_p=100$  && $N_p=50$ & $N_p=100$\\ \midrule
% $c=4$\\
% $e$ & 8.9550 & 8.9615  && 11.4424 & 11.5070 \\
% $u$ & 0.3420 & 0.3428 && 0.3883 & 0.3887 \\
% $v$ & 0.2321 & 0.2318 && 0.2851 & 0.2863 \\
% $c=8$\\
% $e$ & 34.7980 & 35.0126 && 43.7557 & 44.2226\\
% $u$ & 1.0713 & 1.0703 && 1.2443 & 1.2473 \\
% $v$ & 0.8349 & 0.8353 && 1.0562 & 1.0591 \\
% \bottomrule
% \end{tabular}
% \caption{RMSE average in time}
% \label{tab:time_mean_rmse}
% \end{table*}

\section{Conclusions and Future Work}\label{sect:conclusions}

In this paper we introduce a procedure for calibrating stochastic parametrizations. The method can be applied to a large class 
of stochastic parametrizations and it is also agnostic as to the source of data (real or synthetic). It is based on a principal component analysis technique to generate the eigenvectors and the eigenvalues of the covariance matrix of the stochastic parametrization. We test the procedure 
for a stochastic parametrization applied to the rotating shallow water model. The stochastic parametrization tested in this 
paper model the unresolved scales due to model reduction. We calibrate the noise by using the elevation variable of the model, 
as this is an observable easily obtainable in practical application, and use synthetic data as input for the calibration. 
We test the calibration using standard uncertainty quantification tests (RMSE, bias, spread, histograms) and obtain good results. 
 
The work presented in this paper was limited in scope. We introduced the  calibration methodology and used a particular stochastic parametrization for the rotating shallow water model as a test case to assess it. Further work is warranted as a continuation of the project. We enumerate below several possible directions and open problems.  

\begin{itemize}

\item Since we have a generic calibration procedure we may want to compare different stochastic parametrizations. In particular, can we identify the optimal stochastic parametrization for a particular state space model ? For example, can we find the parametrization that requires the least number of sources of noise to explain the model uncertainty? 

%Also it would be useful to compare the calibrated stochastic parametrization with an adhoc fitted stochastic parametrization, but I won't say this here. 

\item We can consider different parameters for the rotating shallow water model, e.g. Rossby number, Reynolds number, rotational Froude number, etc. Also we can explore different initial conditions that lead to more realistic realisations of the model state. For example, we can use an initial condition that generates a meandering eastward jet as initial condition (a synthetic approximation of the Jetstream), or one that comes from real data.  

\item We can explore various stochastic transport parametrizations for the rotating shallow water model, for example parametrizations that account for rotation and non-zero divergence in the stochastic parametrization for $(u,v)$. We can also explore using the velocity field $(u,v)$ as input data instead on the fluid elevation.

\item Now that the calibration step has been completed, we can proceed with the implementation of Data Assimilation methodology with synthetic and real data for the test case presented in this paper. 

\end{itemize}

%%%%%%%%%%%%%%%%%%%%%%%%%%%%%%%%%%%%%%%%%
%\input{appendix.tex}
% OL: Can we please move it at the very end only after we finish? We might modify the appendix as well (i.e. correlate it with other changes in the text etc), and having 2 pages of references in between them is disturbing sometimes :) 
% Ok, this time :)

%%%%%%%%%%%%%%%%%%%%%%%%%%%%%%%%%%%%%%%%%%%%%%%%%%
\clearpage

\paragraph{Acknowledgments.} 
We would like to thank Colin Cotter, Wei Pan, and James Woodfield for their assistance in devising the numerical scheme for the hyperbolic equation. We would like to also thank Long Li, Darryl Holm, Bertrand Chapron, \'{E}tienne Mémin, Baylor Fox-Kemper for many fruitful discussions we had during the preparation of this work.

\paragraph{Funding.} 
Alexander Lobbe, Oana Lang and Dan Crisan were partially supported by the European Research Council (ERC) under the European Union’s Horizon 2020 Research and Innovation Programme (ERC, Grant Agreement No 856408). Peter Jan van Leeuwen was supported through the European Research Council project CUNDA grant 694509 under the Horizon 2020 program, and via NSF grant 1924659. Roland Potthast was supported by the CONTRAILS project, Federal Ministry for Economy and Climate BMWK, Germany, and by German Science Foundation Grant FOR 2589. 

\paragraph{Conflict of interest statement.}
On behalf of all authors, the corresponding author states that there is no conflict of interest.

\bibliographystyle{alpha}
\bibliography{main}

\newcommand{\etalchar}[1]{$^{#1}$}
\begin{thebibliography}{CCH{\etalchar{+}}20b}

\bibitem[AL77]{arakawa1977computational}
Akio Arakawa and Vivian~R Lamb.
\newblock Computational design of the basic dynamical processes of the ucla
  general circulation model.
\newblock {\em General circulation models of the atmosphere}, 17(Supplement
  C):173--265, 1977.

\bibitem[BLBM21]{EtienneLong}
R{\"u}diger Brecht, Long Li, Werner Bauer, and Etienne M{\'e}min.
\newblock Rotating shallow water flow under location uncertainty with a
  structure-preserving discretization.
\newblock {\em Journal of Advances in Modeling Earth Systems},
  13(12):e2021MS002492, 2021.

\bibitem[BMP99]{Buizza}
Roberto Buizza, M~Milleer, and Tim~N Palmer.
\newblock Stochastic representation of model uncertainties in the ecmwf
  ensemble prediction system.
\newblock {\em Quarterly Journal of the Royal Meteorological Society},
  125(560):2887--2908, 1999.

\bibitem[CCH{\etalchar{+}}19]{Wei1}
Colin Cotter, Dan Crisan, Darryl~D Holm, Wei Pan, and Igor Shevchenko.
\newblock Numerically modeling stochastic lie transport in fluid dynamics.
\newblock {\em Multiscale Modeling \& Simulation}, 17(1):192--232, 2019.

\bibitem[CCH{\etalchar{+}}20a]{Wei2}
Colin Cotter, Dan Crisan, Darryl Holm, Wei Pan, and Igor Shevchenko.
\newblock Modelling uncertainty using stochastic transport noise in a 2-layer
  quasi-geostrophic model.
\newblock {\em Foundations of Data Science}, 2(2):173, 2020.

\bibitem[CCH{\etalchar{+}}20b]{Wei3}
Colin Cotter, Dan Crisan, Darryl~D Holm, Wei Pan, and Igor Shevchenko.
\newblock A particle filter for stochastic advection by lie transport: a case
  study for the damped and forced incompressible two-dimensional euler
  equation.
\newblock {\em SIAM/ASA Journal on Uncertainty Quantification},
  8(4):1446--1492, 2020.

\bibitem[CCH{\etalchar{+}}20c]{Weifinal}
Colin Cotter, Dan Crisan, Darryl~D Holm, Wei Pan, and Igor Shevchenko.
\newblock A particle filter for stochastic advection by lie transport: a case
  study for the damped and forced incompressible two-dimensional euler
  equation.
\newblock {\em SIAM/ASA Journal on Uncertainty Quantification},
  8(4):1446--1492, 2020.

\bibitem[Dur10]{Durran}
Dale~R Durran.
\newblock {\em Numerical methods for fluid dynamics: With applications to
  geophysics}, volume~32.
\newblock Springer Science \& Business Media, 2010.

\bibitem[Her20]{longlist}
Hans et.~al. Hersbach.
\newblock The era5 global reanalysis.
\newblock {\em Q. J. R.Meteorol. Soc.}, 146, 2020.

\bibitem[Hol15]{Holm2015}
Darryl~D Holm.
\newblock Variational principles for stochastic fluid dynamics.
\newblock {\em Proceedings of the Royal Society A: Mathematical, Physical and
  Engineering Sciences}, 471(2176):20140963, 2015.

\bibitem[HW65]{HarlowWelch}
Francis~H Harlow and J~Eddie Welch.
\newblock Numerical calculation of time-dependent viscous incompressible flow
  of fluid with free surface.
\newblock {\em The physics of fluids}, 8(12):2182--2189, 1965.

\bibitem[Kal03]{Kalnay}
Eugenia Kalnay.
\newblock {\em Atmospheric modeling, data assimilation and predictability}.
\newblock Cambridge university press, 2003.

\bibitem[KP92]{KloedenPlaten}
Peter~E. Kloeden and Eckhard Platen.
\newblock {\em Numerical solution of stochastic differential equations}.
\newblock Applications of mathematics 23. Springer-Verlag, Berlin, 1992.

\bibitem[LvLCP23]{lucrisanlangmemin}
Oana Lang, Peter~Jan van Leeuwen, Dan Crisan, and Roland Potthast.
\newblock Bayesian inference for fluid dynamics: a case study for the
  stochastic rotating shallow water model.
\newblock {\em Frontiers in Applied Mathematics and Statistics}, 2023.

\bibitem[M{\'e}m14]{Memin2014}
Etienne M{\'e}min.
\newblock Fluid flow dynamics under location uncertainty.
\newblock {\em Geophysical \& Astrophysical Fluid Dynamics}, 108(2):119--146,
  2014.

\bibitem[MTVE01]{MajdaDA}
Andrew~J Majda, Ilya Timofeyev, and Eric Vanden~Eijnden.
\newblock A mathematical framework for stochastic climate models.
\newblock {\em Communications on Pure and Applied Mathematics: A Journal Issued
  by the Courant Institute of Mathematical Sciences}, 54(8):891--974, 2001.

\bibitem[Pal19]{Palmer}
Tim Palmer.
\newblock The ecmwf ensemble prediction system: Looking back (more than) 25
  years and projecting forward 25 years.
\newblock {\em Quarterly Journal of the Royal Meteorological Society},
  145:12--24, 2019.

\bibitem[PvL22]{pathiraja}
S.~Pathiraja and P.J. van Leeuwen.
\newblock Model uncertainty estimation in data assimilation for multiscale
  systems with partially observed resolved variables.
\newblock {\em JAMES}, 2022.

\bibitem[RLJ{\etalchar{+}}21]{Resseguir2021}
Valentin Resseguier, Long Li, Gabriel Jouan, Pierre D{\'e}rian, Etienne
  M{\'e}min, and Bertrand Chapron.
\newblock New trends in ensemble forecast strategy: uncertainty quantification
  for coarse-grid computational fluid dynamics.
\newblock {\em Archives of Computational Methods in Engineering},
  28(1):215--261, 2021.

\bibitem[Val17]{Vallis}
Geoffrey~K Vallis.
\newblock {\em Atmospheric and oceanic fluid dynamics}.
\newblock Cambridge University Press, 2017.

\bibitem[Zei18]{Zeitlinbook}
Vladimir Zeitlin.
\newblock {\em Geophysical fluid dynamics: understanding (almost) everything
  with rotating shallow water models}.
\newblock Oxford University Press, 2018.

\end{thebibliography}
%%%%%%%%%%%%%%%%%%%%%%%%%%%%%%%%%%%%%%%%%%%%%%%%%%

\clearpage
\begin{appendices}\label{sect:appendix}

\section{Discretization}\label{appendix:discretisation}
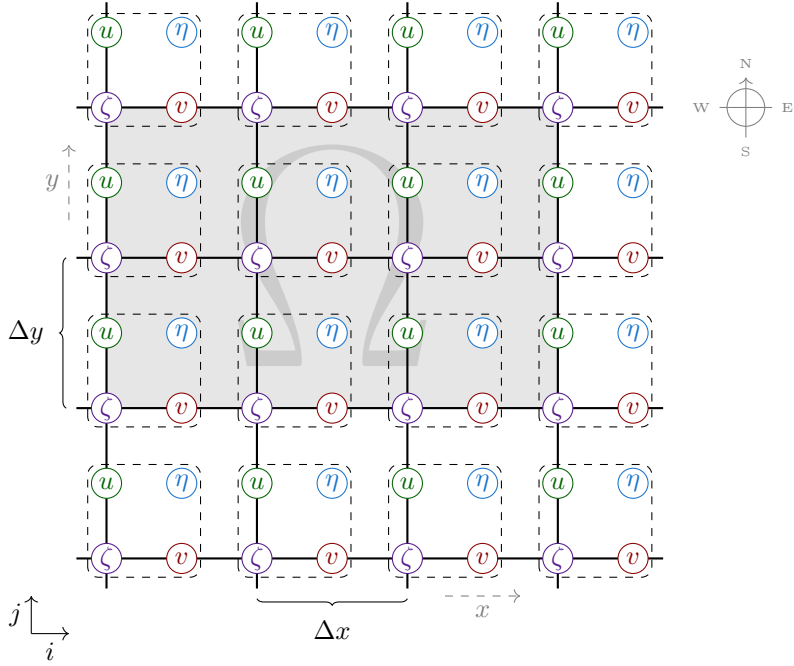
\begin{figure}[h]
    \centering
    \begin{tikzpicture}
	\draw[fill=gray!20] (3,3) rectangle (9,7);
	\node[color=gray!40] at (6,5) {\fontsize{120pt}{144pt}\selectfont $\Omega$};
    %\draw[color=lightgray] (0,0) grid (10,10);
    \foreach \x in {4,6,8,10}{\draw[color=black, thick] (\x-1,0.6) -- (\x-1,8.4);}
     \foreach \x in {2, 4,6,8}{\draw[color=black, thick] (2.6,\x-1) -- (10.4,\x-1);}
    \foreach \x in {3,5,7,9}{
    	\foreach \y in {2,4,6,8}{
    		\filldraw[draw=DodgerBlue3, fill=white] (\x+1,\y) circle [radius=0.2];
    		\node[color=DodgerBlue3] at (\x+1,\y) {$\eta$};
    		\filldraw[draw=DarkGreen, fill=white] (\x,\y) circle [radius=0.2];
    		\node[color=DarkGreen] at (\x,\y) {$u$};
    		\filldraw[draw=DarkRed, fill=white] (\x+1,\y-1) circle [radius=0.2];
    		\node[color=DarkRed] at (\x+1,\y-1) {$v$};
            \filldraw[draw=Purple4, fill=white] (\x,\y-1) circle [radius=0.2];
            \node[color=Purple4] at (\x,\y-1) {$\zeta$};}}
    \foreach \x in {2,4 ,6,8} \foreach \y in {0,2,4,6} \draw[dashed, rounded corners] (\x+0.75, \y+0.75) rectangle +(1.5, 1.5);
    \draw[->] (2, 0) -- (2.5, 0) node[midway, below]{$i$}; 
    \draw[->] (2, 0) -- (2, 0.5) node[midway, left]{$j$}; 
    \draw[->, dashed, color=gray] (2.5, 5.5) -- (2.5, 6.5) node[midway, left]{$y$};
    \draw[->, dashed, color=gray] (7.5, 0.5) -- (8.5, 0.5) node[midway, below]{$x$};
    \draw[color=gray] (11.5,7) circle [radius = 0.25];
    \draw[->,color=gray] (11.5,6.6) node[below]{\tiny S} -- ++(0,0.8)node[above]{\tiny N};
    \draw[color=gray] (11.15, 7) node[left]{\tiny W} -- ++(0.7,0) node[right]{\tiny E};
    \draw[decoration={brace,mirror,raise=1pt},decorate] (5,0.5) -- node[below=6pt] {$\Delta x$} (7,0.5);
    \draw[decoration={brace,raise=1pt},decorate] (2.5,3) -- node[left=6pt] {$\Delta y$} (2.5,5);
\end{tikzpicture}
    \caption{This is a sketch of our computational grid for the SRSW system. The computational domain is $\Omega$ with coordinates $x$ and $y$. The grid spacings are $\Delta x$ and $\Delta y$ in the respective coordinate directions. The fluid surface elevation $\eta$ is located at the cell center, the grids for the zonal velocity $u$ and the meridional velocity $v$ are shifted  half a grid length westward and southward, respectively. The grid for the potential vorticity $\zeta$ is shifted half a grid length in both the west and south directions. The variable $h$ shares the $\eta$-grid, and the variables $f,\xi$ share the $\zeta$-grid. The grid indices are $i$ in the eastward direction and $j$ in the northward direction. A gridbox with indices $(i,j)$ is depicted as a dashed square.}
    \label{fig:staggered-grid}
\end{figure}

%%%%%%%%%%%%%%%%%%%%%%%%%%%%%%%%%%%%%%%%%%%%%%%%%%
%%%%%%%%%%%%%%%%%%%%%%%%%%%%%%%%%%%%%%%%%%%%%%%%%%
    All compressible and incompressible flow simulations can be subject to large discretization errors when all system variables are defined at the same grid points or at the same time levels. More precisely, spurious modes associated with the pressure field can appear when using a collocated mesh due to the intrinsic structure of the difference scheme which generates an odd-even decoupling between pressure and velocity. 
    To cope with the problem, we use a staggered Arakawa C-grid~\cite{arakawa1977computational} (see Figure \ref{fig:staggered-grid}) in which the velocity components are staggered compared to the pressure term i.e. the pressure is stored in the center of the cell, while the velocity components are stored at the cell faces. This generates a reduction of the dispersion errors as well as improvements in the accuracy of the short-wavelength components of the solution\footnote{For more details on this type of numerical issues see e.g. \cite{Durran} or \cite{HarlowWelch}.}. The staggered arrangement, although not easy to implement, is particularly useful in our case as it enables a more accurate implementation of the high-frequency small-scale modes. 
    
For the spatial discretization we use a staggered Arakawa $\mathrm{C}$ grid (see Figure \ref{fig:staggered-grid}), which is especially accurate for motions with length scales of the order of the Rossby radius of deformation $\left(R_d=g h / f_0\right)$ and smaller. This means that $u, v, f, \zeta$, and $h, \eta$ use different grids, shifted from each other. In particular, taking the $(h, \eta)$ grid as reference, the $u$ grid is shifted westward half a grid length, and the $\mathbf{v}$ grid is shifted southward half a grid length, and the $(f, \xi, \zeta)$ grid is shifted both westward and southward half a grid length. %It is best to draw a picture of this grid to understand what follows.
In order to rigorously explain our numerical implementation we denote\footnote{Note that this corresponds precisely to $u=(u^1,u^2)$ in the previous sections.} $\vb{u}=(u,v)$ and rewrite the deterministic model as follows: 
    \begin{equation}
    \begin{aligned}
        \vb{u}_t &= -(\vb{u}\vdot\nabla)\vb{u} + f({\vb{u}}\cross\vu{z}) - g\grad{h}+D\laplacian{\vb{u}}\\
        h_t &= -\div{(h\vb{u})}
    \end{aligned}
    \end{equation}
    So that in component form, using the fact that, due to $b\equiv0$, $h=H+\eta$ and $H$ is constant.
The 1-layer shallow-water model equations are as follows:
\begin{equation}
\begin{aligned}
u_t &=-u u_x-v u_y+f v-g \eta_x+D\left(u_{x x}+u_{y y}\right) \\
v_t &=-u v_x-v v_y-f u-g \eta_y+D\left(v_{x x}+v_{y y}\right) \\
\eta_t &=-(h u)_x-(h v)_y
\end{aligned}
\end{equation}
in which the total height of the layer is $h=H+\eta$, where $H$ is constant in space and time. We assume the $\beta$-plane approximation for the Coriolis parameter such that $f=f_0+\beta y$, in which $f_0$ is constant in space and time, as is $\beta$. The variables $u, v, y$ denote zonal velocity, meridional velocity, and meridional coordinate, respectively.

In absence of diffusion and friction we want the numerical scheme to conserve energy and potential vorticity. To this end we rewrite the advection terms and the Coriolis force as:
\begin{equation}
\begin{aligned}
-u u_x-v u_y+f v &=-\frac{1}{2}\left(u^2+v^2\right)_x+v v_x-v u_y+f v \\
&=-\frac{1}{2}\left(u^2+v^2\right)_x+v(f+\zeta)
\end{aligned}
\end{equation}
in which $\zeta=v_x-u_y$ is the vertical component of the relative vorticity $\boldsymbol{\xi}=\curl{\vb{u}}=\zeta\vu{z}$.
This leads to the equation set:
%todo{Alex: In my calculation the $\xi$ in this and the next equation should be replaced with $\zeta$.}
\begin{equation}
\begin{aligned}
u_t &=-\frac{1}{2}\left(u^2+v^2\right)_x+v(f+\xi)-g \eta_x+D\left(u_{x x}+u_{y y}\right) \\
v_t &=-\frac{1}{2}\left(u^2+v^2\right)_y-u(f+\xi)-g \eta_y+D\left(v_{x x}+v_{y y}\right) \\
\eta_t &=-(h u)_x-(h v)_y
\end{aligned}
\end{equation}
We then write the vorticity terms as potential vorticity advection:
$$
\begin{aligned}
u_t &=-\frac{1}{2}\left(u^2+v^2\right)_x+h v \frac{f+\xi}{h}-g \eta_x+D\left(u_{x x}+u_{y y}\right) \\
v_t &=-\frac{1}{2}\left(u^2+v^2\right)_y-h u \frac{f+\xi}{h}-g \eta_y+D\left(v_{x x}+v_{y y}\right) \\
\eta_t &=-(h u)_x-(h v)_y
\end{aligned}
$$
We can rewrite this further by collecting the gradient terms as:
$$
\begin{aligned}
&u_t=-\left(\frac{1}{2}\left(u^2+v^2\right)+g \eta\right)_x+h v \zeta+D\left(u_{x x}+u_{y y}\right) \\
&v_t=-\left(\frac{1}{2}\left(u^2+v^2\right)+g \eta\right)_y-h u \zeta+D\left(v_{x x}+v_{y y}\right) \\
&\eta_t=-(h u)_x-(h v)_y
\end{aligned}
$$
in which the potential vorticity is defined as
%\todo{Alex: There is now some problem with names, we have two letters for 3 quantities?}
$$
\zeta=\frac{f+\zeta}{h} = \frac{f+v_x-u_y}{h}
$$
In this way the gradient term contains the kinetic and potential energy, and the advection term contains the potential vorticity. In absence of dissipation these two quantities are conserved. This means that a discretization of the equations using this form leads to better conservation properties of both energy and potential vorticity.

    \noindent\textbf{Inner grid points}
    %%%%%%%%%%%%%%%%%%%%%%%%%%%%%%%%%%%%%%%%%%%%%%%%%%
    
Denote $i$ as the eastward grid coordinate and $j$ as the northward grid coordinate. A grid box $i, j$ has 4 grid points:
the upper left corner is an $(h, \eta)$ point, upper right is a $u$ point, lower left is a $v$ point, and lower right is a $(f, \xi, \zeta)$ point. Each of them has index $i, j$.
This means that the pressure gradient can be discretized as:

$$
\begin{aligned}
&{\left[-g \eta_x\right]_{i, j}=-g\left(\eta_{i, j}-\eta_{i-1, j}\right) / \Delta x} \\
&{\left[-g \eta_y\right]_{i, j}=-g\left(\eta_{i, j}-\eta_{i, j-1}\right) / \Delta y}
\end{aligned}
$$

Note that the $i, j$ coordinate on the left-hand side is for $u$ grid and $v$ grid, respectively, and on the right-hand side for the $h, \eta$ grid. We will take $\Delta x=\Delta y$.
The energy term is first interpolated to the $(h, \eta)$ grid points as:
$$
\left[\left(u^2+v^2\right)\right]_{i, j}=\frac{1}{2}( u_{i+1, j} u_{i+1, j} +u_{i, j} u_{i, j}+v_{i, j+1} v_{i, j+1}+v_{i, j} v_{i, j})
$$
and then the zonal derivative is taken similar to the pressure gradient term.

The vorticity term is a bit more complicated. We first determine the total vorticity divided by the layer height at the vorticity grid points, as:
$$
\zeta_{i, j}=\left[\frac{f+\xi}{h}\right]_{i j}=\frac{\left(v_{i, j}-v_{i-1, j}-\left(u_{i, j}-u_{i, j-1}\right)\right) / d x+f_j}{H+\left(e_{i, j}+e_{i-1, j}+e_{i, j-1}+e_{i-1, j-1}\right) / 4}
$$
We then determine the mass fluxes in zonal and meridional direction at the $\mathrm{u}$ and $\mathrm{v}$ grid points, respectively:
$\text{vflux}_{i, j}=v_{i, j}\left(H+\left(e_{i, j}+e_{i, j-1}\right) / 2\right)$
The vorticity terms then become, for the zonal and meridional velocities, respectively (note the interpolation of the $v$ flux to the $u$ grid points, and the interpolation of the $u$ flux to the $v$ grid points):
$$
\begin{aligned}
&\left (vflux_{i, j}+ vflux_{i-1, j}+vflux_{i, j+1}+ vflux_{i-1, j+1}\right)\left(\zeta_{i, j}+\zeta_{i, j+1}\right) / 8 \\
&\left(uflux_{i, j}+ uflux_{i, j-1}+ uflux_{i+1, j}+uflux_{i+1, j-1}\right)\left(\zeta_{i, j}+\zeta_{i+1, j}\right) / 8
\end{aligned}
$$
The dissipation terms are not complicated and follow standard central differences:
$$
\begin{array}{r}
D *\left(u_{i+1, j}+u_{i-1, j}+u_{i, j+1}+u_{i, j-1}-4 * u_{i, j}\right) / d x^2-r * u_{i, j} \\
D *\left(v_{i+1, j}+v_{i-1, j}+v_{i, j+1}+v_{i, j-1}-4 * v_{i, j}\right) / d x^2-r * v_{i, j}
\end{array}
$$
Finally, the continuity equation advection terms are already in flux form and need only to be interpolated to the $(h, \eta)$ grid points:
$$
\begin{array}{r}
\left(uflux_{i+1, j}-\operatorname{uflux}_{i, j}\right) / d x \\
\left(vflux_{i, j+1}-vflux_{i, j}\right) / d x
\end{array}
$$
\noindent\textbf{Boundary points}
    %%%%%%%%%%%%%%%%%%%%%%%%%%%%%%%%%%%%%%%%%%%%%%%%%%
At the northern and southern boundaries we assume no meridional flow, so $v=0$. Since this is a hard condition, we have to fulfill this for mass conservation, the north and south boundaries are positions along the $\mathrm{v}$ grid. This means that the boundaries for $u$ and $\eta$ need an interpolation. For $u$ we use a free-slip boundary condition, so $\partial h u / \partial y=0$ at the norther and southern boundary:
$$
\begin{array}{r}
uflux_{i, 1}=-uflux_{i, 2} \\
uflux_{i, n}=-uflux_{i, n-1}
\end{array}
$$
The model is periodic in the zonal direction, which means zonal derivatives at the boundary have to 'wrap around'. For example:
$$
\left[-g \eta_x\right]_{1, j}=-g\left(\eta_{1, j}-\eta_{m, j}\right) / \Delta x
$$

For the time discretization the model uses as Euler step for the first time step (including the first step after the data assimilation), and Leapfrog for all other time seps.
%Pseudocode\\

\section{Numerical implementation of the SPDE}
%%%%%%%%%%%%%%%%%%%%%%%%%%%%%%%%%%%%%%%%%%%%%%%%%%
%%%%%%%%%%%%%%%%%%%%%%%%%%%%%%%%%%%%%%%%%%%%%%%%%%
%%%%%%%%%%%%%%%%%%%%%%%%%%%%%%%%%%%%%%%%%%%%%%%%%

    The stochastic shallow water system with SALT noise is, formally written,
    \begin{equation}
    \begin{aligned}
        \vb{u}_t &= -(\vb{v}\vdot\nabla)\vb{u} - u\grad{\tilde{u}} - v\grad{\tilde{v}} + f({\vb{u}}\cross\vu{z}) - g\grad{h}+D\laplacian{\vb{u}}\\
        {h}_t &= -\div{(h\vb{v})}
    \end{aligned}
    \end{equation}
    where $\vb{v}=\vb{u}+\vb{\tilde{u}}$ denotes the stochastically perturbed velocity with $\vb{\tilde{u}}=[\tilde{u}, \tilde{v}]$ being the stochastic perturbation.
    The stochastic perturbations are modelled by Brownian noise as $\vb{\tilde{u}} = \sum_i \boldsymbol{\xi}_i \dd{W}_t^i$ so that the correct time scaling is 
    \begin{equation}
    \begin{aligned}
        \dd{\vb{u}} &= (-(\vb{u}\vdot\nabla)\vb{u}  + f({\vb{u}}\cross\vu{z}) - g\grad{h}+D\laplacian{\vb{u}})\dd{t} -((\vb{\tilde{u}}\vdot\nabla)\vb{u} + u\grad{\tilde{u}} + v\grad{\tilde{v}})\sqrt{\dd{t}}\\
        \dd{h} &= -\div{(h\vb{u})}\dd{t} - \div{(h\vb{\tilde{u}})}\sqrt{\dd{t}}.
    \end{aligned}
    \end{equation}
    We get
    \begin{equation}
    \begin{aligned}
        \dd{\vb{u}} &= (-(\vb{u}\vdot\nabla)\vb{u}  + f({\vb{u}}\cross\vu{z}) - g\grad{h}+D\laplacian{\vb{u}})\dd{t} -\sum_i((\boldsymbol{\xi}_i\vdot\nabla)\vb{u} + u\grad{\xi_i^u} + v\grad{\xi_i^{v}}){\dd{W_t^i}}\\
        \dd{h} &= -\div{(h\vb{u})}\dd{t} - \sum_i\div{(h\boldsymbol{\xi}_i)}{\dd{W_t^i}}.
    \end{aligned}
    \end{equation}
    We discretise the Brownian motions as $w_i\sqrt{\dd{t}}$, where $w_i\sim \mathcal{N}(0,1)$ so that
    \begin{equation}
    \begin{aligned}
        \dd{\vb{u}} &= (-(\vb{u}\vdot\nabla)\vb{u}  + f({\vb{u}}\cross\vu{z}) - g\grad{h}+D\laplacian{\vb{u}})\dd{t}\\
        &\phantom{=}-\left( \sum_i((\boldsymbol{\xi}_iw_i\vdot\nabla)\vb{u} + u\grad{\xi_i^uw_i} + v\grad{\xi_i^{v}w_i})\right)\sqrt{\dd{t}}\\
        \dd{h} &= -\div{(h\vb{u})}\dd{t} - \div{\left(\sum_i{h\boldsymbol{\xi}_iw_i}\right)}\sqrt{\dd{t}}.
    \end{aligned}
    \end{equation}
    We discretise the stochastic terms analogously to the deterministic case. In the velocity equation, the added terms wrt. the deterministic system are
    \begin{equation}
        -(\vb{\tilde{u}}\vdot\nabla)\vb{u} - u\grad{\tilde{u}} - v\grad{\tilde{v}},
    \end{equation}
    or, in component form,
    \begin{align}
        &-\tilde{u}u_x-\tilde{v}u_y-u\tilde{u}_x - v\tilde{v}_x \\
        &-\tilde{u}v_x-\tilde{v}v_y-u\tilde{u}_y - v\tilde{v}_y.
    \end{align}
    This can be rewritten as
    \begin{align}
        &-(u\tilde{u}+v\tilde{v})_x + h\tilde{v}\frac{v_x-u_y}{h} \\
        &-(u\tilde{u}+v\tilde{v})_y - h\tilde{u}\frac{v_x-u_y}{h},
    \end{align}
    or, more compactly,
    \begin{equation}
        -\grad{(\vb{u}\vdot\vb{\tilde{u}})} + h\tilde{\zeta}(\vb{\tilde{u}}\cross\vu{z}),
    \end{equation}
    where $$\tilde{\zeta}:=\frac{v_x-u_y}{h}$$ denotes the stochastic vorticity contribution due to the SALT noise and the term $h\vb{\tilde{u}}$ is the stochastic flux. The spatial discretisation of the stochastic equation from here on proceeds along the same lines as for the deterministic equation, see Section \ref{appendix:discretisation}

\section{The Runge-Kutta Approximation Scheme} \label{RK}
%%%%%%%%%%%%%%%%%%%%%%%%%%%%%%%%%%%%%%%%%%%%%%%%%%
%%%%%%%%%%%%%%%%%%%%%%%%%%%%%%%%%%%%%%%%%%%%%%%%%%
%%%%%%%%%%%%%%%%%%%%%%%%%%%%%%%%%%%%%%%%%%%%%%%%%%

For a dynamical system
\begin{equation}
\frac{dX}{dt}= F(X,t)
\end{equation}
the Runge-Kutta scheme of order 4 is given by
\begin{equation}
X^{n+1} = X^n + \frac{\Delta t}{6}(c_1 + 2c_2 + 2c_3 + c_4)
\end{equation}
where $c_i, i \in \{1,2,3,4\}$ are the intermediate steps
\begin{subequations}\label{rk}
\begin{alignat}{5}
& c_1 = F(X^n,t) \\
& c_2 = F\left(X^n + \frac{\Delta t}{2}c_1, t+\frac{\Delta t}{2}\right) \\
& c_3 = F\left(X^n + \frac{\Delta t}{2}c_2, t+\frac{\Delta t}{2}\right) \\
& c_4 = F\left(X^n + \Delta tc_3, t+\Delta t\right)
\end{alignat}
\end{subequations}
with $X^n = X(t_n) = X_{t_n}$ and $t_n = n\Delta t$ for $n \in \{0,1,2, \ldots, N-1\}, t_{N}=T, \Delta t= t_{n+1}-t_n$. Each of these intermediate steps is solved using a standard Euler approximation scheme (see e.g. \cite{KloedenPlaten}). 
The SPDE studied in this paper is of the form
\begin{equation}\label{stratonovich}
dF(X,t) = A(X,t)dt + \displaystyle\sum_iB^i(X, t) \circ dW_t^i.
\end{equation}
We show below that from a numerical perspective we can start by approximating the SPDE in the It\^{o} form, i.e.,
\begin{equation}\label{ito1}
dF(X,t) = A(X,t)dt + \displaystyle\sum_iB^i(X, t)dW_t^i
\end{equation} 
as due to the specific structure of the Runge-Kutta of order 4 scheme, the It\^{o} correction associated with the Stratonovich integral from equation \eqref{stratonovich} appears automatically in the second iteration. The discretized version of equation \eqref{ito1} is given by  
\begin{equation}
F(X^n, t) = A(X^n)\Delta t + \displaystyle\sum_i B^i(X^n) \Delta W^i
\end{equation}
where $\Delta W^i = W_{t_{n+1}}^i-W_{t_n}^i$ and $\Delta t = t_{n+1}-t_n$. Then \footnote{We drop the time dependence for the moment. Also, since in order to see how we retrieve the It\^{o} correction it is enough to compute the first three intermediate steps, we skip also the explicit calculation of $c_4$.}
\begin{subequations}
\begin{alignat}{3}
& c_1 = F(X^n) = A(X^n)\Delta t + \displaystyle\sum_i B^i(X^n)\Delta W^i  = F(\ell_1)\\
& c_2  =F\left(X^n + \frac{c_1}{2}\right)
 = F \left(X^n + \frac{1}{2}\left( A(X^n)\Delta t + \displaystyle\sum_i B^i(X^n)\Delta W^i\right)\right) = F(\ell_2) \\
& c_3 = F \left(X^n + \frac{c_2}{2}\right) = F(\ell_3)
\end{alignat} 
\end{subequations}
where 
\begin{equation*}
\begin{aligned} 
& \ell_1:=X^n \\
& \ell_2 := X^n + \frac{F(\ell_1)}{2} = X^n + \frac{1}{2}\left(A(X^n)\Delta t + \displaystyle\sum_i B^i(X^n)\Delta W^i\right)
\end{aligned}
\end{equation*}
\begin{equation*}
\begin{aligned}
\ell_3 &:= X^n + \frac{F(\ell_2)}{2} = X^n + \frac{1}{2}F\left(X^n + \frac{1}{2} \left(A(X^n)\Delta t + \displaystyle\sum_i B^i(X^n)\Delta W^i \right)\right)  \\
& = X^n + \frac{1}{2}F(X^n) + \frac{1}{4}F\left(A(X^n)\Delta t + \displaystyle\sum_i B^i(X^n)\Delta W^i\right) \\
& = X^n + \frac{1}{2}A(X^n)\Delta t + \frac{1}{2}\displaystyle\sum_i B^i(X^n)\Delta W^i \\ 
& + \frac{1}{4}A^2(X^n)(\Delta t)^2 + \frac{1}{4}\displaystyle\sum_iAB^i(X^n)\Delta W^i\Delta t + \frac{1}{4}\displaystyle\sum_i B^i(A(X^n))\Delta t \Delta W^i \\ &+\frac{1}{4}\displaystyle\sum_i(B^i)^2(X^n)(\Delta W^i)^2 \\
& = X^n + \frac{1}{2} \left( A(X^n)\Delta t + \displaystyle\sum_iB^i(X^n)\Delta W^i + \frac{1}{2} \displaystyle\sum_i(B^i)^2(X^n)(\Delta W^i)^2\right) \\
& + \frac{1}{2} \left( A^2(X^n)(\Delta t)^2 + \frac{1}{2}\displaystyle\sum_iA(B^i(X^n))\Delta W^i \Delta t + \frac{1}{2}\displaystyle\sum_iB^i(A(X^n))\Delta t \Delta W^i\right).
\end{aligned} 
\end{equation*}
Then
\begin{equation*}
\ell_3 = X^n + \frac{1}{2} \left( A(X^n)\Delta t + \displaystyle\sum_iB^i(X^n)\Delta W^i + \frac{1}{2} \displaystyle\sum_i(B^i)^2(X^n)(\Delta W^i)^2\right) + \hbox{higher order terms}
\end{equation*}
and therefore we have recovered the It\^{o} correction $\frac{1}{2} \displaystyle\sum_i(B^i)^2(X^n)(\Delta W^i)^2$. This is known (see \cite{KloedenPlaten}) as \textit{Heun's method}  or the \textit{improved Euler method}. Formally it is based on introducing an auxiliary variable (we neglect the dependence on $i$ here as it is not essential)
\begin{equation*}
X^{\star} := X^n + A(X^n)\Delta t + B(X^n)\Delta W.
\end{equation*}
Then in the Euler scheme used for calculating the intermediate steps $c_i, i\in \{1,2,3,4\}$ we have
\begin{equation*}
\begin{aligned} 
X^{n+1} &= X^n + \frac{\Delta t}{2}(A(X^n)+A(X^{\star}))+\frac{\Delta W}{2}(B(X^n)+B(X^{\star})) \\
& = X^n + A(X^n)\Delta t + B(X^n)\Delta W + + \frac{1}{2} B^2(X^n)(\Delta W)^2\\
&+\frac{1}{2}A^2(X^n)(\Delta t)^2 + \frac{1}{2} A(B(X^n))\Delta t \Delta W  + \frac{1}{2} (B(A)(X^n))\Delta t \Delta W 
\end{aligned} 
\end{equation*}
which is similar to what we had before. 

\end{appendices}
%%%%%%%%%%%%%%%%%%%%%%%%%%%%%%%%%%%%%%%%%%%%%%%%%%

%\listoftodos

\end{document}